\documentclass[a4paper]{amsart}
\usepackage[latin1]{inputenc}
\usepackage{latexsym}
\input amssym.def
\input amssym.tex
\newtheorem{definition}{Definition}[section]

\newtheorem{theorem}[definition]{Theorem}
\newtheorem{proposition}[definition]{Proposition}
\newtheorem{corollary}[definition]{Corollary}
\newtheorem{remark}[definition]{Remark}
\newtheorem{example}[definition]{Example}
\newtheorem{examples}[definition]{Examples}

\def\emph#1{{\bfseries\itshape{#1}}}
 1   
 1  

\def\R{\mathbb{R}}               

\def\lcf{\lbrack\! \lbrack}
\def\rcf{\rbrack\! \rbrack}

\newcount\ancho
\newcount\anchom
\newcount\anchoa
\newcount\anchob
\newcommand{\ltilde}[2]{\ancho=#1 \anchom=\ancho \divide\anchom by 2
            \anchoa=\ancho \divide\anchoa by 4
        \anchob=\anchom \advance\anchob by \anchoa
      $\kern-5pt \begin{array}[b]{c}
                 \begin{picture}(1,1)(\anchom,0)
         \qbezier(0,2)(\anchoa,5)(\anchom,2)
         \qbezier(\anchom,2)(\anchob,-1)(\ancho,4)
         \qbezier(0,2)(\anchoa,4.5)(\anchom,1.8)
         \qbezier(\anchom,1.8)(\anchob,-1.5)(\ancho,4)
      \end{picture} \\[-4pt]
       \mbox{#2}
       \end{array} \kern-9pt$
       }
\begin{document}

\title{Lagrangian submanifolds and dynamics on Lie affgebroids}
\author{D. Iglesias,\;\; J.C. Marrero,\;\; E.  Padr\'{o}n, \;\; D.
Sosa}

\thanks{\noindent Supported in part by BFM2003-01319. D.~Iglesias wishes to thank the Spanish
Ministry of Education and Culture and Fulbright program for a
MECD/Fulbright postdoctoral grant.}



\begin{abstract}
\noindent We introduce the notion of a symplectic Lie affgebroid
and their Lagrangian submanifolds in order to describe the
Lagrangian (Hamiltonian) dynamics on a Lie affgebroid in terms of
this type of structures. Several examples are discussed.
\end{abstract}

\maketitle

\setcounter{section}{0}

\section{Introduction}
\setcounter{equation}{0} Recently, there has been a lot of
interest in the study of Lie algebroids, which can be thought of
as ``generalized tangent bundles", since they generalize Lie
algebras as well as (regular) integrable distributions. From the
Physics point of view, Lie algebroids can be used to give
geometric descriptions of Hamiltonian and Lagrangian Mechanics. In
\cite{We}, A. Weinstein introduces ``Lagrangian systems" on a Lie
algebroid $E$ by means of the linear Poisson structure on the dual
$E^\ast$ and a Legendre-type map from $E$ to $E^*$, associated to
a given Lagrangian function $L$ on $E$, provided that $L$ is
regular. In that paper, he also asks about the possibility to
develop a geometric formalism on Lie algebroids similar to Klein's
formalism in ordinary Lagrangian Mechanics. An answer for this
question was  given by Martinez in \cite{M1} (see also
\cite{M2,PPo}), using the notion of prolongation of a Lie
algebroid over a mapping \cite{HM}.

In \cite{Tul1,Tul2} an interpretation of the Classical Lagrangian
and Hamiltonian dynamics as Lagrangian submanifolds of convenient
special symplectic manifolds is described. In doing so, one
introduces certain canonical isomorphisms which allow to build the
so-called Tulczyjew's triple of Classical Mechanics (see
\cite{Tul1,Tul2}). In \cite{LMM}, this line of research has been
followed. More precisely, the authors introduce the notion of a
symplectic Lie algebroid and their Lagrangian submanifolds in
order to give an interpretation of Lagrangian and Hamiltonian
Mechanics on Lie algebroids in terms of Lagrangian submanifolds of
symplectic Lie algebroids. In addition, as an application, the
authors recover Lagrange-Poincar\'{e} \cite{CMR} (respectively,
Hamilton-Poincar\'{e} \cite{CMPR}) equations associated with a
$G$-invariant Lagrangian (respectively, Hamiltonian) on a
principal $G$-bundle $p: Q\to M$ as the Euler-Lagrange
(respectively, Hamilton) equations on the Atiyah algebroid $TQ/G$.

On the other hand, in \cite{GGrU,MMeS} (see also
\cite{GGU2,Ma3,SMM}) a possible generalization of the concept of a
Lie algebroid to affine bundles is introduced in order to create a
geometric model which provides a natural framework for a
time-dependent version of Lagrange (Hamilton) equations on Lie
algebroids. These new structures are called Lie affgebroid
structures (in the terminology of \cite{GGrU}).

The main aim of this paper is to introduce the notion of a
Lagrangian submanifold of a symplectic Lie affgebroid and, then,
to use this type of geometric objects to describe the Lagrangian
(Hamiltonian) dynamics on Lie affgebroids.

The paper is organized as follows. In Section \ref{Lie-algebroids}
we recall the notion of a Lie algebroid and several constructions
related with them, in particular the prolongation ${\mathcal
L}^fE$ of a Lie algebroid $\tau :E\to M$ over a smooth map
$f:M'\to M$ and the definition of the Lie algebroid structure of
an action Lie algebroid. In Section \ref{sec1.2}, we describe the
concept of a Lie affgebroid structure on an affine bundle $\tau
_A:A\to M$ modelled on a vector bundle $\tau _V:V\to M$. We remark
that if $\tau _{A^+}:A^+=Aff(A,\R)\to M$ is the dual bundle to $A$
and $\widetilde{A}=(A^+)^*$ is the bidual, then a Lie affgebroid
structure on $A$ is equivalent to a Lie algebroid structure on
$\widetilde{A}$ such that the distinguished section $1_A$ of $\tau
_{A^+}:A^+\to M$ (corresponding to the constant function $1$ on
$A$) is a 1-cocycle for the Lie algebroid cohomology.

The description of the Hamiltonian dynamics on a Lie affgebroid is
presented in Section \ref{sec3.1} following \cite{M}. This
geometric framework allows to write, for a Lie affgebroid $\tau_A:
A\to M$ modelled on a vector bundle $\tau_V:V\to M$,  the Hamilton
equations associated with a Hamiltonian section $h:V^\ast \to A^+$
of the canonical projection $\mu:A^+\to V^*$. In doing so, we
introduce a cosymplectic structure $(\Omega _h, \eta)$ on the
prolongation ${\mathcal  L}^{\tau _{V}^*}\widetilde{A}$ of
$\widetilde{A}$ over $\tau_{V}^*:V^*\to M$. Then, the integral
curves of the Reeb section $R_h$ of $(\Omega _h, \eta)$ are just
the solutions of the Hamilton equations. Alternatively, one may
prove that the solutions of the Hamilton equations are just the
integral curves of the Hamiltonian vector field (on $V^*$) of $h$
with respect to the canonical aff-Poisson structure on the line
affine bundle $\mu:A^+\to V^*$ (see Theorem \ref{teor3.2}).
Aff-Poisson structures were introduced in \cite{GGrU} (see also
\cite{GGU2}) as the affine version of standard Poisson structures.
On the other hand, in Section \ref{sec3.2} it is developed the
corresponding Lagrangian formalism for a Lagrangian function
$L:A\to \R$ (see \cite{M}). In this case, we work on the
prolongation ${\mathcal L}^{\tau _A}\widetilde{A}$ of
$\widetilde{A}$ over $\tau_A$. We define the Poincar\'e-Cartan
2-section $\Omega _L$ (as a section of the vector bundle
$\wedge^2({\mathcal  L}^{\tau _A}\widetilde{A})^*\to A$) and the
vertical endomorphism $S$ (as a section of ${\mathcal L}^{\tau
_A}\widetilde{A}\otimes({\mathcal L}^{\tau _A}\widetilde{A})^*\to
A$). These objects allow us to write the Euler-Lagrange equations
for $L$ in an intrinsic way. In the particular case, when $L$ is
regular the pair $(\Omega _L,\phi _0)$ is a cosymplectic structure
on ${\mathcal  L}^{\tau _{A}}\widetilde{A}$, $\phi _0$ being a
certain 1-cocycle of ${\mathcal  L}^{\tau _A}\widetilde{A}$, and
the integral curves of the Reeb section $R_L$ of $(\Omega _L,\phi
_0)$ are just the solutions of the Euler-Lagrange equations for
$L$. In Section \ref{sec3.3} it is stated the equivalence between
both formalisms using the Legendre transformation $leg _L:A\to
V^\ast$, provided that $L$ is hyperregular (that is, $leg _L:A\to
V^\ast$ is a global diffeomorphism). In fact, we may construct a
Hamiltonian section $h_L:V^*\to A^+$ and we have  a Lie algebroid
morphism ${\mathcal L}leg_L: {\mathcal L}^{\tau_A}\widetilde{A}\to
{\mathcal L}^{\tau_{V^*}}\widetilde{A}$ over $leg _L:A\to V^\ast$
which in the hyperregular case is an isomorphism and, in addition,
connects the Lagrangian and Hamiltonian formalism. Conversely, if
$h:V^*\to A^+$ is a Hamiltonian section, one may introduce a map
${\Bbb F}h:V^*\to A$ such that $\tau_A\circ {\Bbb F}h=\tau_{V}^*$.
Furthermore, if $h$ is hyperregular (that is, ${\Bbb F}h:V^*\to A$
is a global diffeomorphism) then there exists a hyperregular
Lagrangian $L:A\to {\Bbb R}$ such that $h_L=h$ and ${\Bbb
F}h=leg_L^{-1}.$

In Sections \ref{sec3} and \ref{sec5}, we extend the construction
of Tulczyjew's triple to the Lie affgebroid setting. More
precisely, for a Lie affgebroid $\tau_A:A\to M$ we consider the
space $${\mathcal  J}^AA=\{(a,v)\in A\times
TA/\rho_A(a)=(T\tau_A)(v)\},$$ where $\rho _A:A\to TM$ is the
anchor map of $A$, and we prove that ${\mathcal  J}^A A$ admits
two Lie aff\-ge\-broids structures. These structures are
isomorphic under the so-called canonical involution $\sigma
_{A}:{\mathcal J}^AA\to {\mathcal  J}^AA$ associated with $A$. In
addition, if $h:V^\ast \to A^+$ is a Hamiltonian section then we
define an affine isomorphism $\flat_{\Omega_h}$ (over the identity
of $V^*$) between the affine bundles $\rho_A^*(TV^*)\to V^*$ and
$({\mathcal L}^{\tau_V^*}V)^*\to V^*$. Here, $\rho_A^*(TV^*)$ is
the pull-back of the vector bundle $T\tau_{V}^*:TV^*\to TM$ over
$\rho_A$ and ${\mathcal  L}^{\tau^*_V}V$ is the prolongation of
$V$ over $\tau_V^*$. The map $\flat_{\Omega_h}$, along with a
canonical vector bundle isomorphism
$A_A:\rho_A^*(TV^*)\rightarrow({\mathcal L}^{\tau_A}V)^*$ (related
with $\sigma _A$) gives us the Tulczyjew's triple associated with
$A$ and $h$.

In Section \ref{sec6}, we introduce the notion of a symplectic Lie
affgebroid as a Lie affgebroid modelled on a symplectic Lie
algebroid (that is, a Lie algebroid $\tau _V :V\to M$ which admits
a non-degenerate 2-cocycle). Moreover, we prove that if
$\tau_A:A\to M$ is a symplectic Lie affgebroid then its
prolongation ${\mathcal  J}^A A$ is also a symplectic Lie
affgebroid. The notion of a Lagrangian Lie subaffgebroid of a
symplectic Lie affgebroid is introduced, in a natural way, in
Section \ref{sec7}. Examples and properties of this type of
objects are discussed in this section.

In Section \ref{sec8}, using the results of Section \ref{sec7}, we
introduce the definition of a Lagrangian submanifold of a
symplectic Lie affgebroid. Then, if $h:V^\ast \to A^+$ is a
Hamiltonian section we deduce that $S_h=R_h (V^\ast)$ is a
Lagrangian submanifold of the symplectic Lie affgebroid
$\rho_A^*(TV^*)$ and, in addition, there exists a bijection
between {\it admissible curves } in $S_h$ and solutions of the
Hamilton equations for $h$. Similarly, given a Lagrangian function
$L:A\to \R$, we prove that $S_L=(A_A^{-1}\circ d^{{\mathcal
L}^{\tau_A}V}L)(A)$ is a Lagrangian submanifold of
$\rho_A^*(TV^*)$ and that there exists a bijection between {\it
admissible curves } in $S_L$ and solutions of the Euler-Lagrange
equations for $L$. When $L$ is hyperregular and $h$ is the
corresponding Hamiltonian section, we deduce that $S_L=S_h$.

Finally, we describe some applications in Section \ref{sec9}. In
fact, in Section \ref{sec9.0}, we prove that if the Lie affgebroid
$A$ is a Lie algebroid then we recover the results obtained in
\cite{LMM} about the relation between Lagrangian submanifolds and
dynamics on Lie algebroids. In addition, in Section \ref{sec9.1},
we apply the results of the paper to the particular case when the
Lie affgebroid $A$ is the $1$-jet bundle $\tau_{1,0}:J^1\tau\to M$
of local sections of a fibration $\tau:M\to \R$. As a consequence,
we deduce that the classical Euler-Lagrange (Hamilton) equations
of time-dependent Mechanics are just the local equations defining
Lagrangian submanifolds of a symplectic Lie affgebroid. On the
other hand, in Section \ref{sec9.2}, we consider a principal
$G$-bundle $p:Q\to M$ such that the base space $M$ is fibred on
$\R$, that is, there exists a fibration $\nu :M\to \R$. Then, if
$\tau=\nu\circ p,$ we have that the quotient affine bundle
$\tau_{1,0}|G:J^1\tau/G\to M$ admits a Lie affgebroid structure in
such way that the bidual Lie algebroid is just the Atiyah
algebroid $\pi_Q|G:TQ/G\to M$ associated with the principal
$G$-bundle $p:Q\to M$ (see \cite{MMeS}). For this reason, the
affine bundle $\tau_{1,0}|G:J^1\tau/G\to M$ is called an Atiyah
affgebroid. We obtain that, in this case, the solutions of the
Euler-Lagrange (Hamilton) equations for a Lagrangian (resp., a
Hamiltonian section) are the solutions of the classical
nonautonomous Lagrange-Poincar\'{e} (resp. Hamilton-Poincar\'{e})
equations for the corresponding $G$-invariant Lagrangian (resp.
$G$-invariant Hamiltonian section). Moreover, all these equations
are reinterpreted as those defining the corresponding Lagrangian
submanifolds of an Atiyah symplectic Lie affgebroid.

Manifolds are real, paracompact and $C^{\infty}$. Maps are
$C^{\infty}$. Sum over crossed repeated indices is understood.

\section{Lie algebroids and Lie affgebroids}\label{preparation}
\subsection{Lie algebroids}\label{Lie-algebroids}
Let $E$ be a vector bundle of rank $n$ over the manifold $M$ of
dimension $m$ and $\tau:E\rightarrow M$ be the vector bundle
projection. Denote by $\Gamma(\tau)$ the $C^{\infty}(M)$-module of
sections of $\tau:E\rightarrow M$. A {\it Lie algebroid structure}
$(\lcf\cdot,\cdot\rcf,\rho)$ on $E$ is a Lie bracket
$\lcf\cdot,\cdot\rcf$ on the space $\Gamma(\tau)$ and a bundle map
$\rho:E\rightarrow TM$, called {\it the anchor map}, such that if
we also denote by $\rho:\Gamma(\tau)\rightarrow\frak{X}(M)$ the
homomorphism of $C^{\infty}(M)$-modules induced by the anchor map
then $\lcf X,fY\rcf=f\lcf X,Y\rcf+\rho(X)(f)Y,$ for
$X,Y\in\Gamma(\tau)$ and $f\in C^{\infty}(M)$. The triple
$(E,\lcf\cdot,\cdot\rcf,\rho)$ is called a {\it Lie algebroid over
$M$} (see \cite{Ma}). In such a case, the anchor map
$\rho:\Gamma(\tau)\rightarrow\frak{X}(M)$ is a homomorphism
between the Lie algebras $(\Gamma(\tau),\lcf\cdot,\cdot\rcf)$ and
$(\frak{X}(M),[\cdot,\cdot])$.



If $(E,\lcf\cdot,\cdot\rcf,\rho)$ is a Lie algebroid, one may
define a cohomology operator, which is called {\it the
differential of $E$},
$d^E:\Gamma(\wedge^k\tau^*)\longrightarrow\Gamma(\wedge^{k+1}\tau^*)$,
as follows
$$
\begin{array}{lcl}
(d^E\mu)(X_0,\dots,X_k)&=&\displaystyle\sum_{i=0}^k(-1)^i\rho(X_i)(\mu(X_0,\dots,\widehat{X_i},\dots,X_k))\\
 &+&\displaystyle\sum_{i<j}(-1)^{i+j}\mu(\lcf
 X_i,X_j\rcf,X_0,\dots,\widehat{X_i},\dots,\widehat{X_j},\dots,X_k),
\end{array}
$$
for $\mu\in\Gamma(\wedge^k\tau^*)$ and
$X_0,\dots,X_k\in\Gamma(\tau)$. Moreover, if $X\in\Gamma(\tau)$
one may introduce, in a natural way, {\it the Lie derivate with
respect to $X$}, as the operator ${\mathcal
L}_X^E:\Gamma(\wedge^k\tau^*)\longrightarrow\Gamma(\wedge^{k}\tau^*)$
given by ${\mathcal  L}_X^E=i_X\circ d^E+d^E\circ i_X.$

If $E$ is the standard Lie algebroid $TM$ then the differential
$d^E=d^{TM}$ is the usual exterior differential associated with
$M$, which we will denote by $d_0$.

 Now, suppose that
$(E,\lcf\cdot,\cdot\rcf,\rho)$ and
$(E',\lcf\cdot,\cdot\rcf',\rho')$ are Lie algebroids over $M$ and
$M'$, respectively, and that $F:E\to E'$ is a vector bundle
morphism over the map $f:M\to M'.$ Then $(F,f)$ is said to be a
{\it Lie algebroid morphism} if
$$d^{E}((F,f)^*\phi')=(F,f)^*(d^{E'}\phi'),\mbox{ for }\phi'\in\Gamma(\wedge^k(\tau')^*) \mbox{ and for all } k.$$
Note that $(F,f)^*\phi'$ is the section of the vector bundle
$\wedge^k E^*\rightarrow M$ defined by
$$((F,f)^*\phi')_{x}(a_1,\dots,a_k)=\phi'_{f(x)}(F(a_1),\dots,F(a_k)),$$
for $x\in M$ and $a_1,\dots,a_k\in E_{x}$. If $(F,f)$ is a Lie
algebroid morphism, $f$ is an injective immersion and
$F_{|E_x}:E_x\rightarrow E'_{f(x)}$ is injective, for all $x\in
M$, then $(E,\lcf\cdot,\cdot\rcf,\rho)$ is said to be a {\it Lie
subalgebroid} of $(E',\lcf\cdot,\cdot\rcf',\rho')$.

\subsubsection{The prolongation of a Lie algebroid over a
smooth map}\label{sec1.1.1}

In this section, we will recall the definition of the Lie
algebroid structure on the prolongation of a Lie algebroid over a
smooth map (see \cite{HM,LMM}).

Let $(E,\lcf\cdot,\cdot\rcf,\rho)$ be a Lie algebroid of rank $n$
over a manifold $M$ of dimension $m$ with vector bundle projection
$\tau:E\rightarrow M$ and $f:M'\rightarrow M$ be a smooth map.

We consider the subset ${\mathcal  L}^{f}E$ of $E\times TM'$ and
the map $\tau^{f}:{\mathcal  L}^{f}E\rightarrow M'$ defined by
$${\mathcal  L}^{f}E=\{(b,v')\in E\times
TM'/\rho(b)=(Tf)(v')\},\;\;\;\;\tau^{f}(b,v')=\pi_{M'}(v'),$$
where $Tf:TM'\rightarrow TM$ is the tangent map to $f$ and
$\pi_{M'}:TM'\rightarrow M'$ is  the canonical projection.

Now, assume that $\rho(E_{f(x')})+(T_{x'}f)(T_{x'}M')=T_{f(x')}M,$
for all $x'\in M'.$ Then, $\tau^{f}:{\mathcal  L}^{f}E\rightarrow
M'$ is a vector bundle over $M'$ of rank $n+\dim M'-m$ which
admits a Lie algebroid structure
$(\lcf\cdot,\cdot\rcf^{f},\rho^{f})$ characterized by
$$\lcf(X\circ f,U'),(Y\circ f,V')\rcf^{f}=(\lcf X,Y\rcf\circ f
,[U',V']),\;\;\rho^{f}(X\circ f,U')=U',$$
for all  $X,Y\in \Gamma(\tau)$  and $U', V'$ $f$-projectable
vector fields to $\rho(X)$ and $\rho(Y)$, respectively.
 $({\mathcal
L}^{f}E,\lcf\cdot,\cdot\rcf^{f},\rho^{f})$ is called {\it the
prolongation of the Lie algebroid $E$ over the map $f$} (for more
details, see \cite{HM,LMM}).

Next, we consider a particular case of the above construction. Let
$E$ be a Lie algebroid over a manifold $M$ with vector bundle
projection $\tau:E\rightarrow M$ and ${\mathcal  L}^{\tau^*}E$ be
the prolongation of $E$ over the projection $\tau^*:E^*\rightarrow
M$. ${\mathcal  L}^{\tau^*}E$ is a Lie algebroid over $E^*$ and we
can define a canonical section $\lambda_E$ of the vector bundle
$({\mathcal  L}^{\tau^*}E)^*\rightarrow E^*$ as follows. If
$a^*\in E^*$ and $(b,v)\in({\mathcal  L}^{\tau^*}E)_{a^*}$ then
\begin{equation}\label{lambdaE}
\lambda_E(a^*)(b,v)=a^*(b).
\end{equation}
$\lambda_E$ is called the {\it Liouville section} associated with
the Lie algebroid $E.$

Now, one may consider the nondegenerate section
$\Omega_E=-d^{{\mathcal  L}^{\tau^*}E}\lambda_E$ of
$\wedge^2({\mathcal  L}^{\tau^*}E)^*\rightarrow E^*$. It is clear
that $d^{{\mathcal  L}^{\tau^*}E}\Omega_E=0$. In other words,
$\Omega_E$ is a symplectic section. $\Omega_E$ is called {\it the
canonical symplectic section} associated with the Lie algebroid
$E$. Using the symplectic section $\Omega_E$ one may introduce a
linear Poisson structure $\Pi_{E^*}$ on $E^*$, with linear Poisson
bracket $\{\cdot, \cdot\}_{E^*}$ given by
\[\{F,G\}_{E^*}=-\Omega_E(X_F,Y_G),\;\;\;\; \mbox{ for } F,G\in
C^\infty(E^*),
\]
where $X_F$ and $X_G$ are the Hamiltonian sections associated with
$F$ and $G$, that is, $i_{X_F}\Omega_E=d^{{\mathcal
L}^{\tau^*}E}F$ and $i_{X_G}\Omega_E=d^{{\mathcal L}^{\tau^*}E}G.$

  Suppose that $(x^i)$
are local coordinates on an open subset $U$ of $M$ and that
$\{e_{\alpha}\}$ is a local basis of sections of the vector bundle
$\tau^{-1}(U)\rightarrow U$ such that
$$\rho(e_{\alpha})=\rho_{\alpha}^i\frac{\partial}{\partial
x^i}\;,\;\;\;\lcf
e_{\alpha},e_{\beta}\rcf=C_{\alpha\beta}^{\gamma}e_{\gamma}\;.$$
Then, $\{\tilde{e}_{\alpha},\bar{e}_{\alpha}\}$ is a local basis
of sections of the vector bundle
$(\tau^{\tau^*})^{-1}((\tau^*)^{-1}(U))\rightarrow
(\tau^*)^{-1}(U)$, where $\tau^{\tau^*}:{\mathcal
L}^{\tau^*}E\rightarrow E^*$ is the vector bundle projection and
$$\tilde{e}_{\alpha}(a^*)=(e_{\alpha}(\tau^*(a^*)),\rho_{\alpha}^i\displaystyle\frac{\partial}{\partial
x^i}_{|a^*}),\;\;\;\bar{e}_{\alpha}(a^*)=(0,\displaystyle\frac{\partial}{\partial
y_{\alpha}}_{|a^*}).
$$
Here, $(x^i,y_{\alpha})$ are the local coordinates on $E^*$
induced by the local coordinates $(x^i)$ and the dual basis
$\{e^{\alpha}\}$ of $\{e_{\alpha}\}$. Moreover, we have that
\begin{equation}\label{coretilde}
\begin{array}{c}
\lcf\tilde{e}_{\alpha},\tilde{e}_{\beta}\rcf^{\tau^*}\kern-8pt=C_{\alpha\beta}^{\gamma}\tilde{e}_{\gamma},\;\;
\lcf\tilde{e}_{\alpha},\bar{e}_{\beta}\rcf^{\tau^*}\kern-8pt=\lcf\bar{e}_{\alpha},\bar{e}_{\beta}\rcf^{\tau^*}\kern-8pt=0,\;\;
\rho^{\tau^*}(\tilde{e}_{\alpha})\kern-3pt=\rho_{\alpha}^i\displaystyle\frac{\partial}{\partial
x^i},\;\;\rho^{\tau^*}(\bar{e}_{\alpha})\kern-3pt=\displaystyle\frac{\partial}{\partial
y_{\alpha}},
\end{array}\end{equation}
and
\begin{equation}\label{formas}
\lambda_E(x^i,y_{\alpha})=y_{\alpha}\tilde{e}^{\alpha},\;\;\;\Omega_E(x^i,y_{\alpha})=\tilde{e}^{\alpha}\wedge\bar{e}^{\alpha}+
\displaystyle\frac{1}{2}C_{\alpha\beta}^{\gamma}y_{\gamma}\tilde{e}^{\alpha}\wedge\tilde{e}^{\beta},
\end{equation}
\begin{equation}\label{eq2.3'}
\Pi_{E^*}=\displaystyle\frac{1}{2}C_{\alpha\beta}^\gamma
y_\gamma\frac{\partial}{\partial
y_\alpha}\wedge\frac{\partial}{\partial
y_\beta}+\rho_\alpha^i\frac{\partial}{\partial
y_\alpha}\wedge\frac{\partial}{\partial x^i},
\end{equation}
(for more details, see \cite{LMM,M2}).

\subsubsection{Action Lie algebroids}\label{sec1.1.2}

In this section, we will recall the definition of the Lie
algebroid structure of an action Lie algebroid (see
\cite{HM,LMM}).

Let $(E,\lcf\cdot,\cdot\rcf,\rho)$ be a Lie algebroid over a
manifold $M$ and $f:M'\rightarrow M$ be a smooth map. Then, the
pull-back of $E$ over $f$, $f^*E=\{(x',a)\in M'\times
E/f(x')=\tau(a)\},$ is a vector bundle over $M'$ whose vector
bundle projection is the restriction to $f^*E$ of the first
canonical projection $pr_1:M'\times E\rightarrow M'$. However,
$f^*E$ is not, in general, a Lie algebroid.

Now, suppose that $\Psi:\Gamma(\tau)\rightarrow\frak X(M')$ is an
action of $E$ on $f$, that is, $\Psi$ is a $\R$-linear map which
satisfies the following conditions
$$\Psi(hX)=(h\circ f)\Psi X,\;\;\Psi\lcf X,Y\rcf=[\Psi X,\Psi
Y],\;\;\Psi X(h\circ f)=\rho(X)(h)\circ f,$$ for
$X,Y\in\Gamma(\tau)$ and $h\in C^{\infty}(M)$. Then, one may
introduce a Lie algebroid structure
$(\lcf\cdot,\cdot\rcf_{\Psi},\rho_{\Psi})$ on the vector bundle
$f^*E\rightarrow M'$ which is characterized by the following
conditions
\begin{equation}\label{2.3'}\lcf X\circ f,Y\circ
f\rcf_{\Psi}=\lcf X,Y\rcf\circ f,\;\;\rho_{\Psi}(X\circ
f)=\Psi(X),\makebox[1cm]{for} X,Y\in\Gamma(\tau).
\end{equation}

The resultant Lie algebroid is denoted by $E\ltimes f$ and we call
it {\it an action Lie algebroid}.

Next, we will apply the above construction to a particular case.
First of all, we recall that there is a one-to-one correspondence
between linear functions on a vector bundle $E$ and sections of
$E^*$. In this paper, we don't distinguish between a section of
$E^*$ and its associated linear function on $E$. Let
$(E,\lcf\cdot,\cdot\rcf,\rho)$ be a Lie algebroid with vector
bundle projection $\tau:E\rightarrow M$ and $X$ be a section of
$\tau:E\rightarrow M$. Then, we define {\it the vertical lift} of
$X$ as the vector field on $E$ given by $X^v(a)=X(\tau(a))_a^v,$
for $a\in E,$ where $\;_a^v:E_{\tau(a)}\rightarrow
T_a(E_{\tau(a)})$ is the canonical vector space isomorphism.
In addition, there exists a unique vector field $X^c$ on $E$, {\it
the complete lift} of $X$, which satisfies the following
conditions: $i)$ $X^c$ is $\tau$-projectable on $\rho(X)$ and
$ii)$ $X^c({\alpha})={{\mathcal  L}_X^E\alpha}$, for all
$\alpha:E\to \R$ section of $E^*$ (for more details, see
\cite{GU1,GU2,LMM}).

On the other hand, it is well-known (see, for instance, \cite{Go})
that the tangent bundle to $E$, $TE$, is a vector bundle over $TM$
with vector bundle projection the tangent map to $\tau$,
$T\tau:TE\rightarrow TM$. Moreover, the tangent map to $X$,
$TX:TM\rightarrow TE$, is a section of the vector bundle
$T\tau:TE\rightarrow TM$. We may also consider  the section
$\hat{X}:TM\rightarrow TE$ of $T\tau:TE\rightarrow TM$ given by
\begin{equation}\label{Xtil}
\hat{X}(u)=(T_x0)(u)+X(x)_{0(x)}^v, \end{equation}
 for $u\in
T_xM$, where $0:M\rightarrow E$ is the zero section of $E$ and
$\;_{0(x)}^v:E_x\rightarrow T_{0(x)}(E_x)$ is the canonical
isomorphism between $E_x$ and $T_{0(x)}(E_x)$.

If $\{e_{\alpha}\}$ is a local basis of $\Gamma(\tau)$ then $\{
Te_{\alpha},\widehat{e}_{\alpha}\}$ is a local basis of
$\Gamma(T\tau)$.

The vector bundle $T\tau:TE\rightarrow TM$ admits a Lie algebroid
structure with anchor map $\rho^T$ given by
$\rho^T=\sigma_{TM}\circ T\rho$, $\sigma_{TM}:T(TM)\rightarrow
T(TM)$ being the canonical involution of the double tangent
bundle. The Lie bracket $\lcf\cdot,\cdot\rcf^T$ on the space
$\Gamma(T\tau)$ is characterized by the following equalities
\begin{equation}\label{LieT}
\lcf TX,TY\rcf^T=T\lcf X,Y\rcf,\;\;\lcf
TX,\hat{Y}\rcf^T=\widehat{\lcf X,Y\rcf},\;\;\lcf
\hat{X},\hat{Y}\rcf^T=0, \end{equation}
 for $X,Y\in\Gamma(\tau)$
(see \cite{LMM,MX}).

Furthermore, there exists a unique action
$\Psi:\Gamma(T\tau)\rightarrow\frak X(E)$ of the Lie algebroid
$(TE,\lcf\cdot,\cdot\rcf^T,$ $\rho^T)$ over the anchor map
$\rho:E\rightarrow TM$ such that $\Psi(TX)=X^c$ and
$\Psi(\hat{X})=X^v,$ for $X\in\Gamma(\tau)$ (see \cite{LMM}).
Thus, the vector bundle $\rho^*(TE)$ is a Lie algebroid with Lie
algebroid structure
$(\lcf\cdot,\cdot\rcf_{\Psi}^T,\rho_{\Psi}^T)$, which is
characterized as in (\ref{2.3'}).

\subsection{Lie affgebroids}\label{sec1.2}

Let $\tau_A:A\rightarrow M$ be an affine bundle with associated
vector bundle $\tau_V:V\rightarrow M$. Denote by
$\tau_{A^+}:A^+=Aff(A,\R)\rightarrow M$ the dual bundle whose
fibre over $x\in M$ consists of affine functions on the  fibre
$A_x$. Note that this bundle has a distinguished section
$1_A\in\Gamma(\tau_{A^+})$ corresponding to the constant function
$1$ on $A$. We also consider the bidual bundle
$\tau_{\widetilde{A}}:\widetilde{A}\rightarrow M$ whose fibre at
$x\in M$ is the vector space $\widetilde{A}_x=(A_x^+)^*$. Then,
$A$ may be identified with an affine subbundle of $\widetilde{A}$
via the inclusion $i_A:A\rightarrow\widetilde{A}$ given by
$i_A(a)(\varphi)=\varphi(a)$, which is injective affine map whose
associated linear map is denoted by
$i_V:V\rightarrow\widetilde{A}$. Thus, $V$ may be identified with
a vector subbundle of $\widetilde{A}$. Using these facts, one can
prove that there is a one-to-one correspondence between affine
functions on $A$ and linear functions on $\widetilde{A}$. On the
other hand, there is an obvious one-to-one correspondence between
affine functions on $A$ and sections of $A^+$.

A {\it
 Lie affgebroid structure} on $A$ consists of a Lie algebra structure
 $\lcf\cdot,\cdot\rcf_V$ on the space
 $\Gamma(\tau_V)$ of the sections of $\tau_V:V\rightarrow M$, a $\R$-linear action
 $D:\Gamma(\tau_A)\times\Gamma(\tau_V)\rightarrow\Gamma(\tau_V)$ of
 the sections of $A$ on $\Gamma(\tau_V)$ and an affine map
 $\rho_A:A\rightarrow TM$, the {\it anchor map}, satisfying the following
 conditions:

\smallskip

  \begin{enumerate}
\item[$ \bullet$]
$D_X\lcf\bar{Y},\bar{Z}\rcf_V=\lcf
D_X\bar{Y},\bar{Z}\rcf_V+\lcf\bar{Y},D_X\bar{Z}\rcf_V,$
\item[$\bullet$]
$D_{X+\bar{Y}}\bar{Z}=D_X\bar{Z}+\lcf\bar{Y},\bar{Z}\rcf_V,$
\item[$\bullet$] $D_X(f\bar{Y})=fD_X\bar{Y}+\rho_A(X)(f)\bar{Y},$
\end{enumerate}

\smallskip

\noindent for $X\in\Gamma(\tau_A)$,
$\bar{Y},\bar{Z}\in\Gamma(\tau_V)$ and $f\in C^{\infty}(M)$ (see
\cite{GGrU,MMeS}).

If $(\lcf\cdot,\cdot\rcf_V,D,\rho_A)$ is a Lie affgebroid
structure on an affine bundle $A$ then
$(V,\lcf\cdot,\cdot\rcf_V,\rho_V)$ is a Lie algebroid, where
$\rho_V:V\rightarrow TM$ is the vector bundle map associated with
the affine morphism $\rho_A:A\rightarrow TM$.

A Lie affgebroid structure on an affine bundle
$\tau_A:A\rightarrow M$ induces a Lie algebroid structure
$(\lcf\cdot,\cdot\rcf_{\widetilde{A}},\rho_{\widetilde{A}})$ on
the bidual bundle $\widetilde{A}$ such that
$1_A\in\Gamma(\tau_{A^+})$ is a $1$-cocycle in the corresponding
Lie algebroid cohomology, that is, $d^{\widetilde{A}}1_A=0$.
Indeed, if $X_0\in\Gamma(\tau_A)$ then for every section
$\widetilde{X}$ of $\widetilde{A}$ there exists a function $f\in
C^{\infty}(M)$ and a section $\bar{X}\in\Gamma(\tau_V)$ such that
$\widetilde{X}=fX_0+\bar{X}$ and
\begin{equation}\label{*}
\begin{array}{rcl}
\rho_{\widetilde{A}}(fX_0+\bar{X})&=&f\rho_A(X_0)+\rho_V(\bar{X}),\\
\lcf
fX_0+\bar{X},gX_0+\bar{Y}\rcf_{\widetilde{A}}&=&(\rho_V(\bar{X})(g)-\rho_V(\bar{Y})(f)+f\rho_A(X_0)(g)\\&&-g\rho_A(X_0)(f))X_0
+\lcf\bar{X},\bar{Y}\rcf_V+fD_{X_0}\bar{Y}-gD_{X_0}\bar{X}.
\end{array}
\end{equation}

Conversely, let $(U,\lcf\cdot,\cdot\rcf_U,\rho_U)$ be a Lie
algebroid over $M$ and $\phi:U\rightarrow\R$ be a $1$-cocycle of
$(U,\lcf\cdot,\cdot\rcf_U,\rho_U)$ such that $\phi_{|U_x}\neq 0$,
for all $x\in M$. Then, $A=\phi^{-1}\{1\}$  is an affine bundle
over $M$ which admits a Lie affgebroid structure in such a way
that $(U,\lcf\cdot,\cdot\rcf_U,\rho_U)$ may be identified with the
bidual Lie algebroid
$(\widetilde{A},\lcf\cdot,\cdot\rcf_{\widetilde{A}},\rho_{\widetilde{A}})$
to $A$ and, under this identification, the $1$-cocycle
$1_A:\widetilde{A}\rightarrow\R$ is just $\phi$. The affine bundle
$\tau_A:A\rightarrow M$ is modelled on the vector bundle
$\tau_V:V=\phi^{-1}\{0\}\rightarrow M$. In fact, if
$i_V:V\rightarrow U$ and $i_A:A\rightarrow U$ are the canonical
inclusions, then
\begin{equation}\label{+}
\begin{array}{rcl} i_V\circ\lcf\bar{X},\bar{Y}\rcf_V&=&\lcf
i_V\circ\bar{X},i_V\circ\bar{Y}\rcf_U,\;\;i_V\circ D_X\bar{Y}=\lcf
i_A\circ X,i_V\circ\bar{Y}\rcf_U,\\
\rho_A(X)&=&\rho_{U}(i_A\circ X),
\end{array}
\end{equation}
 for $\bar{X},\bar{Y}\in\Gamma(\tau_V)$ and $X\in\Gamma(\tau_A).$

A trivial example of a Lie affgebroid may be constructed as
follows. Let $\tau:M\to\R$ be a fibration and
$\tau_{1,0}:J^1\tau\to M$ be the $1$-jet bundle of local sections
of $\tau:M\to\R$. It is well known that $\tau_{1,0}:J^1\tau\to M$
is an affine bundle modelled on the vector bundle
$\pi=(\pi_M)_{|V\tau}:V\tau\to M$, where $V\tau$ is the vertical
bundle of $\tau:M\to\R$. Moreover, if $t$ is the usual coordinate
on $\R$ and $\eta$ is the closed $1$-form on $M$ given by
$\eta=\tau^*(dt)$ then we have the following identification
$$J^1\tau\cong\{v\in TM/\eta(v)=1\}$$
(see, for instance, \cite{Sa}). Note that $V\tau=\{v\in
TM/\eta(v)=0\}.$ Thus, the bidual bundle $\widetilde{J^1\tau}$ to
the affine bundle $\tau_{1,0}:J^1\tau\to M$ may be identified with
the tangent bundle $TM$ to $M$ and, under this identification, the
Lie algebroid structure on $\pi_M:TM\to M$ is the standard Lie
algebroid structure and the $1$-cocycle $1_{J^1\tau}$ on
$\pi_M:TM\to M$ is just the $1$-form $\eta$.

Let $\tau_A:A\to M$ be a Lie affgebroid modelled on the Lie
algebroid $\tau_V:V\to M$. Suppose that $(x^i)$ are local
coordinates on an open subset $U$ of $M$ and that
$\{e_0,e_{\alpha}\}$ is a local basis of sections of
$\tau_{\widetilde{A}}:\widetilde{A}\to M$ in $U$ which is adapted
to the $1$-cocycle $1_A$, i.e., such that $1_A(e_0)=1$ and
$1_A(e_{\alpha})=0,$ for all $\alpha.$ Note that if
$\{e^0,e^{\alpha}\}$ is the dual basis of $\{e_0,e_{\alpha}\}$
then $e^0=1_A$. Denote by $(x^i,y^0,y^{\alpha})$ the corresponding
local coordinates on $\widetilde{A}$. Then, the local equation
defining the affine subbundle $A$ (respectively, the vector
subbundle $V$) of $\widetilde{A}$ is $y^0=1$ (respectively,
$y^0=0$). Thus, $(x^i,y^{\alpha})$ may be considered as local
coordinates on $A$ and $V$.

Now, let $\tau_A:A\rightarrow M$ (respectively,
$\tau_{A'}:A'\rightarrow M'$) be an affine bundle with associated
vector bundle $\tau_V:V\rightarrow M$ (respectively,
$\tau_{V'}:V'\rightarrow M'$) and $F:A\to A'$ be an affine bundle
morphism over the map $f:M\to M'$ with associated morphism
$(F^l,f)$ between the vector bundles $\tau_V:V\rightarrow M$ and
$\tau_{V'}:V'\rightarrow M'.$

Then, a direct computation proves that the map
$\widetilde{F}:\widetilde{A}\rightarrow\widetilde{A'}$  given by
$$\widetilde{F}(\widetilde{a})(\varphi')=\widetilde{a}(\varphi'\circ
F),\;\;\;\mbox{ for $\widetilde{a}\in\widetilde{A}_x$ and
$\varphi'\in(A')^+_{f(x)}$, with $x\in M$,}$$ defines a morphism
between the vector bundles $\widetilde{A}$ and $\widetilde{A'}$
over $f$ and, moreover, $(\widetilde{F},f)^*1_{A'}=1_A.$

Conversely, suppose that $\tau_U:U\rightarrow M$ and
$\tau_{U'}:U'\rightarrow M'$ are vector bundles and that $\phi$
and $\phi'$ are sections of the vector bundles
$\tau^*_U:U^*\rightarrow M$ and $\tau^*_{U'}:(U')^*\rightarrow M'$
such that $\phi(x)\neq 0$, for all $x\in M$, and $\phi'(x')\neq
0$, for all $x'\in M'$. Assume also that the pair
$(\widetilde{F},f)$ is a morphism between the vector bundles
$\tau_U:U\rightarrow M$ and $\tau_{U'}:U'\rightarrow M'$ such that
$(\widetilde{F},f)^*\phi'=\phi$ and denote by $A$ and $V$
(respectively, $A'$ and $V'$) the subsets of $U$ (respectively,
$U'$) defined by $A=\phi^{-1}\{1\}$ and $V=\phi^{-1}\{0\}$
(respectively, $A'=(\phi')^{-1}\{1\}$ and $V'=(\phi')^{-1}\{0\}$).
Then, it is easy to prove that $\widetilde{F}(A)\subseteq A'$ and
$\widetilde{F}(V)\subseteq V'$. Thus, the pair $(F,f)$ is a
morphism between the affine bundles
$\tau_A=(\tau_U)_{|A}:A\rightarrow M$ and
$\tau_{A'}=(\tau_{U'})_{|A'}:A'\rightarrow M'$, where
$F:A\rightarrow A'$ is the restriction of $\widetilde{F}$ to $A$.
The corresponding morphism between the vector bundles
$\tau_V=(\tau_U)_{|V}:V\rightarrow M$ and
$\tau_{V'}=(\tau_{U'})_{|V'}:V'\rightarrow M'$ is the pair
$(F^l,f)$, $F^l:V\rightarrow V'$ being the restriction of
$\widetilde{F}$ to $V$.

Now, suppose that $(\tau_A:A\rightarrow M$, $\tau_V:V\rightarrow
M, (\lcf\cdot,\cdot\rcf_V,D,\rho_A))$ and
$(\tau_{A'}:A'\rightarrow M'$, $\tau_{V'}:V'\rightarrow M',
(\lcf\cdot,\cdot\rcf_{V'},D',\rho_{A'}))$ are two Lie affgebroids
and that $((F,f),(F^l,f))$ is a morphism between the affine
bundles ($\tau_A:A\rightarrow M$, $\tau_V:V\rightarrow M$) and
($\tau_{A'}:A'\rightarrow M'$, $\tau_{V'}:V'\rightarrow M'$).
Then, the pair $((F,f),(F^l,f))$ is said to be a {\it Lie
affgebroid morphism} if:
\medskip

\noindent $(i)$ The pair $(F^l,f)$ is a morphism between the Lie
algebroids $(V,\lcf\cdot,\cdot\rcf_V,$  $\rho_V)$ and
$(V',\lcf\cdot,\cdot\rcf_{V'},$ $\rho_{V'})$,

\medskip

\noindent $(ii)$ $Tf\circ\rho_A=\rho_{A'}\circ F$ and

\medskip

\noindent $(iii)$ If $X$ (respectively, $X'$) is a section of
$\tau_A:A\rightarrow M$ (respectively, $\tau_{A'}:A'\rightarrow
M'$) and $\bar{Y}$ (respectively, $\bar{Y}'$) is a section of
$\tau_V:V\rightarrow M$ (respectively, $\tau_{V'}:V'\rightarrow
M'$) such that $X'\circ f=F\circ X$ and $\bar{Y}'\circ
f=F^l\circ\bar{Y}$ then $F^l\circ
D_X\bar{Y}=(D'_{X'}\bar{Y}')\circ f.$

\medskip

Now, using (\ref{*}), one may deduce the following result.
\begin{proposition} Suppose that $(\tau_A:A\rightarrow M, \tau_V:V\rightarrow
M)$ and $(\tau_{A'}:A'\rightarrow M', \tau_{V'}:V'\rightarrow M')$
are Lie affgebroids.
 If
$((F,f),(F^l,f))$ is a Lie affgebroid morphism and
$\widetilde{F}:\widetilde{A}\rightarrow\widetilde{A'}$ is the
corresponding morphism between the bidual vector bundles
$\widetilde{A}$ and $\widetilde{A'},$  then the pair
$(\widetilde{F},f)$ is a morphism between the Lie algebroids
$\widetilde{A}$ and $\widetilde{A'}$.
\end{proposition}

Finally, using (\ref{+}), one may prove
\begin{proposition}\label{prop2.2} Suppose that $\tau_U:U\rightarrow M$ and $\tau_{U'}:U'\rightarrow
M'$ are Lie algebroids and that $\phi\in\Gamma(\tau_U^*)$ and
$\phi'\in\Gamma(\tau_{U'}^*)$ are $1$-cocycles of
$\tau_U:U\rightarrow M$ and $\tau_{U'}:U'\rightarrow M'$,
respectively, such that $\phi(x)\neq 0$, for all $x\in M$, and
$\phi'(x')\neq 0$, for all $x'\in M'$. Then, if the pair
$(\widetilde{F},f)$ is a Lie algebroid morphism between the Lie
algebroids $\tau_U:U\rightarrow M$ and $\tau_{U'}:U'\rightarrow
M'$ satisfying $(\widetilde{F},f)^*\phi'=\phi$, we have that the
corresponding morphism $((F,f),(F^l,f))$ between the Lie
affgebroids ($\tau_A=(\tau_U)_{|A}:A=\phi^{-1}\{1\}\rightarrow M$,
$\tau_V=(\tau_U)_{|V}:V=\phi^{-1}\{0\}\rightarrow M$) and
($\tau_{A'}=(\tau_{U'})_{|A'}:A'=(\phi')^{-1}\{1\}\rightarrow M'$,
$\tau_{V'}=(\tau_U')_{|V'}:V'=(\phi')^{-1}\{0\}\rightarrow M'$) is
a Lie affgebroid morphism.
\end{proposition}


 \setcounter{equation}{0}
\section{Hamiltonian and Lagrangian formalism on Lie
affgebroids}

\subsection{The Hamiltonian formalism}\label{sec3.1}

In this section, we will develop a geometric framework, which
allows to write the Hamilton equations associated with a
Hamiltonian section on a Lie affgebroid (see \cite{M}).

Suppose that $(\tau_A:A\rightarrow M, \tau_V:V\rightarrow M,
(\lcf\cdot,\cdot\rcf_V,D,\rho_A))$ is a Lie affgebroid. Now, let
$({\mathcal
L}^{\tau_{V}^*}\widetilde{A},\lcf\cdot,\cdot\rcf_{\widetilde{A}}^{\tau_{V}^*},\rho_{\widetilde{A}}^{\tau_{V}^*})$
be the prolongation of the bidual Lie algebroid
$(\widetilde{A},\lcf\cdot,\cdot\rcf_{\widetilde{A}},$ $
\rho_{\widetilde{A}})$ over the fibration
$\tau_{V}^*:V^*\rightarrow M$.

Let $(x^i)$ be local coordinates on an open subset $U$ of $M$ and
$\{e_0,e_{\alpha}\}$ be a local basis of sections of the vector
bundle $\tau^{-1}_{\widetilde{A}}(U)\rightarrow U$ adapted to
$1_A$  and
$$
\begin{array}{rclrclccrclrcl}
\lcf
e_0,e_{\alpha}\rcf_{\widetilde{A}}=C_{0\alpha}^{\gamma}e_{\gamma},\;\;\lcf
e_{\alpha},e_{\beta}\rcf_{\widetilde{A}}=C_{\alpha\beta}^{\gamma}e_{\gamma},\;\;
\rho_{\widetilde{A}}(e_0)=\rho_0^i\displaystyle\frac{\partial}{\partial
x^i},\;\;
\rho_{\widetilde{A}}(e_{\alpha})=\rho_{\alpha}^i\displaystyle\frac{\partial}{\partial
x^i}.
\end{array}
$$

Denote by $(x^i,y^0,y^{\alpha})$ the corresponding local
coordinates on $\widetilde{A}$ and by $(x^i,y_0,y_{\alpha})$ the
dual coordinates on the dual vector bundle $\tau_{A^+}:A^+\to M$
to $\widetilde{A}$. Then, $(x^i,y_{\alpha})$ are local coordinates
on $V^*$ and $\{\tilde{e}_0,\tilde{e}_{\alpha},\bar{e}_{\alpha}\}$
is a local basis of sections of the vector bundle
$\tau^{\tau_{V}^*}_{\widetilde{A}}:{\mathcal
L}^{\tau_{V}^*}\widetilde{A}\rightarrow V^*$ , where
$$\tilde{e}_0(\psi)=(e_0(\tau_V^*(\psi)),\rho_0^i\displaystyle\frac{\partial}{\partial
x^i}_{|\psi}),\;\;\tilde{e}_{\alpha}(\psi)=(e_{\alpha}(\tau_V^*(\psi)),\rho_{\alpha}^i\displaystyle\frac{\partial}{\partial
x^i}_{|\psi}),\;\;\bar{e}_{\alpha}(\psi)=(0,\displaystyle\frac{\partial}{\partial
y_{\alpha}}_{|\psi}).$$ Using this local basis one may introduce
local coordinates $(x^i,y_{\alpha};z^0,$ $z^{\alpha},v_{\alpha})$
on ${\mathcal  L}^{\tau_{V}^*}\widetilde{A}$.

Let $\mu:A^+\rightarrow V^*$ be the canonical projection given by
$\mu(\varphi)=\varphi^l$, for $\varphi\in A^+_x$, with $x\in M$,
where $\varphi^l\in V^*_x$ is the linear map associated with the
affine map $\varphi$ and $h:V^*\rightarrow A^+$ be a Hamiltonian
section of $\mu$.

Now, we consider the Lie algebroid prolongation ${\mathcal
L}^{\tau_{A^+}}\widetilde{A}$ of the Lie algebroid $\widetilde{A}$
over $\tau_{A^+}:A^+\to M$ with vector bundle projection
$\tau^{\tau_{A^+}}_{\widetilde{A}}:{\mathcal
L}^{\tau_{A^+}}\widetilde{A}\rightarrow A^+$ (see Section
\ref{sec1.1.1}). Then, we may introduce the map ${\mathcal
L}h:{\mathcal  L}^{\tau_V^*}\widetilde{A}\rightarrow{\mathcal
L}^{\tau_{A^+}}\widetilde{A}$ defined by ${\mathcal
L}h(\tilde{a},X_{\alpha})=(\tilde{a},(T_{\alpha}h)(X_{\alpha})),$
for $(\tilde{a},X_{\alpha})\in({\mathcal
L}^{\tau_V^*}\widetilde{A})_{\alpha}$, with $\alpha\in V^*.$ It is
easy to prove that the pair $({\mathcal  L}h,h)$ is a Lie
algebroid morphism between the Lie algebroids
$\tau_{\widetilde{A}}^{\tau_V^*}:{\mathcal
L}^{\tau_V^*}\widetilde{A}\rightarrow V^*$ and
$\tau_{\widetilde{A}}^{\tau_{A^+}}:{\mathcal
L}^{\tau_{A^+}}\widetilde{A}\rightarrow A^+$.

\noindent Next, denote by $\lambda_h$ and $\Omega_h$ the sections
of the vector bundles $({\mathcal
L}^{\tau_V^*}\widetilde{A})^*\rightarrow V^*$ and
$\Lambda^2({\mathcal  L}^{\tau_V^*}\widetilde{A})^*\rightarrow
V^*$ given by
\begin{equation}\label{Omegah}
\lambda_h=({\mathcal
L}h,h)^*(\lambda_{\widetilde{A}}),\;\;\Omega_h=({\mathcal
L}h,h)^*(\Omega_{\widetilde{A}}), \end{equation}

where $\lambda_{\widetilde{A}}$ and $\Omega_{\widetilde{A}}$ are
the Liouville section and  the canonical symplectic section,
respectively, associated with the Lie algebroid $\widetilde{A}.$
Note that $\Omega_h=-d^{{\mathcal
L}^{\tau_V^*}\widetilde{A}}\lambda_h.$

On the other hand, let $\eta:{\mathcal
L}^{\tau_V^*}\widetilde{A}\rightarrow\R$ be the section of
$({\mathcal  L}^{\tau_V^*}\widetilde{A})^*\rightarrow V^*$ defined
by
\begin{equation}\label{eta}
\eta(\tilde{a},X_{\alpha})=1_A(\tilde{a}),
\end{equation}
 for
$(\tilde{a},X_{\alpha})\in({\mathcal
L}^{\tau_V^*}\widetilde{A})_{\alpha}$, with $\alpha\in V^*$. Note
that if $pr_1:{\mathcal  L}^{\tau_V^*}\widetilde{A}\to
\widetilde{A}$ is the canonical projection on the first factor
then $(pr_1,\tau_V^*)$ is a morphism between the Lie algebroids
$\tau_{\widetilde{A}}^{\tau_V^*}:{\mathcal
L}^{\tau_V^*}\widetilde{A}\rightarrow V^*$ and
$\tau_{\widetilde{A}}:\widetilde{A}\to M$ and
$(pr_1,\tau_V^*)^*(1_A)=\eta$. Thus, since $1_A$ is a $1$-cocycle
of $\tau_{\widetilde{A}}:\widetilde{A}\rightarrow M$, we deduce
that $\eta$ is a $1$-cocycle of the Lie algebroid
$\tau_{\widetilde{A}}^{\tau_V^*}:{\mathcal
L}^{\tau_V^*}\widetilde{A}\rightarrow V^*.$

Suppose that
$h(x^i,y_{\alpha})=(x^i,-H(x^j,y_{\beta}),y_{\alpha})$ and that
$\{\tilde{e}^0,\tilde{e}^{\alpha},\bar{e}^{\alpha}\}$ is the dual
basis of $\{\tilde{e}_0,\tilde{e}_{\alpha},$ $\bar{e}_{\alpha}\}$.
Then $\eta=\tilde{e}^0$ and, from (\ref{formas}), (\ref{Omegah})
and the definition of the map ${\mathcal  L}h$, it follows that
\begin{equation}\label{Omegahloc}
\Omega_h=\tilde{e}^{\gamma}\wedge\bar{e}^{\gamma}+\frac{1}{2}C_{\gamma\beta}^{\alpha}y_{\alpha}\tilde{e}^{\gamma}\wedge\tilde{e}^{\beta}+(\rho_{\gamma}^i\frac{\partial
H}{\partial
x^i}-C_{0\gamma}^{\alpha}y_{\alpha})\tilde{e}^{\gamma}\wedge\tilde{e}^0+\frac{\partial
H }{\partial y_{\gamma}}\bar{e}^{\gamma}\wedge\tilde{e}^0.
\end{equation}

Thus, it is easy to prove that the pair $(\Omega_h,\eta)$ is a
cosymplectic structure on the Lie algebroid
$\tau_{\widetilde{A}}^{\tau_V^*}:{\mathcal
L}^{\tau_V^*}\widetilde{A}\rightarrow V^*$, that is,
$$
\begin{array}{c}
\{\eta\wedge\Omega_h\wedge\dots^{(n}
\dots\wedge\Omega_h\}(\alpha)\neq 0,\makebox[1.5cm]{for
all}\alpha\in V^*,\\
d^{{\mathcal  L}^{\tau_V^*}\widetilde{A}}\eta=0,\;\;\;d^{{\mathcal
L}^{\tau_V^*}\widetilde{A}}\Omega_h=0.
\end{array}
$$

\begin{remark} {\em Let ${\mathcal  L}^{\tau_V^*}V$ be the prolongation of
the Lie algebroid $V$ over the projection $\tau_V^*:V^*\to M$.
Denote by $\lambda_V$ and $\Omega_V$ the Liouville section and the
canonical symplectic section, respectively, of $V$ and by
$(i_V,Id):{\mathcal  L}^{\tau_V^*}V\to{\mathcal
L}^{\tau_V^*}\widetilde{A}$ the canonical inclusion. Then, using
(\ref{lambdaE}), (\ref{Omegah}), (\ref{eta}) and the fact that
$\mu\circ h=Id$, we obtain that
\begin{equation}\label{3.2'}
(i_V,Id)^*(\lambda_h)=\lambda_V,\;\;\;(i_V,Id)^*(\eta)=0.
\end{equation}
Thus, since $(i_V,Id)$ is a Lie algebroid morphism, we also deduce
that
\begin{equation}\label{inclusion2}
(i_V,Id)^*(\Omega_h)=\Omega_V.
\end{equation}
\hfill$\diamondsuit$}
\end{remark}

Now, let $R_h\in\Gamma(\tau_{\widetilde{A}}^{\tau_V^*})$ be the
Reeb section of the cosymplectic structure $(\Omega_h,\eta)$
characterized by the following conditions
\begin{equation}\label{reeb}
i_{R_h}\Omega_h=0\makebox[1cm]{and}i_{R_h}\eta=1.
\end{equation}
With respect to the basis
$\{\tilde{e}_0,\tilde{e}_{\alpha},\bar{e}_{\alpha}\}$ of
$\Gamma(\tau_{\widetilde{A}}^{\tau_V^*})$, $R_h$ is locally
expressed as follows:
\begin{equation}\label{rh}
R_h=\tilde{e}_0+\frac{\partial H}{\partial
y_{\alpha}}\tilde{e}_{\alpha}-(C_{\alpha\beta}^{\gamma}y_{\gamma}\frac{\partial
H}{\partial y_{\beta}}+\rho^i_{\alpha}\frac{\partial H}{\partial
x^i}-C_{0\alpha}^{\gamma}y_{\gamma})\bar{e}_{\alpha}.
\end{equation}

Thus, the vector field $\rho_{\widetilde{A}}^{\tau_{V}^*}(R_h)$ is
locally given by
\begin{equation}\label{eq3.7'}
\rho_{\widetilde{A}}^{\tau_{V}^*}(R_h)=(\rho_0^i+\frac{\partial H
}{\partial y_{\alpha}}\rho_{\alpha}^i)\frac{\partial}{\partial
x^i}+\big(-\rho_{\alpha}^i\frac{\partial H}{\partial
x^i}+y_\gamma(C_{0\alpha}^{\gamma}+C_{\beta\alpha}^{\gamma}\frac{\partial
H}{\partial y_{\beta}})\big)\frac{\partial}{\partial y_{\alpha}}
\end{equation}
and the integral curves of $R_h$ (i.e., the integral curves of
$\rho_{\widetilde{A}}^{\tau_{V}^*}(R_h)$) are just {\it the
solutions of the Hamilton equations for $h$},
\begin{equation}\label{ech}
\frac{dx^i}{dt}=\rho_0^i+\frac{\partial H}{\partial
y_{\alpha}}\rho_{\alpha}^i,\;\;\;\frac{dy_{\alpha}}{dt}=-\rho_{\alpha}^i\frac{\partial
H }{\partial
x^i}+y_{\gamma}(C_{0\alpha}^{\gamma}+C_{\beta\alpha}^{\gamma}\frac{\partial
H }{\partial y_{\beta}}),
\end{equation}
for $i\in\{1,\dots,m\}$ and $\alpha\in\{1,\dots,n\}$.

Next, we will present an alternative approach in order to obtain
the Hamilton equations. For this purpose, we will use the notion
of an aff-Poisson structure on an AV-bundle which was introduced
in \cite{GGrU} (see also \cite{GGU2}).

Let $\tau_Z:Z\to M$ be an affine bundle of rank $1$ modelled on
the trivial vector bundle $\tau_{M\times\R}:M\times\R\to M$, that
is, $\tau_Z:Z\to M$ is an {\em AV-bundle } in the terminology of
\cite{GGU2}.

Then, we have an action of $\R$ on the fibers of $Z$. This action
induces a vector field $X_Z$ on $Z$ which is vertical respect to
the projection $\tau_Z:Z\to M$.

On the other hand, there exists a one-to-one correspondence
between the space of sections of $\tau_Z:Z\to M$,
$\Gamma(\tau_Z)$, and the set
$$\{F_h\in C^{\infty}(Z)/X_Z(F_h)=-1\}.$$
In fact, if $h\in\Gamma(\tau_Z)$ and $(x^i,s)$ are local fibred
coordinates on $Z$ such that
$X_Z=\displaystyle\frac{\partial}{\partial s}$ then $h$ may be
considered as a local function $H$ on $M$, $x^i\to H(x^i)$, and
the function $F_h$ on $Z$ is locally given by
\begin{equation}\label{eq3.81}
F_h(x^i,s)=-H(x^i)-s,
\end{equation}
(for more details, see \cite{GGU2}).

Now, an {\em aff-Poisson structure } on the affine bundle
$\tau_Z:Z\to M$ is a bi-affine map
$$\{\cdot,\cdot\}:\Gamma(\tau_Z)\times\Gamma(\tau_Z)\to
C^{\infty}(M)$$ which satisfies the following properties:
\begin{enumerate}
\item[i)] Skew-symmetric: $\{h_1,h_2\}=-\{h_2,h_1\}$.
\item[ii)] Jacobi identity:
$$\{h_1,\{h_2,h_3\}\}_V^2+\{h_2,\{h_3,h_1\}\}_V^2+\{h_3,\{h_1,h_2\}\}_V^2=0,$$
where $\{\cdot,\cdot\}_V^2$ is the affine-linear part of the
bi-affine bracket.
\item[iii)] If $h\in\Gamma(\tau_Z)$ then
$$\{h,\cdot\}:\Gamma(\tau_Z)\to
C^{\infty}(M),\;\;\;h'\mapsto\{h,h'\},$$ is an affine derivation.
\end{enumerate}
Condition $iii)$ implies that the linear part
$\{h,\cdot\}_V:C^{\infty}(M)\to C^{\infty}(M)$ of the affine map
$\{h,\cdot\}:\Gamma(\tau_Z)\to C^{\infty}(M)$ defines a vector
field on $M$, which is called {\it the Hamiltonian vector field of
$h$} (see \cite{GGU2}).

In \cite{GGU2}, the authors proved that there is a one-to-one
correspondence between aff-Poisson brackets $\{\cdot,\cdot\}$ on
$\tau_Z:Z\to M$ and Poisson brackets $\{\cdot,\cdot\}_{\Pi}$ on
$Z$ which are $X_Z$-invariant, i.e., which are associated with
Poisson 2-vectors $\Pi$ on $Z$ such that ${\mathcal
L}_{X_Z}\Pi=0$. This correspondence is determined by
$$\{h,h'\}\circ\tau_Z=\{F_h,F_{h'}\}_\Pi,\makebox[1cm]{for}h,h'\in\Gamma(\tau_Z).$$
Note that the function $\{F_h,F_{h'}\}_\Pi$ on $Z$ is
$\tau_Z$-projectable, i.e., ${\mathcal L}_{X_Z}\{F_h,F_{h'}\}_\Pi$
$=0$ (because the Poisson 2-vector $\Pi$ is $X_Z$-invariant).

Using this correspondence we will prove the following result.

\begin{theorem}\label{teor3.2} Let $\tau_A:A\to M$ be a Lie affgebroid modelled
on the vector bundle $\tau_V:V\to M$. Denote by $\tau_{A^+}:A^+\to
M$ (resp., $\tau_V^*:V^*\to M$) the dual vector bundle to $A$
(resp., to $V$) and by $\mu:A^+\to V^*$ the canonical projection.
Then:
\begin{enumerate}
\item[i)] $\mu:A^+\to V^*$ is an AV-bundle which admits an
aff-Poisson structure.
\item[ii)] If $h:V^*\to A^+$ is a Hamiltonian section (that is,
$h\in\Gamma(\mu)$) then the Hamiltonian vector field of $h$ with
respect to the aff-Poisson structure is a vector field on $V^*$
whose integral curves are just the solutions of the Hamilton
equations for $h$.
\end{enumerate}
\end{theorem}
\begin{proof}$i)$ It is clear that $\mu:A^+\to V^*$ is an
AV-bundle. In fact, if $a^+\in A^+_x$, with $x\in M$, and $t\in\R$
then
$$a^++t=a^++t1_A(x).$$
Thus, the $\mu$-vertical vector field $X_{A^+}$ on $A^+$ is just
the vertical lift $1_A^V$ of the section
$1_A\in\Gamma(\tau_{A^+})$. Moreover, one may consider the Lie
algebroid $\tau_{\widetilde{A}}:\widetilde{A}=(A^+)^*\to M$ and
the corresponding linear Poisson 2-vector $\Pi_{A^+}$ on $A^+$.
Then, using the fact that $1_A$ is a 1-cocycle of
$\tau_{\widetilde{A}}:\widetilde{A}=(A^+)^*\to M$, it follows that
the Poisson 2-vector $\Pi_{A^+}$ is $X_{A^+}$-invariant.
Therefore, $\Pi_{A^+}$ induces an aff-Poisson structure
$\{\cdot,\cdot\}$ on $\mu:A^+\to V^*$ which is characterized by
the condition
\begin{equation}\label{eq3.82}
\{h,h'\}\circ\mu=\{F_h,F_{h'}\}_{\Pi_{A^+}},\makebox[1cm]{for}h,h'\in\Gamma(\mu).
\end{equation}
One may also prove this first part of the theorem using the
relation between special Lie affgebroid structures on an affine
bundle $A'$ and aff-Poisson structures on the AV-bundle
$AV((A')^\sharp)$ (see Theorem $23$ in \cite{GGU2}).

$ii)$ From (\ref{eq3.81}) and (\ref{eq3.82}), we deduce that the
linear map $\{h,\cdot\}_V:C^\infty(V^*)\to C^\infty(V^*)$
associated with the affine map $\{h,\cdot\}:\Gamma(\mu)\to
C^\infty(V^*)$ (that is, the Hamiltonian vector field of $h$) is
given by
\begin{equation}\label{eq3.83}
\{h,\cdot\}_V(\varphi)\circ\mu=\{F_h,\varphi\circ\mu\}_{\Pi_{A^+}},\makebox[1cm]{for}\varphi\in
C^\infty(V^*).
\end{equation}
Now, suppose that the local expression of $h$ is
\begin{equation}\label{eq3.84}
h(x^i,y_\alpha)=(x^i,-H(x^j,y_\beta),y_\alpha).
\end{equation}
On the other hand, using (\ref{eq2.3'}), we have that
\begin{equation}\label{eq3.85}
\Pi_{A^+}=\displaystyle\frac{1}{2}C_{\alpha\beta}^\gamma
y_\gamma\frac{\partial}{\partial
y_\alpha}\wedge\frac{\partial}{\partial
y_\beta}+C_{0\alpha}^\gamma y_\gamma\frac{\partial}{\partial
y_0}\wedge\frac{\partial}{\partial
y_\alpha}+\rho_0^i\displaystyle\frac{\partial}{\partial
y_0}\wedge\frac{\partial}{\partial
x^i}+\rho_\alpha^i\displaystyle\frac{\partial}{\partial
y_\alpha}\wedge\frac{\partial}{\partial x^i}.
\end{equation}
Thus, from (\ref{eq3.83}), (\ref{eq3.84}) and (\ref{eq3.85}), we
conclude that the Hamiltonian vector field of $h$ is locally given
by
$$(\rho_0^i+\frac{\partial H
}{\partial y_{\alpha}}\rho_{\alpha}^i)\frac{\partial}{\partial
x^i}+(-\rho_{\alpha}^i\frac{\partial H}{\partial
x^i}+y_\gamma(C_{0\alpha}^{\gamma}+C_{\beta\alpha}^{\gamma}\frac{\partial
H}{\partial y_{\beta}}))\frac{\partial}{\partial y_{\alpha}}$$
which proves our result (see (\ref{eq3.7'})).
\end{proof}

\subsection{The Lagrangian formalism}\label{sec3.2}

In this section, we will develop a geometric framework, which
allows to write the Euler-Lagrange equations associated with a
Lagrangian function $L$ on a Lie affgebroid $A$ in an intrinsic
way (see \cite{MMeS}).

Suppose that $(\tau_A:A\rightarrow M, \tau_V:V\rightarrow
M,(\lcf\cdot,\cdot\rcf_V,D,\rho_A))$  is a Lie affgebroid on $M$.
Then, the bidual bundle $\tau_{\widetilde{A}}:\widetilde{A}\to M$
to $A$ admits a Lie algebroid structure
$(\lcf\cdot,\cdot\rcf_{\widetilde{A}},\rho_{\widetilde{A}})$ in
such a way that the section $1_A$ of the dual bundle $A^+$ is a
$1$-cocycle.

Now, we consider the Lie algebroid prolongation (${\mathcal
L}^{\tau_{A}}\widetilde{A},$
$\lcf\cdot,\cdot\rcf_{\widetilde{A}}^{\tau_A},\rho_{\widetilde{A}}^{\tau_A})$
of the Lie algebroid
$(\widetilde{A},\lcf\cdot,\cdot\rcf_{\widetilde{A}},$ $
\rho_{\widetilde{A}})$ over the fibration $\tau_{A}:A\rightarrow
M$ with vector bundle projection
$\tau^{\tau_{A}}_{\widetilde{A}}:{\mathcal
L}^{\tau_{A}}\widetilde{A}\rightarrow A$.

If $(x^i)$ are local coordinates on an open subset $U$ of $M$ and
$\{e_0,e_{\alpha}\}$ is a local basis of sections of the vector
bundle $\tau^{-1}_{\widetilde{A}}(U)\rightarrow U$ adapted to
$1_A$, then
$\{\tilde{T}_0,\tilde{T}_{\alpha},\tilde{V}_{\alpha}\}$ is a local
basis of sections of the vector bundle
$(\tau^{\tau_{A}}_{\widetilde{A}})^{-1}(\tau^{-1}_{\widetilde{A}}(U))\rightarrow\tau^{-1}_{\widetilde{A}}(U)$,
where
\begin{equation}\label{tildeT}
\tilde{T}_0(a)=(e_0(\tau_A(a)),\rho_0^i\displaystyle\frac{\partial}{\partial
x^i}_{|a}),\;\;\;
\tilde{T}_{\alpha}(a)=(e_{\alpha}(\tau_A(a)),\rho_{\alpha}^i\displaystyle\frac{\partial}{\partial
x^i}_{|a}),\;\;\;
\tilde{V}_{\alpha}(a)=(0,\displaystyle\frac{\partial}{\partial
y^{\alpha}}_{|a}),
\end{equation}
$(x^i,y^{\alpha})$ are the local coordinates on $A$ induced by the
local coordinates $(x^i)$ and the basis $\{e_{\alpha}\}$ and
$\rho_0^i,\;\rho_{\alpha}^i$ are the components of the anchor map
$\rho_{\widetilde{A}}$. Therefore, we have that
\begin{equation}\label{corTVtilde}
\begin{array}{c}
\lcf\tilde{T}_0,\tilde{T}_{\alpha}\rcf_{\widetilde{A}}^{\tau_{A}}=C_{0\alpha}^{\gamma}\tilde{T}_{\gamma},\;\;\;
\lcf\tilde{T}_{\alpha},\tilde{T}_{\beta}\rcf_{\widetilde{A}}^{\tau_{A}}=C_{\alpha\beta}^{\gamma}\tilde{T}_{\gamma},\\[8pt]
\lcf\tilde{T}_0,\tilde{V}_{\alpha}\rcf_{\widetilde{A}}^{\tau_{A}}=\lcf\tilde{T}_{\alpha},\tilde{V}_{\beta}\rcf_{\widetilde{A}}^{\tau_{A}}=\lcf\tilde{V}_{\alpha},\tilde{V}_{\beta}\rcf_{\widetilde{A}}^{\tau_{A}}=0,\\[10pt]
\rho_{\widetilde{A}}^{\tau_{A}}(\tilde{T}_0)=\rho_0^i\displaystyle\frac{\partial}{\partial
x^i},\;\;\;\rho_{\widetilde{A}}^{\tau_{A}}(\tilde{T}_{\alpha})=\rho_{\alpha}^i\displaystyle\frac{\partial}{\partial
x^i},\;\;\;\rho_{\widetilde{A}}^{\tau_{A}}(\tilde{V}_{\alpha})=\displaystyle\frac{\partial}{\partial
y^{\alpha}},
\end{array}
\end{equation}
where $C_{0\beta}^{\gamma}$ and $C_{\alpha\beta}^{\gamma}$ are the
structure functions of the Lie bracket
$\lcf\cdot,\cdot\rcf_{\widetilde{A}}$ with respect to the basis
$\{e_0,e_{\alpha}\}$.  Note that, if
$\{\tilde{T}^0,\tilde{T}^{\alpha},\tilde{V}^{\alpha}\}$ is the
dual basis of
$\{\tilde{T}_0,\tilde{T}_{\alpha},\tilde{V}_{\alpha}\}$, then
$\tilde{T}^0$ is globally defined and it is a $1$-cocycle. We will
denote by $\phi_0$ the $1$-cocycle $\tilde{T}^0$. Thus, we have
that
\begin{equation}\label{T0}
\phi_0(a)(\tilde{b},X_a)=1_A(\tilde{b}),\mbox{ for
}(\tilde{b},X_a)\in({\mathcal  L}^{\tau_A}\widetilde{A})_a.
\end{equation}

One may also consider {\it the vertical endomorphism} $S:{\mathcal
L}^{\tau_A}\widetilde{A}\rightarrow{\mathcal
L}^{\tau_A}\widetilde{A}$, as a section of the vector bundle
${\mathcal  L}^{\tau_A}\widetilde{A}\otimes({\mathcal
L}^{\tau_A}\widetilde{A})^*\to A$, whose local expression is (see
\cite{MMeS})
\begin{equation}\label{verend}
S=(\tilde{T}^{\alpha}-y^{\alpha}\tilde{T}^0)\otimes\tilde{V}_{\alpha}.
\end{equation}

Now, a curve $\gamma:I\subseteq\R\rightarrow A$ in $A$ is said to
be {\it admissible} if $\rho_{\widetilde{A}}\circ
i_A\circ\gamma=\dot{\widehat{(\tau_A\circ\gamma)}}$
or,equivalently, $(i_A(\gamma(t)),\dot{\gamma}(t))\in({\mathcal
L}^{\tau_A}\widetilde{A})_{\gamma(t)}$, for all $t\in I$ ,
$i_A:A\rightarrow\widetilde{A}$ being the canonical inclusion.
Thus, if $\gamma(t)=(x^i(t),y^{\alpha}(t)),$ for all $t\in I$,
then $\gamma$ is an admissible curve if and only if
$$\displaystyle\frac{dx^i}{dt}=\rho_0^i+\rho_{\alpha}^iy^{\alpha},\makebox[1cm]{for}i\in\{1,\dots,m\}.$$

A section $\xi$ of $\tau^{\tau_{A}}_{\widetilde{A}}:{\mathcal
L}^{\tau_{A}}\widetilde{A}\rightarrow A$ is said to be a {\it
second order differential equation} (SODE) on $A$ if the integral
curves of $\xi$, that is, the integral curves of the vector field
$\rho^{\tau_A}_{\widetilde{A}}(\xi)$, are admissible.

If $\xi\in\Gamma(\tau^{\tau_A}_{\widetilde{A}})$ is a SODE then
$\xi=\tilde{T}_0+y^{\alpha}\tilde{T}_{\alpha}+\xi^{\alpha}\tilde{V}_{\alpha},$
where $\xi^{\alpha}$ are arbitrary local functions on $A$, and
$$\rho^{\tau_A}_{\widetilde{A}}(\xi)=(\rho_0^i+y^{\alpha}\rho_{\alpha}^i)\displaystyle\frac{\partial}{\partial
x^i}+\xi^{\alpha}\displaystyle\frac{\partial}{\partial
y^{\alpha}}.$$

On the other hand, let $L:A\rightarrow\R$ be a Lagrangian
function. Then, we introduce {\it the Poincar\'{e}-Cartan $1$-section}
$\Theta_L\in\Gamma((\tau^{\tau_A}_{\widetilde{A}})^*)$ and {\it
the Poincar\'{e}-Cartan $2$-section}
$\Omega_L\in\Gamma(\wedge^2(\tau_{\widetilde{A}}^{\tau_A})^*)$
a\-sso\-cia\-ted with $L$ defined by
\begin{equation}\label{PC}
\Theta_L=L\phi_0+(d^{{\mathcal  L}^{\tau_A}\widetilde{A}}L)\circ
S,\;\;\;\Omega_L=-d^{{\mathcal  L}^{\tau_A}\widetilde{A}}\Theta_L.
\end{equation}

From (\ref{corTVtilde}), (\ref{verend}) and (\ref{PC}), we obtain
that
\begin{equation}\label{PCloc}
\begin{array}{lcl}
\Theta_L&=&(L-y^{\alpha}\displaystyle\frac{\partial L}{\partial
y^{\alpha}})\tilde{T}^0+\displaystyle\frac{\partial L}{\partial
y^{\alpha}}\tilde{T}^{\alpha},\\[8pt]
\Omega_L&=&(i_{\xi_0}(d^{{\mathcal
L}^{\tau_A}\widetilde{A}}(\displaystyle\frac{\partial L}{\partial
y^{\alpha}}))-\displaystyle\frac{\partial L}{\partial
y^{\gamma}}(C_{0\alpha}^{\gamma}+C_{\beta\alpha}^{\gamma}y^{\beta})-\rho_{\alpha}^i\displaystyle\frac{\partial
L}{\partial
x^i})\theta^{\alpha}\wedge\tilde{T}^0\\[8pt]
&+&\displaystyle\frac{\partial^2 L}{\partial y^{\alpha}\partial
y^{\beta}}\theta^{\alpha}\wedge\psi^{\beta}
+\displaystyle\frac{1}{2}(\rho_{\beta}^i\displaystyle\frac{\partial^2
L}{\partial x^i\partial
y^{\alpha}}-\rho_{\alpha}^i\displaystyle\frac{\partial^2
L}{\partial x^i\partial y^{\beta}}+\displaystyle\frac{\partial
L}{\partial
y^{\gamma}}C_{\alpha\beta}^{\gamma})\theta^{\alpha}\wedge\theta^{\beta},
\end{array}
\end{equation}
where $\theta^{\alpha}=\tilde{T}^{\alpha}-y^{\alpha}\tilde{T}^0$,
$\psi^{\alpha}=\tilde{V}^{\alpha}-\xi_0^{\alpha}\tilde{T}^0$ and
$\xi_0=\tilde{T}_0+y^{\alpha}\tilde{T}_{\alpha}+\xi_0^{\alpha}\tilde{V}_{\alpha}$
is an arbitrary SODE.

Now, a curve $\gamma:I=(-\epsilon,\epsilon)\subseteq\R\rightarrow
A$ in $A$ is {\it a solution of the Euler-Lagrange equations}
associated with $L$ if and only if $\gamma$ is admissible and
$i_{(i_A(\gamma(t)),\dot{\gamma}(t))}\Omega_L(\gamma(t))=0$, for
all $t$.

If $\gamma(t)=(x^i(t),y^{\alpha}(t))$ then $\gamma$ is a solution
of the Euler-Lagrange equations if and only if
\begin{equation}\label{EL}
\displaystyle\frac{dx^i}{dt}=\rho_0^i+\rho_{\alpha}^iy^{\alpha},\;\;\;\displaystyle\frac{d}{dt}(\displaystyle\frac{\partial
L}{\partial
y^{\alpha}})=\rho_{\alpha}^i\displaystyle\frac{\partial
L}{\partial
x^i}+(C_{0\alpha}^{\gamma}+C_{\beta\alpha}^{\gamma}y^{\beta})\displaystyle\frac{\partial
L}{\partial y^{\gamma}},
\end{equation}
for $i\in\{1,\dots,m\}$ and $\alpha\in\{1,\dots,n\}$.

The Lagrangian $L$ is {\it regular} if and only if the matrix
$(W_{\alpha\beta})=(\displaystyle\frac{\partial^2L}{\partial
y^{\alpha}\partial y^{\beta}})$ is regular or, in other words, the
pair $(\Omega_L,\phi_0)$ is a cosymplectic structure on ${\mathcal
L}^{\tau_A}\widetilde{A}$.

If the Lagrangian $L$ is regular, then the Reeb section $R_L$ of
$(\Omega_L,\phi_0)$ is the unique Lagrangian SODE associated with
$L$, that is, the integral curves of the vector field
$\rho^{\tau_A}_{\widetilde{A}}(R_L)$ are solutions of the
Euler-Lagrange equations associated with $L$. In such a case,
$R_L$ is called {\it the Euler-Lagrange section associated with
$L$} and its local expression is
\begin{equation}\label{RLloc}
R_L=\tilde{T}_0+y^{\alpha}\tilde{T}_{\alpha}+W^{\alpha\beta}(\rho_{\beta}^i\displaystyle\frac{\partial
L}{\partial
x^i}-(\rho_0^i+y^{\gamma}\rho_{\gamma}^i)\frac{\partial^2L}{\partial
x^i\partial
y^{\beta}}+(C_{0\beta}^{\gamma}+y^{\mu}C_{\mu\beta}^{\gamma})\frac{\partial
L}{\partial y^{\gamma}})\tilde{V}_{\alpha},
\end{equation}
where $(W^{\alpha\beta})$ is the inverse matrix of
$(W_{\alpha\beta})$.

\subsection{The Legendre transformation and the equivalence
between the Lagrangian and Hamiltonian formalisms}\label{sec3.3}

Let $L:A\to \R$ be a Lagrangian function and $\Theta_L\in
\Gamma((\tau^{\tau_A}_{\widetilde{A}})^*)$ be the Poincar\'{e}-Cartan
$1$-section associated with $L$. We introduce {\it the extended
Legendre transformation associated with $L$} as the smooth map
$Leg_L:A\to A^+$ defined by
$Leg_L(a)(b)=\Theta_L(a)(z),$
for $a,b\in A_x,$ where $z$ is a point in the fibre of ${\mathcal
L}^{\tau_A}\widetilde{A}$ over the point $a$ such that
$pr_1(z)=i_A(b),$ $pr_1:{\mathcal
L}^{\tau_A}\widetilde{A}\to\widetilde{A}$ being the restriction to
${\mathcal  L}^{\tau_A}\widetilde{A}$ of the first canonical
projection $pr_1:\widetilde{A}\times TA\to\widetilde{A}$.

The map $Leg_L$ is well-defined and its local expression in fibred
coordinates on $A$ and  $A^+$ is
\begin{equation}\label{locLegL}
Leg_L(x^i,y^\alpha)=(x^i,L-\frac{\partial L}{\partial
y^{\alpha}}y^{\alpha},\frac{\partial L}{\partial y^{\alpha}}).
\end{equation}

Thus, we can define {\it the Legendre transformation associated
with $L$}, $leg_L:A\to V^*$, by $leg_L=\mu\circ Leg_L.$ From
(\ref{locLegL}) and since
$\mu(x^i,y_0,y_{\alpha})=(x^i,y_{\alpha})$, we have that
\begin{equation}\label{loclegL}
leg_L(x^i,y^{\alpha})=(x^i,\displaystyle\frac{\partial L}{\partial
y^\alpha}).
\end{equation}

The maps $Leg_L$ and $leg_L$ induce the maps ${\mathcal
L}{Leg_L}:{\mathcal  L}^{\tau_A}\widetilde{A}\to {\mathcal
L}^{\tau_{A^+}}\widetilde{A}$ and ${\mathcal  L}{leg_L}:{\mathcal
L}^{\tau_A}\widetilde{A}\to {\mathcal  L}^{\tau_V^*}\widetilde{A}$
defined by
\begin{equation}\label{LLegL}
({\mathcal
L}{Leg_L})(\tilde{b},X_a)=(\tilde{b},(T_aLeg_L)(X_a)),\;\;({\mathcal
L}{leg_L})(\tilde{b},X_a)=(\tilde{b},(T_aleg_L)(X_a)),
\end{equation}
for $a\in A$ and $(\tilde{b},X_a)\in ({\mathcal
L}^{\tau_A}\widetilde{A})_a$.

Now, let $\{\tilde{T}_0,\tilde{T}_{\alpha},\tilde{V}_{\alpha}\}$
(respectively,
$\{\tilde{e}_0,\tilde{e}_{\alpha},\bar{e}_0,\bar{e}_{\alpha}\}$)
be a local basis of $\Gamma(\tau_{\widetilde{A}}^{\tau_A})$ as in
Section \ref{sec3.2} (respectively, of
$\Gamma(\tau_{\widetilde{A}}^{\tau_{A^+}})$ as in Section
\ref{sec1.1.1}) and denote by
$(x^i,y^{\alpha};z^0,z^{\alpha},v^{\alpha})$ (respectively,
$(x^i,y_0,y_{\alpha};z^0,z^{\alpha},v_0,v_{\alpha})$) the
corresponding local coordinates on ${\mathcal
L}^{\tau_A}\widetilde{A}$ (respectively, ${\mathcal
L}^{\tau_{A^+}}\widetilde{A}$). In addition, suppose that
$(x^i,y_{\alpha};z^0,z^{\alpha},v_{\alpha})$ are local coordinates
on ${\mathcal  L}^{\tau_{V^*}}\widetilde{A}$ as in Section
\ref{sec3.1}. Then, using (\ref{locLegL}), (\ref{loclegL}) and
(\ref{LLegL}), we deduce that the local expression of the maps
${\mathcal  L}{Leg_L}$ and ${\mathcal  L}leg_L$ is
\begin{equation}\label{LLegloc}
\begin{array}{rcl}
{\mathcal
L}{Leg_L}(x^i,y^\alpha;z^0,z^\alpha,v^\alpha)\kern-10pt&=&\kern-10pt(x^i,L-\displaystyle\frac{\partial
L}{\partial y^{\alpha}}y^{\alpha},\frac{\partial L}{\partial
y^\alpha};z^0,z^\alpha,z^0\rho_0^i(\frac{\partial L}{\partial
x^i}-\frac{\partial^2L}{\partial x^i\partial
y^{\alpha}}y^{\alpha})\\[8pt]&&
\kern-20pt+z^\alpha\rho_\alpha^i(\displaystyle\frac{\partial
L}{\partial x^i}-\frac{\partial^2 L}{\partial x^i\partial
y^\beta}y^{\beta})- v^\alpha\frac{\partial^2 L}{\partial
y^\alpha\partial y^\beta}y^{\beta}, z^0\rho_0^i\frac{\partial^2
L}{\partial x^i\partial
y^{\alpha}}\\[8pt]&&\kern-20pt+z^{\beta}\rho_{\beta}^i\displaystyle\frac{\partial^2L}{\partial
x^i\partial
y^{\alpha}}+v^{\beta}\displaystyle\frac{\partial^2L}{\partial
y^{\alpha}\partial y^{\beta}} ),
\end{array}\end{equation}
\begin{equation}\label{Llegloc}
\begin{array}{rcl}
 {\mathcal
L}{leg_L}(x^i,y^\alpha;z^0,z^\alpha,v^\alpha)&=&(x^i,\displaystyle\frac{\partial
L}{\partial
y^\alpha};z^0,z^{\alpha},z^0\rho_0^i\displaystyle\frac{\partial^2L}{\partial
x^i\partial
y^{\alpha}}+z^{\beta}\rho_{\beta}^i\frac{\partial^2L}{\partial
x^i\partial
y^{\alpha}}\\[8pt] &&+v^{\beta}\displaystyle\frac{\partial^2L}{\partial
y^{\alpha}\partial y^{\beta}}).  \end{array}\end{equation} Thus,
using (\ref{formas}), (\ref{PCloc}), (\ref{LLegL}) and
(\ref{LLegloc}), we can prove the following result.
\begin{theorem}\label{t3.2}
The pair $({\mathcal  L}{Leg_L},Leg_L)$ is a morphism between the
Lie algebroids $({\mathcal  L}^{\tau_A}\widetilde{A}, \linebreak
\lcf\cdot,\cdot\rcf^{\tau_A}_{\widetilde{A}},\rho^{\tau_A}_{\widetilde{A}})$
and $({\mathcal
L}^{\tau_{A^+}}\widetilde{A},\lcf\cdot,\cdot\rcf^{\tau_{A^+}}_{\widetilde{A}},\rho^{\tau_{A^+}}_{\widetilde{A}}).$
Moreover, if $\Theta_L$ and $\Omega_L$ (respectively,
$\lambda_{\widetilde{A}}$ and $\Omega_{\widetilde{A}})$ are the
Poincar\'{e}-Cartan $1$-section and $2$-section associated with $L$
(respectively, the Liouville $1$-section and the canonical
symplectic section associated with $\widetilde{A}$) then
\begin{equation}\label{pullback}
({\mathcal
L}{Leg_L},Leg_L)^*(\lambda_{\widetilde{A}})=\Theta_L,\;\;\;
({\mathcal  L}{Leg_L},Leg_L)^*(\Omega_{\widetilde{A}})=\Omega_L.
\end{equation}
\end{theorem}
From (\ref{loclegL}), it follows
\begin{proposition}
The Lagrangian $L$ is regular if and only if the Legendre
transformation $leg_L:A\to V^*$ is a local diffeomorphism.
\end{proposition}

Next, we will assume that $L$ is {\it hyperregular}, that is,
$leg_L$ is a global diffeomorphism. Then, from (\ref{LLegL}) and
Theorem \ref{t3.2}, we conclude that the pair $({\mathcal
L}leg_L,leg_L)$ is a Lie algebroid isomorphism. Moreover, we can
consider the Hamiltonian section $h_L:V^*\to A^+$ defined by
\begin{equation}\label{h}
h_L=Leg_L\circ leg_L^{-1}, \end{equation} the corresponding
cosymplectic structure $(\Omega_{h_L},\eta)$ on the Lie algebroid
$\tau_{\widetilde{A}}^{\tau_V^*}:{\mathcal
L}^{\tau_V^*}\widetilde{A}\to V^*$ and the Hamiltonian section
$R_{h_L}\in\Gamma(\tau_{\widetilde{A}}^{\tau_V^*})$.

Using (\ref{Omegah}), (\ref{Llegloc}), (\ref{pullback}), (\ref{h})
and Theorem \ref{t3.2}, we deduce that

\begin{theorem}\label{t3.4}
If the Lagrangian $L$ is hyperregular then the Euler-Lagrange
section $R_L$ asso\-cia\-ted with $L$ and the Hamiltonian section
$R_{h_L}$ associated with $h_L$ satisfy the following relation
\begin{equation}\label{Related}
R_{h_L}\circ leg_L={\mathcal  L}leg_L\circ R_L.
\end{equation}

Moreover, if $\gamma:I\to A$ is a solution of the Euler-Lagrange
equations associated with $L$, then $leg_L\circ \gamma:I\to V^*$
is a solution of the Hamilton equations associated with $h_L$ and,
conversely, if $\bar{\gamma}:I\to V^*$ is a solution of the
Hamilton equations for $h_L$ then
$\gamma=leg_L^{-1}\circ\bar{\gamma}$ is a solution of the
Euler-Lagrange equations for $L$.
\end{theorem}

Now, we will analyze the local expression of the transformation
$leg_L^{-1}:V^*\to A$. Suppose that
\begin{equation}\label{3.24-1}
leg_L^{-1}(x^i,y_\alpha)=(x^i,y^\alpha(x^j,y_\beta))
\end{equation}
and $h_L(x^i,y_\alpha)=(x^i,-H_L(x^j,y_\beta),y_\alpha).$ Then,
from (\ref{locLegL}) and (\ref{h}), it follows that
\[
H_L(x^i,y_\alpha)=y^\alpha(x^j,y_\beta)y_\alpha-L(x^i,y^\alpha(x^j,y_\beta)).
\]
Thus, we obtain that
\[
\frac{\partial H_L}{\partial y_\alpha}=y^\alpha + \frac{\partial
y^\beta}{\partial y_\alpha}y_\beta -\frac{\partial L}{\partial
y^\beta}\frac{\partial y^\beta}{\partial y_\alpha}
\]
and, using (\ref{loclegL}) and (\ref{3.24-1}), we deduce that
$\displaystyle\frac{\partial H_L}{\partial y_\alpha}=y^\alpha.$
Therefore, we conclude that
\begin{equation}\label{3.24-2}
leg_L^{-1}(x^i,y_\alpha)=(x^i,\frac{\partial H_L}{\partial
y_\alpha}).
\end{equation}
Note that, since $leg_L^{-1}:V^*\to A$ is a diffeomorphism, it
follows that the matrix $\big(\displaystyle\frac{\partial^2
H_L}{\partial y_\alpha\partial y_\beta}\big)$ is regular.

Next, we will introduce the notion of a hyperregular Hamiltonian
section and we will prove that given a hyperregular Hamiltonian
section $h:V^*\to A^+$ then one may construct a hyperregular
Lagrangian function $L:A\to \R$ and $h_L=h.$

Let $h:V^*\to A^+$ be a Hamiltonian section. If ${\mathcal R}$ is
an arbitrary section of $A$, we may consider the real
$C^\infty$-function $H_{\mathcal R}:V^*\to \R$ on $V^*$ given by
\begin{equation}\label{3.24-3}
H_{\mathcal R}(\alpha)=h(\alpha)({\mathcal
R}(\tau_{V^*}(\alpha))), \mbox{ for } \alpha \in V^*.
\end{equation}
Using the function $H_{\mathcal R}$ we may define the map $({\Bbb
F}h)_{\mathcal R}:V^*\to V$ as follows
\[
\alpha\in V_x^*, \mbox{ with } x\in M \Rightarrow ({\Bbb
F}h)_{\mathcal R}(\alpha)\in V_x
\]
and
\begin{equation}\label{3.24-4}
\beta(({\Bbb F}h)_{\mathcal
R}(\alpha))=\frac{d}{dt}_{|t=0}H_{\mathcal R}(\alpha + t \beta),
\mbox{ for } \beta\in V_x^*.
\end{equation}
Now, we introduce the map ${\Bbb F}h:V^*\to A$ given by
\begin{equation}\label{3.24-5}
({\Bbb F}h)(\alpha)={\mathcal R}(\tau_{V^*}(\alpha))+({\Bbb
F}h)_{\mathcal R}(\alpha).
\end{equation}
If the local expressions of $h$ and ${\mathcal R}$ are
\[
h(x^i,y_\alpha)=(x^i, -H(x^j,y_\beta), y_\alpha),\;\;\;\;\;
{\mathcal R}(x^i)=(x^i,{\mathcal R}^\alpha(x^i)),
\]
then, from (\ref{3.24-3}), (\ref{3.24-4}) and (\ref{3.24-5}), we
obtain that
\begin{equation}\label{3.24-6}
H_{\mathcal R}(x^i,y_\alpha)=-H(x^i,y_\alpha)+{\mathcal
R}^\alpha(x^i)y_\alpha,
\end{equation}
\begin{equation}\label{3.24-7}
({\Bbb F}h)_{\mathcal R}(x^i,y_\alpha)=(x_i,-\frac{\partial
H}{\partial y_\alpha}(x^j,y_\beta) + {\mathcal R}^\alpha(x^i)),
\end{equation}
\begin{equation}\label{3.24-8}
({\Bbb F}h)(x^i,y_\alpha)=(x^i,\frac{\partial H}{\partial
y_\alpha}).
\end{equation}

Thus, the map ${\Bbb F}h$ doesn't depend on the chosen section
${\mathcal R}.$

The Hamiltonian section $h$ is said to be {\it regular } if the
map ${\Bbb F}h:V^*\to A$ is a local diffeormorphism or
equivalently if the matrix $\big(\displaystyle\frac{\partial^2
H}{\partial y_\alpha\partial y_\beta}\big)$ is regular. $h$ is
said to be {\it hyperregular } if ${\Bbb F}h:V^*\to A$ is a global
diffeomorphism.

It is clear that if $L:A\to \R$ is a hyperregular Lagrangian
function and $h_L:V^*\to A^+$ is the Hamiltonian section
associated with $L$ then $h_L$ is hyperregular. In fact, from
(\ref{3.24-2}) and (\ref{3.24-8}), it follows that ${\Bbb F}h_L$
is a diffeomorphism and ${\Bbb F}h_L=leg_L^{-1}.$

Now, we will prove that the converse is also true.

\begin{theorem}
If $h:V^*\to A^+$ is a hyperregular Hamiltonian section then there
exists a hyperregular Lagrangian function $L:A\to \R$ such that
the Hamiltonian section associated with $L$ is just $h$. In other
words, $h_L=h.$
\end{theorem}

\begin{proof} We define the Lagrangian function $L:A\to \R$ by
\[
L(a)=({\Bbb F}h)^{-1}(a)(a-{\mathcal R}(\tau_A(a)))+ H_{\mathcal
R}(({\Bbb F}h)^{-1}(a)),\;\;\; \mbox{ for $a\in A$}.
\]
The function $L$ doesn't depend on the chosen section ${\mathcal
R}$. In fact, if $({\Bbb F}
h)^{-1}(x^i,y^\alpha)=(x^i,y_\alpha(x^j,y^\beta))$ then, using
(\ref{3.24-6}) and (\ref{3.24-8}), we have that
\begin{equation}\label{3.24-9}
L(x^i,y^\alpha)=y^\alpha
y_\alpha(x^j,y^\beta)-H(x^i,y_\alpha(x^j,y^\beta)).
\end{equation}
Therefore,  we deduce that
\[
\frac{\partial L}{\partial y^\alpha}=y_\alpha +
y^\beta\frac{\partial y_\beta}{\partial y^\alpha}-\frac{\partial
H}{\partial y_\beta}\frac{\partial y_\beta}{\partial y^\alpha} \]
and using that $\displaystyle\frac{\partial H}{\partial
y_\beta}=y^\beta$ (see (\ref{3.24-8})), it follows that
\begin{equation}\label{3.24-10}
\frac{\partial L}{\partial y^\alpha}=y_\alpha.
\end{equation}
This implies that (see (\ref{loclegL})) $leg_L=({\Bbb F} h)^{-1}$
and consequently, from (\ref{locLegL}), (\ref{3.24-9}),
(\ref{3.24-10}) and since $h_L=Leg_L\circ leg_L^{-1},$ we conclude
that $h_L=h.$
\end{proof}

\setcounter{equation}{0}
\section{The canonical involution associated with a Lie
affgebroid}\label{sec3}

Let $(\tau_A:A\rightarrow M, \tau_V:V\rightarrow M)$ be a Lie
affgebroid.  Denote by $\rho_V:V\rightarrow TM$ the anchor map of
the Lie algebroid $\tau_V:V\rightarrow M$, by
$(\lcf\cdot,\cdot\rcf_{\widetilde{A}},\rho_{\widetilde{A}})$ the
Lie algebroid structure on the bidual bundle
$\tau_{\widetilde{A}}:\widetilde{A}\rightarrow M$ to $A$ and by
$1_A:\widetilde{A}\rightarrow\R$ the distinguished $1$-cocycle on
$\widetilde{A}$.

We consider the subset ${\mathcal J}^AA$ of the product manifold
$A\times TA$ defined by
$${\mathcal J}^AA=\{(a,v)\in A\times TA/\rho_A(a)=(T\tau_A)(v)\}.$$
Next, we will see that ${\mathcal J}^AA$ admits two Lie affgebroid
structures. We will also see that these Lie affgebroid structures
are isomorphic under the so-called canonical involution associated
with $A$.

$i)$ {\it The first structure}:

Let $({\mathcal
L}^{\tau_A}\widetilde{A},\lcf\cdot,\cdot\rcf^{\tau_A}_{\widetilde{A}},\rho^{\tau_A}_{\widetilde{A}})$
be the prolongation of the Lie algebroid
$(\widetilde{A},\lcf\cdot,\cdot\rcf_{\widetilde{A}},\rho_{\widetilde{A}})$
over the fibration $\tau_A:A\rightarrow M$ and $\phi_0:{\mathcal
L}^{\tau_A}\widetilde{A}\rightarrow\R$ be the $1$-cocycle of the
Lie algebroid cohomology complex of $({\mathcal
L}^{\tau_A}\widetilde{A},\lcf\cdot,\cdot\rcf^{\tau_A}_{\widetilde{A}},\rho^{\tau_A}_{\widetilde{A}})$
given by (\ref{T0}). Using (\ref{T0}) and the fact that
$(1_A)_{|\widetilde{A}_x}\neq 0$, for all $x\in M$, we deduce that
$(\phi_0)_{|({\mathcal L}^{\tau_A}\widetilde{A})_a}\neq 0$, for
all $a\in A$.

Moreover, we have that
\begin{equation}\label{1atil1}
(\phi_0)^{-1}\{1\}=\{(\tilde{a},v)\in\widetilde{A}\times
TA/\rho_{\widetilde{A}}(\tilde{a})=(T\tau_A)(v),1_A(\tilde{a})=1\}={\mathcal
J}^AA.
\end{equation}
In addition, if ${\mathcal L}^{\tau_A}V$ is the prolongation of
the Lie algebroid $(V,\lcf\cdot,\cdot\rcf_V,\rho_V)$ over the
fibration $\tau_A:A\rightarrow M$ then, it is easy to prove that
\begin{equation}\label{1atil2}
(\phi_0)^{-1}\{0\}={\mathcal L}^{\tau_A}V.
\end{equation}

We will denote by
$(\lcf\cdot,\cdot\rcf^{\tau_A}_V,\rho^{\tau_A}_V)$ the Lie
algebroid structure on $\tau^{\tau_A}_V:{\mathcal
L}^{\tau_A}V\rightarrow A$.

From (\ref{1atil1}), we conclude that ${\mathcal J}^AA$ is an
affine bundle over $A$ with affine bundle projection
$\tau^{\tau_A}_A:{\mathcal J}^AA\rightarrow A$ defined by
$$\tau^{\tau_A}_A(a,v)=\pi_A(v),$$
and, moreover, the affine bundle $\tau^{\tau_A}_A:{\mathcal
J}^AA\rightarrow A$ admits a Lie affgebroid structure in such a
way that the bidual Lie algebroid to $\tau^{\tau_A}_A:{\mathcal
J}^AA\rightarrow A$ is just $({\mathcal
L}^{\tau_A}\widetilde{A},\lcf\cdot,\cdot\rcf^{\tau_A}_{\widetilde{A}},\rho^{\tau_A}_{\widetilde{A}})$.
Finally, using (\ref{1atil2}), it follows that the Lie affgebroid
$\tau^{\tau_A}_A:{\mathcal J}^AA\rightarrow A$ is modelled on the
Lie algebroid $({\mathcal
L}^{\tau_{A}}V,\lcf\cdot,\cdot\rcf^{\tau_A}_V,\rho^{\tau_A}_V)$.

\begin{remark}{\em Let ${\mathcal L}^{\tau_{\widetilde{A}}}\widetilde{A}$ be
the prolongation of the Lie algebroid
$(\widetilde{A},\lcf\cdot,\cdot\rcf_{\widetilde{A}},\rho_{\widetilde{A}})$
over the fibration $\tau_{\widetilde{A}}:\widetilde{A}\rightarrow
M$. Denote by $(Id,Ti_A):{\mathcal
L}^{\tau_{A}}\widetilde{A}\rightarrow{\mathcal
L}^{\tau_{\widetilde{A}}}\widetilde{A}$ the inclusion defined by
$$(Id,Ti_A)(\tilde{a},v_b)=(\tilde{a},(T_bi_A)(v_b)),$$
for $(\tilde{a},v_b)\in({\mathcal L}^{\tau_{A}}\widetilde{A})_b$,
with $b\in A$. Then, it is easy to prove that the pair
$((Id,Ti_A),i_A)$ is a Lie algebroid morphism. Thus, $({\mathcal
L}^{\tau_A}\widetilde{A},\lcf\cdot,\cdot\rcf^{\tau_A}_{\widetilde{A}},\rho^{\tau_A}_{\widetilde{A}})$
is a Lie subalgebroid of  $({\mathcal
L}^{\tau_{\widetilde{A}}}\widetilde{A},\lcf\cdot,\cdot\rcf^{\tau_{\widetilde{A}}}_{\widetilde{A}},$
$ \rho^{\tau_{\widetilde{A}}}_{\widetilde{A}})$. \vspace{-5pt}

\hfill$\diamondsuit$}
\end{remark}

$ii)$ {\it The second structure}:

As we know, the tangent bundle to $\widetilde{A}$,
$T\widetilde{A}$, is a Lie algebroid over $TM$ with vector bundle
projection $T\tau_{\widetilde{A}}:T\widetilde{A}\rightarrow TM$.

Now, we consider the subset $d_0(1_A)^0$ of $T\widetilde{A}$ given
by \begin{equation}\label{d1A0} d_0(1_A)^0=\{\tilde{v}\in
T\widetilde{A}/d_0(1_A)(\tilde{v})=0\}=\{\tilde{v}\in
T\widetilde{A}/\tilde{v}(1_A)=0\}.
\end{equation}
$d_0(1_A)^0$ is the total space of a vector subbundle of
$T\tau_{\widetilde{A}}:T\widetilde{A}\rightarrow TM$. More
precisely, suppose that $\tilde{X}\in\Gamma(\tau_{\widetilde{A}})$
and denote by $T\tilde{X}:TM\rightarrow T\widetilde{A}$ the
tangent map to $\tilde{X}$ and by $\hat{\tilde{X}}:TM\rightarrow
T\widetilde{A}$ the section of
$T\tau_{\widetilde{A}}:T\widetilde{A}\rightarrow TM$ defined by
(\ref{Xtil}). Then, using (\ref{d1A0}), we deduce the following
facts:
\begin{enumerate}
\item[$i)$] If $1_A(\tilde{X})=c$, with $c\in\R$, we have that
$T\tilde{X}(TM)\subseteq d_0(1_A)^0$ and, thus,
$T\tilde{X}:TM\rightarrow d_0(1_A)^0$ is a section of the vector
bundle
$(T\tau_{\widetilde{A}})_{|d_0(1_A)^0}:d_0(1_A)^0\rightarrow TM$.
\item[$ii)$] If $1_A(\tilde{X})=0$ it follows that
$\hat{\tilde{X}}(TM)\subseteq d_0(1_A)^0$ and, therefore,
$\hat{\tilde{X}}:TM\rightarrow d_0(1_A)^0$ is a section of the
vector bundle
$(T\tau_{\widetilde{A}})_{|d_0(1_A)^0}:d_0(1_A)^0\rightarrow TM$.
\end{enumerate}

In fact, if $\{e_0,e_{\alpha}\}$ is a local basis of
$\Gamma(\tau_{\widetilde{A}})$ adapted to $1_A$, then
$\{Te_0,Te_{\alpha},\hat{e}_{\alpha}\}$ is a local basis of
$\Gamma((T\tau_{\widetilde{A}})_{|d_0(1_A)^0})$. Consequently, the
canonical inclusion $i:d_0(1_A)^0\rightarrow T\widetilde{A}$ is a
monomorphism (over the identity of $TM$) between the vector
bundles
$(T\tau_{\widetilde{A}})_{|d_0(1_A)^0}:d_0(1_A)^0\rightarrow TM$
and $T\tau_{\widetilde{A}}:T\widetilde{A}\rightarrow TM$.
Moreover, using (\ref{LieT}) and the fact that
$1_A:\widetilde{A}\rightarrow\R$ is a $1$-cocycle of the Lie
algebroid
$(\widetilde{A},\lcf\cdot,\cdot\rcf_{\widetilde{A}},\rho_{\widetilde{A}})$,
we deduce that the Lie bracket
$\lcf\cdot,\cdot\rcf^T_{\widetilde{A}}$ on
$\Gamma(T\tau_{\widetilde{A}})$ restricts to a Lie bracket on the
space $\Gamma((T\tau_{\widetilde{A}})_{|d_0(1_A)^0})$. Therefore,
we have proved the following result.

\begin{proposition} The vector bundle $(T\tau_{\widetilde{A}})_{|d_0(1_A)^0}:d_0(1_A)^0\rightarrow
TM$ is a Lie algebroid and the canonical inclusion
$i:d_0(1_A)^0\rightarrow T\widetilde{A}$ is a monomorphism between
the Lie algebroids
$(T\tau_{\widetilde{A}})_{|d_0(1_A)^0}:d_0(1_A)^0\rightarrow TM$
and $T\tau_{\widetilde{A}}:T\widetilde{A}\rightarrow TM$. Thus,
$(T\tau_{\widetilde{A}})_{|d_0(1_A)^0}:d_0(1_A)^0\rightarrow TM$
is a Lie subalgebroid of
$T\tau_{\widetilde{A}}:T\widetilde{A}\rightarrow TM$.
\end{proposition}
Next, we consider the pull-back of the vector bundle
$(T\tau_{\widetilde{A}})_{|d_0(1_A)^0}:d_0(1_A)^0\rightarrow TM$
over the anchor map $\rho_A:A\rightarrow TM$, that is,
$$\rho_A^*(d_0(1_A)^0)=\{(a,\tilde{v})\in A\times
d_0(1_A)^0/\rho_A(a)=(T\tau_{\widetilde{A}})_{|d_0(1_A)^0}(\tilde{v})\}.$$
$\rho_A^*(d_0(1_A)^0)$ is a vector bundle over $A$ with vector
bundle projection
$$pr_1:\rho_A^*(d_0(1_A)^0)\rightarrow A,\;\;(a,\tilde{v})\mapsto
a.$$ On the other hand, we will denote by
$(i_A,i):\rho_A^*(d_0(1_A)^0)\rightarrow\rho_{\widetilde{A}}^*(T\widetilde{A})$
the monomorphism (over the canonical inclusion
$i_A:A\rightarrow\widetilde{A}$) between the vector bundles
$\rho_A^*(d_0(1_A)^0)\rightarrow A$ and
$\rho_{\widetilde{A}}^*(T\widetilde{A})\rightarrow\widetilde{A}$
defined by
$$(i_A,i)(a,\tilde{v})=(i_A(a),\tilde{v}),\makebox[1cm]{for}(a,\tilde{v})\in\rho_A^*(d_0(1_A)^0).$$
We recall that the vector bundle
$\rho_{\widetilde{A}}^*(T\widetilde{A})\rightarrow\widetilde{A}$
is an action Lie algebroid (see Section \ref{sec1.1.2}).
Furthermore, we have
\begin{proposition}\label{prop3.2}
$i)$ The vector bundle $\rho_A^*(d_0(1_A)^0)\rightarrow A$ is a
Lie algebroid over $A$ and the pair $((i_A,i),i_A)$ is a
monomorphism between the Lie algebroids
$\rho_A^*(d_0(1_A)^0)\rightarrow A$ and
$\rho_{\widetilde{A}}^*(T\widetilde{A})\rightarrow\widetilde{A}$.

\medskip

$ii)$ If $\pi_{\widetilde{A}}:T\widetilde{A}\to\widetilde{A}$ is
the canonical projection and
$\varphi_0:\rho_A^*(d_0(1_A)^0)\rightarrow\R$ is the linear map
given by
\begin{equation}\label{phi0}
\varphi_0(a,\tilde{v})=1_A(\pi_{\widetilde{A}}(\tilde{v})),
\end{equation}
 then $\varphi_0$ is
a $1$-cocycle of the Lie algebroid
$\rho_A^*(d_0(1_A)^0)\rightarrow A$ and
$\varphi_{0|\rho_A^*(d_0(1_A)^0)_a}\neq 0$, for all $a\in A$.

\end{proposition}
\begin{proof} $i)$ Let $\tilde{X}$ be a
section of $\tau_{\widetilde{A}}:\widetilde{A}\rightarrow M$ and
$\tilde{X}^c$ (respectively, $\tilde{X}^v$) be the complete lift
(respectively, the vertical lift) of $\tilde{X}$. If
$1_A(\tilde{X})=c$, with $c\in\R$, it follows that
$\tilde{X}^c_{|A}(1_A)=0$ and, thus, the restriction of
$\tilde{X}^c$ to $A$ is tangent to $A$. In addition, if
$1_A(\tilde{X})=0$ we obtain that $\tilde{X}^v_{|A}(1_A)=0$ and,
therefore, the restriction of $\tilde{X}^v$ to $A$ is tangent to
$A$.

Now, proceeding as in the proof of Theorem 4.4 in \cite{LMM}, we
deduce that there exists a unique action $\Psi_0$ of the Lie
algebroid
$(T\tau_{\widetilde{A}})_{|d_0(1_A)^0}:d_0(1_A)^0\rightarrow TM$
over the anchor map $\rho_A:A\rightarrow TM$ such that
$\Psi_0(T\tilde{X})=\tilde{X}^c_{|A},$ for
$\tilde{X}\in\Gamma(\tau_{\widetilde{A}})$ with $1_A(\tilde{X})=c$
and $c\in\R$, and $\Psi_0(\hat{\tilde{X}})=\tilde{X}^v_{|A},$ for
$\tilde{X}\in\Gamma(\tau_{\widetilde{A}})$ such that
$1_A(\tilde{X})=0$.

Note that if
$\Psi:\Gamma(T\tau_{\widetilde{A}})\rightarrow\mathfrak
X(\widetilde{A})$ is the usual action of the Lie algebroid
$T\tau_{\widetilde{A}}:T\widetilde{A}\rightarrow TM$ over the
anchor map $\rho_{\widetilde{A}}:\widetilde{A}\rightarrow TM$ (see
Section \ref{sec1.1.2}) and if
$j:\Gamma((T\tau_{\widetilde{A}})_{|d_0(1_A)^0})\to\Gamma(T\tau_{\widetilde{A}})$
is the canonical inclusion then
$$(\Psi(j(Z_0)))\circ i_A=(Ti_A)(\Psi_0(Z_0)),\makebox[1cm]{for}Z_0\in\Gamma((T\tau_{\widetilde{A}})_{|d_0(1_A)^0}).$$

Consequently, the pair $((i_A,i),i_A)$ is a monomorphism between
the Lie algebroids $\rho_A^*(d_0(1_A)^0)$ $\rightarrow A$ and
$\rho_{\widetilde{A}}^*(T\widetilde{A})\rightarrow\widetilde{A}$.

$ii)$ Let $\{e_0,e_{\alpha}\}$ be a local basis of
$\Gamma(\tau_{\widetilde{A}})$ adapted to $1_A.$ Then,
$$\{T^{\rho_A}e_0=Te_0\circ\rho_A,T^{\rho_A}e_{\alpha}=Te_{\alpha}\circ\rho_A,\hat{e}_{\alpha}^{\rho_A}=\hat{e}_{\alpha}\circ\rho_A\}$$
is a local basis of sections of the vector bundle
$\rho_A^*(d_0(1_A)^0)\rightarrow A$. Moreover, if
$((\lcf\cdot,\cdot\rcf_{\widetilde{A}})_{\Psi_0}^T,$ $
(\rho_{\widetilde{A}})_{\Psi_0}^T)$ is the Lie algebroid structure
on $\rho_A^*(d_0(1_A)^0)\rightarrow A$, we have that
\begin{equation}\label{Est}
\begin{array}{c} (\lcf
T^{\rho_A}e_0,T^{\rho_A}e_{\alpha}\rcf_{\widetilde{A}})_{\Psi_0}^T=T\lcf
e_0,e_{\alpha}\rcf_{\widetilde{A}}\circ\rho_A,\\[8pt] (\lcf
T^{\rho_A}e_{\alpha},\;T^{\rho_A}e_{\beta}\rcf_{\widetilde{A}})_{\Psi_0}^T=T\lcf
e_{\alpha},e_{\beta}\rcf_{\widetilde{A}}\circ\rho_A, \\[8pt]
(\lcf
T^{\rho_A}e_0,\hat{e}_{\alpha}^{\rho_A}\rcf_{\widetilde{A}})_{\Psi_0}^T=(\lcf
T^{\rho_A}e_{\alpha},\hat{e}_{\beta}^{\rho_A}\rcf_{\widetilde{A}})_{\Psi_0}^T=(\lcf
\hat{e}_{\alpha}^{\rho_A},\hat{e}_{\beta}^{\rho_A}\rcf_{\widetilde{A}})_{\Psi_0}^T=0,
\end{array}
\end{equation}
and
$$(\rho_{\widetilde{A}})_{\Psi_0}^T(T^{\rho_A}e_0)=(e_0)^c_{|A},\;\;
(\rho_{\widetilde{A}})_{\Psi_0}^T(T^{\rho_A}e_{\alpha})=(e_{\alpha})^c_{|A},\;\;
(\rho_{\widetilde{A}})_{\Psi_0}^T(\hat{e}_{\alpha}^{\rho_A})=(e_{\alpha})^v_{|A}.$$

On the other hand, if $\tilde{X}\in\Gamma(\tau_{\widetilde{A}})$
then, from (\ref{Xtil}) and (\ref{phi0}),we obtain that
\begin{equation}\label{phi1}
\begin{array}{rcl}
1_A(\tilde{X})=c\in\R&\Rightarrow&\varphi_0(T\tilde{X}\circ\rho_A)=c,\\
1_A(\tilde{X})=0&\Rightarrow&\varphi_0(\hat{\tilde{X}}\circ\rho_A)=0.
\end{array}\end{equation}

Thus, using (\ref{Est}), (\ref{phi1}) and the fact that $1_A$ is a
$1$-cocycle of
$(\widetilde{A},\lcf\cdot,\cdot\rcf_{\widetilde{A}},\rho_{\widetilde{A}})$,
we conclude that $\varphi_0$ is a $1$-cocycle of the Lie algebroid
$\rho_A^*(d_0(1_A)^0)\rightarrow A$.

\end{proof}

Now, from (\ref{phi0}), it follows that
\begin{equation}\label{phi0m}
\varphi_0^{-1}\{1\}=\{(a,v)\in A\times
TA/\rho_A(a)=(T\tau_A)(v)\}={\mathcal J}^AA.
\end{equation}
Therefore, we deduce that ${\mathcal J}^AA$ is an affine bundle
over $A$ with affine bundle projection $pr_1:{\mathcal
J}^AA\rightarrow A$ defined by $pr_1(a,v)=a$ and, moreover, the
affine bundle $pr_1:{\mathcal J}^AA\rightarrow A$ admits a Lie
affgebroid structure in such a way that the bidual Lie algebroid
to $pr_1:{\mathcal J}^AA\rightarrow A$ is just
$(\rho_A^*(d_0(1_A)^0),$ $
(\lcf\cdot,\cdot\rcf_{\widetilde{A}})_{\Psi_0}^T,(\rho_{\widetilde{A}})_{\Psi_0}^T)$.

On the other hand, using (\ref{phi0}), we obtain that
$$\varphi_0^{-1}\{0\}=\{(a,u)\in A\times
TV/\rho_A(a)=(T\tau_V)(u)\}=\rho_A^*(TV).$$ Consequently, the
affine bundle $pr_1:{\mathcal J}^AA\rightarrow A$ is modelled on
the vector bundle $pr_1:\rho_A^*(TV)\rightarrow A$. Furthermore,
using (\ref{Est}), we deduce that the corresponding Lie algebroid
structure is induced by an action $\Psi_V$ of the Lie algebroid
$(TV,\lcf\cdot,\cdot\rcf_V^T,\rho_V^T)$ over the anchor map
$\rho_A:A\rightarrow TM$. For this action, we have that
$$\Psi_V(T\bar{X})=(i_V\circ\bar{X})^c_{|A},\;\;\Psi_V(\hat{\bar{X}})=
(i_V\circ\bar{X})^v_{|A},\makebox[1cm]{for}\bar{X}\in\Gamma(\tau_V).$$

{\it The canonical involution:}

Let ${\mathcal L}^{\tau_{\widetilde{A}}}\widetilde{A}$ be the
prolongation of the Lie algebroid
$(\widetilde{A},\lcf\cdot,\cdot\rcf_{\widetilde{A}},\rho_{\widetilde{A}})$
over the fibration $\tau_{\widetilde{A}}:\widetilde{A}\rightarrow
M$ and $\rho^*_{\widetilde{A}}(T\widetilde{A})\equiv {\mathcal
L}^{\tau_{\widetilde{A}}}\widetilde{A}$ be the pull-back of the
Lie algebroid $T\tau_{\widetilde{A}}:T\widetilde{A}\rightarrow TM$
over the anchor map $\rho_{\widetilde{A}}:\widetilde{A}\rightarrow
TM$. If $(\tilde{a},\tilde{v}_{\tilde{b}})\in({\mathcal
L}^{\tau_{\widetilde{A}}}\widetilde{A})_{\tilde{b}}$, with
$\tilde{b}\in\widetilde{A}_x$ and $x\in M$, then there exists a
unique tangent vector $\tilde{u}_{\tilde{a}}\in
T_{\tilde{a}}\widetilde{A}$ such that:
\begin{equation}\label{Defsig}
\tilde{u}_{\tilde{a}}(f\circ\tau_{\widetilde{A}})=(d^{\widetilde{A}}f)(x)(\tilde{b}),\;\;\tilde{u}_{\tilde{a}}({\theta})=\tilde{v}_{\tilde{b}}({\theta})+(d^{\widetilde{A}}\theta)(x)(\tilde{b},\tilde{a}),
\end{equation}
for $f\in C^{\infty}(M)$ and $\theta:\widetilde{A}\to \R
\in\Gamma(\tau_{A^+})$ (see \cite{LMM}). Thus, one may define the
map $\sigma_{\widetilde{A}}:{\mathcal
L}^{\tau_{\widetilde{A}}}\widetilde{A}\rightarrow\rho^*_{\widetilde{A}}(T\widetilde{A})$
as follows
\begin{equation}\label{sigma}
\sigma_{\widetilde{A}}(\tilde{a},\tilde{v}_{\tilde{b}})=(\tilde{b},\tilde{u}_{\tilde{a}}),\makebox[1cm]{for}(\tilde{a},\tilde{v}_{\tilde{b}})\in({\mathcal
L}^{\tau_{\widetilde{A}}}\widetilde{A})_{\tilde{b}}.
\end{equation}
$\sigma_{\widetilde{A}}$ is an isomorphism (over the identity
$Id:\widetilde{A}\rightarrow\widetilde{A}$) between the Lie
algebroids $({\mathcal
L}^{\tau_{\widetilde{A}}}\widetilde{A},\lcf\cdot,\cdot\rcf_{\widetilde{A}}^{\tau_{\widetilde{A}}},$
$ \rho_{\widetilde{A}}^{\tau_{\widetilde{A}}})$ and
$(\rho^*_{\widetilde{A}}(T\widetilde{A}),(\lcf\cdot,\cdot\rcf_{\widetilde{A}})^T_{\Psi},(\rho_{\widetilde{A}})^T_{\Psi})$
and, moreover, $\sigma_{\widetilde{A}}^2=Id$.
$\sigma_{\widetilde{A}}$ is called {\it the canonical involution
a\-sso\-ciated with the Lie algebroid
$(\widetilde{A},\lcf\cdot,\cdot\rcf_{\widetilde{A}},\rho_{\widetilde{A}})$}
(for more details, see \cite{LMM}).

\begin{theorem} The restriction of $\sigma_{\widetilde{A}}$ to
${\mathcal J}^AA$ induces an isomorphism $\sigma_{A}:{\mathcal
J}^AA\rightarrow{\mathcal J}^AA$ between the Lie affgebroids
$\tau_A^{\tau_A}:{\mathcal J}^AA\rightarrow A$ and $pr_1:{\mathcal
J}^AA\rightarrow A$ and, moreover, $\sigma_A^2=Id$. The
corresponding Lie algebroid isomorphism $\sigma_A^l:{\mathcal
L}^{\tau_A}V\rightarrow\rho_A^*(TV)$ between the Lie algebroids
$\tau_V^{\tau_A}:{\mathcal L}^{\tau_A}V\rightarrow A$ and
$pr_1:\rho_A^*(TV)\rightarrow A$ is the restriction of
$\sigma_{\widetilde{A}}$ to ${\mathcal L}^{\tau_A}V$, that is,
$\sigma_A^l=(\sigma_{\widetilde{A}})_{|{\mathcal L}^{\tau_A}V}$.
\end{theorem}
\begin{proof} Suppose that $b$ is a point of $A$ and that
$(\tilde{a},v_b)\in({\mathcal
L}^{\tau_A}\widetilde{A})_b\subseteq({\mathcal
L}^{\tau_{\widetilde{A}}}\widetilde{A})_{i_A(b)}$. Then, from
(\ref{Defsig}) and (\ref{sigma}), it follows that
$\sigma_{\widetilde{A}}(\tilde{a},v_b)\in\rho_A^*(d_0(1_A)^0)_b$.

Furthermore, if
$(b,\tilde{u}_{\tilde{a}})\in\rho_A^*(d_0(1_A)^0)_b\subseteq({\mathcal
L}^{\tau_{\widetilde{A}}}\widetilde{A})_{i_A(b)}$ then, using
again (\ref{Defsig}) and (\ref{sigma}), we deduce that
$\sigma_{\widetilde{A}}(b,\tilde{u}_{\tilde{a}})\in({\mathcal
L}^{\tau_A}\widetilde{A})_b.$ Thus, since $\sigma_{\widetilde{A}}$
is an involution, we conclude that the restriction of
$\sigma_{\widetilde{A}}$ to the prolongation ${\mathcal
L}^{\tau_A}\widetilde{A}$ induces an isomorphism
$\widetilde{\sigma_A}:{\mathcal
L}^{\tau_A}\widetilde{A}\rightarrow\rho_A^*(d_0(1_A)^0)$ (over the
identity $Id:A\rightarrow A$) between the vector bundles
$\tau_{\widetilde{A}}^{\tau_A}:{\mathcal
L}^{\tau_A}\widetilde{A}\rightarrow A$ and
$pr_1:\rho_A^*(d_0(1_A)^0)\rightarrow A$.

On the other hand, as we know, ${\mathcal
L}^{\tau_A}\widetilde{A}$ (res\-pec\-tively,
$\rho_A^*(d_0(1_A)^0)$) is a Lie subalgebroid of ${\mathcal
L}^{\tau_{\widetilde{A}}}\widetilde{A}$ (res\-pec\-tively,
$\rho_{\widetilde{A}}^*(T\widetilde{A})$). Therefore, using that
$\sigma_{\widetilde{A}}:{\mathcal
L}^{\tau_{\widetilde{A}}}\widetilde{A}\rightarrow\rho^*_{\widetilde{A}}(T\widetilde{A})$
is a Lie algebroid isomorphism, we obtain that
$\widetilde{\sigma_A}:{\mathcal
L}^{\tau_A}\widetilde{A}\rightarrow\rho_A^*(d_0(1_A)^0)$ is also a
Lie algebroid isomorphism.

Now, denote by $\phi_0$ (respectively, $\varphi_0$) the
$1$-cocycle of the Lie algebroid
$\tau_{\widetilde{A}}^{\tau_A}:{\mathcal
L}^{\tau_A}\widetilde{A}\rightarrow A$ (res\-pec\-ti\-vely,
$pr_1:\rho_A^*(d_0(1_A)^0)\rightarrow A$) given by (\ref{T0})
(respectively, (\ref{phi0})). From (\ref{T0}), (\ref{phi0}),
(\ref{Defsig}) and (\ref{sigma}), it follows that
$\varphi_0\circ\widetilde{\sigma_A}=\phi_0.$ Consequently, using
(\ref{1atil1}), (\ref{phi0m}), Proposition \ref{prop2.2} and the
fact that the bidual Lie algebroid to the Lie affgebroid
$\tau_A^{\tau_A}:{\mathcal J}^AA\rightarrow A$ (respectively,
$pr_1:{\mathcal J}^AA\rightarrow A$) is ${\mathcal
L}^{\tau_A}\widetilde{A}$ (respectively, $\rho_A^*(d_0(1_A)^0)$),
we prove the result.
\end{proof}

\begin{definition} The map $\sigma_A:{\mathcal J}^AA\rightarrow{\mathcal
J}^AA$ is called the canonical involution associated with the Lie
affgebroid $A$.
\end{definition}

Suppose that $(x^i)$ are local coordinates on an open subset $U$
of $M$ and that $\{e_0,e_{\alpha}\}$ is a local basis of sections
of $\tau_{\widetilde{A}}:\widetilde{A}\to M$ in $U$ adapted to
$1_A$. Denote by $(x^i,y^0,y^{\alpha})$ the co\-rres\-pon\-ding
local coordinates on $\widetilde{A}$ and by $\rho_0^i$ and
$\rho_{\alpha}^i$ the components of the anchor map
$\rho_{\widetilde{A}}$. Then,
$\{\tilde{T}_0,\tilde{T}_{\alpha},\tilde{V}_0,\tilde{V}_{\alpha}\}$
is a local basis of sections of
$\tau_{\widetilde{A}}^{\tau_{\widetilde{A}}}:{\mathcal
L}^{\tau_{\widetilde{A}}}\widetilde{A}\to\widetilde{A}$, where
$$
\begin{array}{lcllcl}
\tilde{T}_0(\tilde{a})&=&(e_0(\tau_{\widetilde{A}}(\tilde{a})),\rho_0^i\displaystyle\frac{\partial}{\partial
x^i}_{|\tilde{a}}),&\;\;\tilde{T}_{\alpha}(\tilde{a})&=&(e_{\alpha}(\tau_{\widetilde{A}}(\tilde{a})),\rho_{\alpha}^i\displaystyle\frac{\partial}{\partial
x^i}_{|\tilde{a}}),\\[8pt]
\tilde{V}_0(\tilde{a})&=&(0,\displaystyle\frac{\partial}{\partial
y^0}_{|\tilde{a}}),&\;\;\tilde{V}_{\alpha}(\tilde{a})&=&(0,\displaystyle\frac{\partial}{\partial
y^{\alpha}}_{|\tilde{a}}).
\end{array}
$$

This local basis induces a system of local coordinates
$(x^i,y^0,y^{\alpha};z^0,z^{\alpha},v^0,v^{\alpha})$ on ${\mathcal
L}^{\tau_{\widetilde{A}}}\widetilde{A}\equiv\rho^*_{\widetilde{A}}(T\widetilde{A})$.
The local expression of the canonical involution
$\sigma_{\widetilde{A}}:{\mathcal
L}^{\tau_{\widetilde{A}}}\widetilde{A}\to\rho^*_{\widetilde{A}}(T\widetilde{A})$
in these coordinates is (see \cite{LMM})
$$\sigma_{\widetilde{A}}(x^i,y^0,y^{\alpha};z^0,z^{\alpha},v^0,v^{\alpha})=(x^i,z^0,z^{\alpha};y^0,y^{\alpha},v^0,v^{\alpha}+C_{0\gamma}^{\alpha}(z^0y^{\gamma}-z^{\gamma}y^0)+C_{\beta\gamma}^{\alpha}z^{\beta}y^{\gamma}).$$
Here, $C_{0\gamma}^{\alpha}$ and $C_{\beta\gamma}^{\alpha}$ are
the structure functions of the Lie bracket
$\lcf\cdot,\cdot\rcf_{\widetilde{A}}$ with respect to the basis
$\{e_0,e_{\alpha}\}$.

On the other hand, the local equations defining the affine
subbundle ${\mathcal J}^AA$ of ${\mathcal
L}^{\tau_{\widetilde{A}}}\widetilde{A}$ are $y^0=1$, $z^0=1$,
$v^0=0$. Thus, $(x^i,y^{\alpha};z^{\alpha},v^{\alpha})$ may be
considered as local coordinates on ${\mathcal J}^AA$. Using these
coordinates, we deduce that the local expression of
$\sigma_A:{\mathcal J}^AA\to{\mathcal J}^AA$ is
$$\sigma_{{A}}(x^i,y^{\alpha};z^{\alpha},v^{\alpha})=(x^i,z^{\alpha};y^{\alpha},v^{\alpha}+C_{0\gamma}^{\alpha}(y^{\gamma}-z^{\gamma})+C_{\beta\gamma}^{\alpha}z^{\beta}y^{\gamma}).$$

Finally, the local equations defining the vector subbundles
${\mathcal L}^{\tau_A}V$ and $\rho_A^*(TV)$ of ${\mathcal
L}^{\tau_{\widetilde{A}}}\widetilde{A}$ and
$\rho^*_{\widetilde{A}}(T\widetilde{A})$, respectively, are
$y^0=1$, $z^0=0$, $v^0=0$ and $z^0=1$, $y^0=0$, $v^0=0.$
Therefore, $(x^i,y^{\alpha};z^{\alpha},v^{\alpha})$ may be
considered as local coordinates on ${\mathcal L}^{\tau_A}V$ and
$\rho_A^*(TV)$. Using these coordinates, we obtain that the local
expression of $\sigma_A^l:{\mathcal L}^{\tau_A}V\to\rho_A^*(TV)$
is
\begin{equation}\label{sigmaAl}
\sigma_A^l(x^i,y^{\alpha};z^{\alpha},v^{\alpha})=(x^i,z^{\alpha};y^{\alpha},v^{\alpha}-C_{0\gamma}^{\alpha}z^{\gamma}+C_{\beta\gamma}^{\alpha}z^{\beta}y^{\gamma}).
\end{equation}

\setcounter{equation}{0}
\section{Tulczyjew's triple associated with a Lie
affgebroid and a Hamiltonian section}\label{sec5}

Let $\tau_A:A\rightarrow M$ be a Lie affgebroid modelled on the
Lie algebroid $\tau_V:V\rightarrow M$. Denote by
$(\lcf\cdot,\cdot\rcf_V,D,\rho_A)$ the Lie affgebroid structure on
$A$, by $\rho_A^*(TV)$ (respectively, $\rho_A^*(TV^*)$) the
pull-back of the vector bundle $T\tau_V:TV\rightarrow TM$
(respectively, $T\tau_V^*:TV^*\rightarrow TM$) over the anchor map
$\rho_A:A\rightarrow TM$ and by ${\mathcal
L}^{\tau_V^*}\widetilde{A}$ (respectively, ${\mathcal
L}^{\tau_A}V$ and ${\mathcal L}^{\tau_V^*}V$) the prolongation of
the Lie algebroid $\tau_{\widetilde{A}}:\widetilde{A}\to M$
(respectively, $\tau_V:V\rightarrow M$) over the projection
$\tau_V^*:V^*\rightarrow M$ (respectively, $\tau_A:A\to M$ and
$\tau_V^*:V^*\to M$).

Suppose that $(x^i)$ are local coordinates on an open subset $U$
of $M$ and that $\{e_0,e_{\alpha}\}$ is a local basis of
$\Gamma(\tau_{\widetilde{A}})$ adapted to  $1_A.$ Denote by
$(x^i,y^0,y^{\alpha})$ (respectively, $(x^i,y^{\alpha})$) the
corresponding local coordinates on $\widetilde{A}$ (respectively,
$V$ and $A$). Then, we may consider local coordinates
$(x^i,y^{\alpha};z^{\alpha},v^{\alpha})$ of $\rho_A^*(TV)$ and
${\mathcal L}^{\tau_A}V$ as in Section \ref{sec3} and the
corresponding dual coordinates
$(x^i,y^{\alpha};z_{\alpha},v_{\alpha})$ of $({\mathcal
L}^{\tau_A}V)^*$.

Now, let $(x^i,y_{\alpha};z^0,z^{\alpha},v_{\alpha})$ be local
coordinates on ${\mathcal
L}^{\tau_V^*}\widetilde{A}\equiv\rho_{\widetilde{A}}^*(TV^*)$ as
in Section \ref{sec3.1}. Then, the local equations defining the
affine subbundle $\rho_A^*(TV^*)\to V^*$ of ${\mathcal
L}^{\tau_V^*}\widetilde{A}\to V^*$ and the vector subbundle
${\mathcal L}^{\tau_V^*}V\to V^*$ of ${\mathcal
L}^{\tau_V^*}\widetilde{A}\to V^*$ are $z^0=1$ and $z^0=0$,
respectively. Thus, $(x^i,y_{\alpha};z^{\alpha},v_{\alpha})$ may
be considered as a system of local coordinates on $\rho_A^*(TV^*)$
and ${\mathcal L}^{\tau_V^*}V$.

Next, we will introduce the so-called Tulczyjew's triple
associated with the Lie affgebroid $A$ and a Hamiltonian section.

For this purpose, we will proceed in two steps.

{\it First step}: In this first step, we will introduce a
canonical isomorphism $A_A:\rho_A^*(TV^*)\rightarrow({\mathcal
L}^{\tau_A}V)^*$, over the identity of $A$, between the vector
bundles $\rho_A^*(TV^*)\rightarrow A$ and $({\mathcal
L}^{\tau_A}V)^*\rightarrow A$.

Let $\langle\cdot,\cdot\rangle:V\times_MV^*\rightarrow\R$ be the
natural pairing given by
$$\langle u,\alpha\rangle=\alpha(u),\makebox[1cm]{for}(u,\alpha)\in V_x\times
V_x^*,$$ with $x\in M$. If $b\in A$, $(b,X_u)\in\rho_A^*(TV)_b$
and $(b,X_{\alpha})\in\rho_A^*(TV^*)_b$ then
$$(X_u,X_{\alpha})\in
T_{(u,\alpha)}(V\times_MV^*)=\{(X'_u,X'_{\alpha})\in T_uV\times
T_{\alpha}V^*/(T_u\tau_V)(X'_u)=(T_{\alpha}\tau^*_V)(X'_{\alpha})\}$$
and we may consider the map
$\widetilde{T\langle\cdot,\cdot\rangle}:\rho_A^*(TV)\times_A\rho_A^*(TV^*)\rightarrow\R$
defined by
$$\widetilde{T\langle\cdot,\cdot\rangle}((b,X_u),(b,X_{\alpha}))=dt_{\langle u,\alpha\rangle}((T_{(u,\alpha)}\langle\cdot,\cdot\rangle)(X_u,X_{\alpha})),$$
where $t$ is the usual coordinate on $\R$. The local expression of
the map $\widetilde{T\langle\cdot,\cdot\rangle}$ is
$$\widetilde{T\langle\cdot,\cdot\rangle}((x^i,y^{\alpha};z^{\alpha},v^{\alpha}),(x^i,y_{\alpha};z^{\alpha},v_{\alpha}))=y^{\alpha}v_{\alpha}+v^{\alpha}y_{\alpha}.$$
Thus, $\widetilde{T\langle\cdot,\cdot\rangle}$ is also a
non-singular pairing and it induces an isomorphism (over the
identity of $A$) between the vector bundles
$\rho_A^*(TV)\rightarrow A$ and $\rho_A^*(TV^*)^*\rightarrow A$
which we also denote by $\widetilde{T\langle\cdot,\cdot\rangle}$,
that is,
$\widetilde{T\langle\cdot,\cdot\rangle}:\rho_A^*(TV)\rightarrow\rho_A^*(TV^*)^*.$
Note that
\begin{equation}\label{Tpartil}
\widetilde{T\langle\cdot,\cdot\rangle}(x^i,y^{\alpha};z^{\alpha},v^{\alpha})=(x^i,v^{\alpha};z^{\alpha},y^{\alpha}).
\end{equation}

Next, we consider the isomorphism of vector bundles
$A_A^*:{\mathcal L}^{\tau_A}V\rightarrow\rho_A^*(TV^*)^*$ given by
\begin{equation}\label{Aas}
A_A^*=\widetilde{T\langle\cdot,\cdot\rangle}\circ\sigma_A^l,
\end{equation}
$\sigma_A^l:{\mathcal L}^{\tau_A}V\rightarrow\rho_A^*(TV)$ being
the linear part of the canonical involution associated with the
Lie affgebroid $A$ (see Section \ref{sec3}). Then, the isomorphism
$A_A:\rho_A^*(TV^*)\rightarrow({\mathcal L}^{\tau_A}V)^*$ is just
the dual map to $A_A^*:{\mathcal
L}^{\tau_A}V\rightarrow\rho_A^*(TV^*)^*$.

From (\ref{sigmaAl}), (\ref{Tpartil}) and (\ref{Aas}), it follows
that
\[
A_A^*(x^i,y^{\alpha};z^{\alpha},v^{\alpha})=(x^i,v^{\gamma}-
C_{0\beta}^{\gamma}z^{\beta}+C_{\alpha\beta}^{\gamma}z^{\alpha}y^{\beta};
y^{\gamma},z^{\gamma}),
\]
and therefore
\begin{equation}\label{AA}
A_A(x^i,y_{\alpha};z^{\alpha},v_{\alpha})=
(x^i,z^{\gamma};v_{\gamma}-C_{0\gamma}^{\beta}y_{\beta}
-C_{\alpha\gamma}^{\beta}z^{\alpha}y_{\beta},y_{\gamma}).
\end{equation}

 {\it Second step}:
 Let $\tau_{A^+}:A^+\to M$ be the dual vector bundle to the affine
 bundle $\tau_A:A\to M$, $\mu:A^+\rightarrow V^*$ be
the canonical projection and $h:V^*\to A^+$ be a Hamiltonian
section, that is, $h$ is a section of $\mu:A^+\rightarrow V^*$.

Denote by $(\Omega_h,\eta)$ the cosymplectic structure on
${\mathcal L}^{\tau_V^*}\widetilde{A}$ given by (\ref{Omegah}) and
(\ref{eta}). Then, a direct computation, using (\ref{eta}), proves
that
$$\eta^{-1}\{1\}=\rho_A^*(TV^*)\makebox[1cm]{and}\eta^{-1}\{0\}={\mathcal
L}^{\tau_V^*}V.$$ Thus, $\rho_A^*(TV^*)$ is an affine bundle over
$V^*$ with affine bundle projection
$\widetilde{\pi_{V^*}}:\rho_A^*(TV^*)\rightarrow V^*$ defined by
$\widetilde{\pi_{V^*}}(a,X)=\pi_{V^*}(X)$ and, furthermore, this
affine bundle admits a Lie affgebroid structure in such a way that
the bidual Lie algebroid is
$\tau_{\widetilde{A}}^{\tau_V^*}:{\mathcal
L}^{\tau_V^*}\widetilde{A}\rightarrow V^*$. In addition, the Lie
affgebroid $\widetilde{\pi_{V^*}}:\rho_A^*(TV^*)\rightarrow V^*$
is modelled on the Lie algebroid $\tau_V^{\tau_V^*}:{\mathcal
L}^{\tau_V^*}V\rightarrow V^*$.

In this step, we will introduce an affine isomorphism
$\flat_{\Omega_h}:\rho_A^*(TV^*)\rightarrow({\mathcal
L}^{\tau_V^*}V)^*$, over the identity of $V^*$, between the affine
bundle $\widetilde{\pi_{V^*}}:\rho_A^*(TV^*)\rightarrow V^*$ and
the vector bundle $(\tau_V^{\tau_V^*})^*:({\mathcal
L}^{\tau_V^*}V)^*\rightarrow V^*$.

The map $\flat_{\Omega_h}$ is defined as follows. If $\alpha\in
V^*$ and $(a,X_{\alpha})\in\rho_A^*(TV^*)_{\alpha}$ then
\begin{equation}\label{b0mh}
\{\flat_{\Omega_h}(\alpha)(a,X_{\alpha})\}(u,Y_{\alpha})=\Omega_h(\alpha)((i_A(a),X_{\alpha}),(i_V(u),Y_{\alpha})),
\end{equation}
for $(u,Y_{\alpha})\in({\mathcal L}^{\tau_V^*}V)_{\alpha}$.

On the other hand, let $\flat_{\Omega_V}:{\mathcal
L}^{\tau_V^*}V\rightarrow({\mathcal L}^{\tau_V^*}V)^*$ be the
canonical isomorphism over the identity of $V^*$ induced by the
canonical symplectic section $\Omega_V$ associated with the Lie
algebroid $\tau_V:V\to M$, that is,
\begin{equation}\label{b0mV}
\{\flat_{\Omega_V}(\alpha)(u,Y_{\alpha})\}(v,Z_{\alpha})=\Omega_V(\alpha)((u,Y_{\alpha}),(v,Z_{\alpha})),
\end{equation}
for $(u,Y_{\alpha}),(v,Z_{\alpha})\in({\mathcal
L}^{\tau_V^*}V)_{\alpha}$, with $\alpha\in V^*$. Then, using
(\ref{inclusion2}), (\ref{b0mh}) and (\ref{b0mV}), it follows that
$\flat_{\Omega_h}$ is an affine isomorphism over the identity of
$V^*$ and the corresponding linear isomorphism between the vector
bundles $\tau_V^{\tau_V^*}:{\mathcal L}^{\tau_V^*}V\rightarrow
V^*$ and $(\tau_V^{\tau_V^*})^*:({\mathcal
L}^{\tau_V^*}V)^*\rightarrow V^*$ is just the map
$\flat_{\Omega_V}$.

In conclusion, we have the following commutative diagram

\begin{picture}(455,90)(120,20)
\put(245,20){\makebox(0,0){$A$}} \put(310,70){\vector(-3,-2){60}}
 \put(480,70){\vector(-3,-2){60}}
\put(417,20){\makebox(0,0){$V^*$}}
\put(170,40){$(\tau_V^{\tau_A})^*$}
\put(360,40){$\widetilde{\pi_{V^*}}$}
\put(460,40){$(\tau_V^{\tau_V^*})^*$}\put(180,70){\vector(3,-2){60}}
\put(290,40){$pr_1$} \put(350,70){\vector(3,-2){60}}
\put(180,80){\makebox(0,0){$({\mathcal L}^{\tau_A}V)^* $}}
\put(250,85){$A_A$}\put(400,85){$\flat_{\Omega_h}$}\put(300,80){\vector(-1,0){95}}
\put(333,80){\makebox(0,0){$\rho_A^*(TV^*)$}}
\put(357,80){\vector(1,0){95}}
\put(480,80){\makebox(0,0){$({\mathcal
L}^{\tau_V^*}V)^*$}}\end{picture}

\vspace{0.7cm} This diagram will be called {\it Tulczyjew's triple
associated with the Lie affgebroid $A$ and the Hamiltonian section
$h:V^*\to A^+$}.

\setcounter{equation}{0}
\section{The prolongation of a symplectic Lie affgebroid}\label{sec6}

Let $\tau_A:A\rightarrow M$ be a Lie affgebroid modelled on the
Lie algebroid $\tau_V:V\rightarrow M.$ Denote by
$(\lcf\cdot,\cdot\rcf_{\widetilde{A}},\rho_{\widetilde{A}})$ the
Lie algebroid structure on the bidual bundle
$\tau_{\widetilde{A}}:\widetilde{A}\to M$ to $A$.

If $f\in C^{\infty}(M)$ one may define {\it the complete and
vertical lift $f^c$ and $f^v$ of $f$ to $A$} as the function on
$A$ given by
\begin{equation}\label{6.0}
f^c(a)=\rho_A(a)(f),\;\;\;f^v(a)=f(\tau_A(a)),\;\;\; \mbox{ for
$a\in A$}. \end{equation} A direct computation proves that
\begin{equation}\label{cvf}
(fg)^c=f^cg^v+f^vg^c,\;\;\;(fg)^v=f^vg^v,\makebox[1cm]{for}f,g\in
C^{\infty}(M).
\end{equation}

Now, if $\tilde{X}$ is a section of
$\tau_{\widetilde{A}}:\widetilde{A}\rightarrow M$, we may consider
the vertical and complete lift $\tilde{X}^v$ and $\tilde{X}^c$ of
$\tilde{X}$ as the vector fields on $\widetilde{A}$ defined in
Section \ref{sec1.1.2}. Using these constructions, we may
introduce the sections $\tilde{X}^{\textbf{v}}$ and
$\tilde{X}^{\textbf{c}}$ of the prolongation
$\tau_{\widetilde{A}}^{\tau_{\widetilde{A}}}:{\mathcal
L}^{\tau_{\widetilde{A}}}\widetilde{A}\to\widetilde{A}$ given by
\begin{equation}\label{lcvs}
\tilde{X}^{\textbf{v}}(\tilde{a})=(0(\tau_{\widetilde{A}}(\tilde{a})),\tilde{X}^v(\tilde{a})),
\;\;\tilde{X}^{\textbf{c}}(\tilde{a})=(\tilde{X}(\tau_{\widetilde{A}}(\tilde{a})),\tilde{X}^c(\tilde{a})),
\end{equation}
for $\tilde{a}\in\widetilde{A}$. If $\{e_0,e_{\alpha}\}$ is a
local basis of $\Gamma(\tau_{\widetilde{A}})$ adapted to
$1_A:\widetilde{A}\to \R$, then
$\{e_0^{\textbf{c}},e_{\alpha}^{\textbf{c}},e_0^{\textbf{v}},e_{\alpha}^{\textbf{v}}\}$
is a local basis of
$\Gamma(\tau_{\widetilde{A}}^{\tau_{\widetilde{A}}})$ (see
\cite{LMM,M1}).

Next, denote by $i_A:A\to\widetilde{A}$, $i_V:V\to\widetilde{A}$
the canonical inclusions  and by $\tau_V^{\tau_A}:{\mathcal
L}^{\tau_A}V\to A$ the prolongation of the Lie algebroid
$\tau_V:V\to M$ over the projection $\tau_A:A\to M$.

If $X$ is a section of $\tau_V:V\to M$ then $i_V\circ X\in
\Gamma(\tau_{\widetilde{A}})$ and we have that the restriction to
$A$ of the vector fields $(i_V\circ X)^v$ and $(i_V\circ X)^c$ are
tangent to $A$. Thus,
\begin{equation}\label{lvcn}
X^{\textbf{v}}=(i_V\circ X)^{\textbf
v}_{|A}\in\Gamma(\tau_V^{\tau_A}),\;\;X^{\textbf{c}}=(i_V\circ
X)^{\textbf c}_{|A}\in\Gamma(\tau_V^{\tau_A}).
\end{equation}

Moreover, using the properties of the complete and vertical lifts
(see \cite{LMM,M1}), we deduce that
\begin{equation}\label{cvl2}
\begin{array}{c}
(fX)^{\textbf{c}}=f^cX^{\textbf{v}}+f^vX^{\textbf{c}},\;\;(fX)^{\textbf{v}}=f^vX^{\textbf{v}},\\
\lcf X^{\textbf{c}},Y^{\textbf{c}}\rcf_V^{\tau_A}=\lcf
X,Y\rcf_V^{\textbf{c}},\;\;\lcf
X^{\textbf{c}},Y^{\textbf{v}}\rcf_V^{\tau_A}=\lcf
X,Y\rcf_V^{\textbf{v}},\;\;\lcf
X^{\textbf{v}},Y^{\textbf{v}}\rcf_V^{\tau_A}=0,\\
\rho_V^{\tau_A}(X^{\textbf{c}})=(i_V\circ
X)^c_{|A},\;\;\;\rho_V^{\tau_A}(X^{\textbf{v}})=(i_V\circ
X)^v_{|A},
\end{array}
\end{equation}
for $f\in C^{\infty}(M)$ and $X,Y\in\Gamma(\tau_V)$, where
$(\lcf\cdot,\cdot\rcf_V,\rho_V)$ and
$(\lcf\cdot,\cdot\rcf_V^{\tau_A},\rho_V^{\tau_A})$ are the Lie
algebroids structures on $V$ and ${\mathcal L}^{\tau_A}V$,
respectively. In addition, if $\{e_{\alpha}\}$ is a local basis of
$\Gamma(\tau_V)$ then
$\{e_{\alpha}^{\textbf{c}},e_{\alpha}^{\textbf{v}}\}$ is a local
basis of $\Gamma(\tau_A^{\tau_V})$ (for more details, see
\cite{MMeS}).

On the other hand, if  $pr_1:V\times TA\to V$ is the canonical
projection on the first factor then the pair $(pr_{1|{\mathcal
L}^{\tau_A}V},\tau_A)$ is a morphism between the Lie algebroids
$({\mathcal
L}^{\tau_A}V,\lcf\cdot,\cdot\rcf_V^{\tau_A},\rho_V^{\tau_A})$ and
$(V,\lcf\cdot,\cdot\rcf_V,\rho_V)$. Thus, if
$\alpha\in\Gamma(\wedge^k \tau _V^*)$ we may consider the section
$\alpha^{\textbf{v}}$ of the vector bundle $\wedge^k({\mathcal
L}^{\tau_A}V)^*\to A$ defined by
\begin{equation}\label{lvf}
\alpha^{\textbf{v}}=(pr_{1|{\mathcal
L}^{\tau_A}V},\tau_A)^*(\alpha).
\end{equation}

$\alpha^{\textbf{v}}$ is called {\it the vertical lift to
${\mathcal L}^{\tau_A}V$ of $\alpha$} and it is clear that
\begin{equation}\label{lvfd}
d^{{\mathcal
L}^{\tau_A}V}\alpha^{\textbf{v}}=(d^V\alpha)^{\textbf{v}}.
\end{equation}

Furthermore, we have that
$(\alpha_1\wedge\dots\wedge\alpha_k)^{\textbf{v}}=\alpha_1^{\textbf{v}}\wedge\dots\wedge\alpha_k^{\textbf{v}},\makebox[1cm]{for}\alpha_i\in\Gamma(\wedge^{k_i}V^*).$

Now, we define the complete lift of $\alpha$ as follows.
\begin{proposition}\label{prop6.1}
If $\alpha$ is a section of the vector bundle $\wedge^kV^*\to M$,
then there exists a unique section $\alpha^{\textbf{c}}$ of the
vector bundle $\wedge^k({\mathcal L}^{\tau_A}V)^*\to A$ such that
$$
\begin{array}{rcl}
\alpha^{\textbf{\em c}}(X_1^{\textbf{\em c}},\dots,X_k^{\textbf{\em c}})&=&\alpha(X_1,\dots,X_k)^c,\\
\alpha^{\textbf{\em c}}(X_1^{\textbf{\em v}},X_2^{\textbf{\em c}},\dots,X_k^{\textbf{\em c}})&=&\alpha(X_1,X_2,\dots,X_k)^v,\\
\alpha^{\textbf{\em c}}(X_1^{\textbf{\em
v}},\dots,X_s^{\textbf{\em v}},X_{s+1}^{\textbf{\em
c}},\dots,X_k^{\textbf{\em c}})&=&0,\makebox[1cm]{if}2\leq s\leq
k,
\end{array}
$$
for $X_1,\dots,X_k\in\Gamma(\tau_V)$. Moreover, $d^{{\mathcal
L}^{\tau_A}V}\alpha^{\textbf{\em c}}=(d^V\alpha)^{\textbf{\em
c}}.$
\end{proposition}
\begin{proof}
Using (\ref{cvf}), (\ref{cvl2}), (\ref{lvf}), (\ref{lvfd}) and
proceeding as in the proof of Proposition 6.3 in \cite{LMM}, we
deduce the result.
\end{proof}

The section $\alpha^{\textbf{c}}$ of the vector bundle
$\wedge^k({\mathcal L}^{\tau_A}V)^*\to A$ is called {\it the
complete lift of $\alpha$}.

\begin{remark}{\em $i)$ If $\{e_{\alpha}\}$ is a
basis of $\Gamma(\tau_V)$ and $\{e^{\alpha}\}$ is the dual basis
to $\{e_{\alpha}\}$ then
$$\begin{array}{c}
(e^{\alpha})^{\textbf{c}}((e_{\beta})^{\textbf{v}})=(e^{\alpha})^{\textbf{v}}((e_{\beta})^{\textbf{c}})
=\delta_{\alpha\beta},\;\;\;\;\;(e^{\alpha})^{\textbf{c}}((e_{\beta})^{\textbf{c}})
=(e^{\alpha})^{\textbf{v}}((e_{\beta})^{\textbf{v}})=0.
\end{array}
$$
Therefore,
$\{(e^{\alpha})^{\textbf{v}},(e^{\alpha})^{\textbf{c}}\}$ is the
dual basis to the local basis
$\{(e_{\alpha})^{\textbf{c}},(e_{\alpha})^{\textbf{v}}\}$ of
$\Gamma(\tau_V^{\tau_A})$.

\medskip

$ii)$ If $\alpha_i\in\Gamma(\wedge^{k_i}\tau _V^*)$,
$i=1,\dots,k$, then
\begin{equation}\label{lcf}
(\alpha_1\wedge\dots\wedge\alpha_k)^{\textbf{c}}=
\displaystyle\sum_{i=1}^k\alpha_1^{\textbf{v}}\wedge\dots\wedge\alpha_i^{\textbf{c}}\wedge\dots\wedge\alpha_k^{\textbf{v}}.
\end{equation}
\hfill$\diamondsuit$}
\end{remark}

Next, we will introduce the notion of a symplectic Lie affgebroid.

\begin{definition}
A Lie affgebroid $\tau_A:A\to M$ modelled on the Lie algebroid
$\tau_V:V\to M$ is said to be symplectic if $\tau_V:V\to M$ admits
a symplectic section $\Omega$, that is, $\Omega$ is a section of
the vector bundle $\wedge^2V^*\to M$ such that:
\begin{enumerate}
\item[(i)] For all $x\in M$, the $2$-form $\Omega(x):V_x\times
V_x\to\R$ on the vector space $V_x$ is non-degenerate and
\item[(ii)] $\Omega$ is a $2$-cocycle, i.e., $d^V\Omega=0$.
\end{enumerate}
\end{definition}

\begin{examples}{\em $(i)$ Let $\tau:M\to\R$ be a fibration. Then, as we know (see
Section \ref{sec1.2}), the $1$-jet bundle $\tau_{1,0}:J^1\tau\to
M$ is a Lie affgebroid modelled on the Lie algebroid
$(\pi_M)_{|V\tau}:V\tau\to M$. Now, suppose that $M$ has odd
dimension $2n+1$ and that $(\Omega,\eta)$ is a cosymplectic
structure on $M$, with $\eta=\tau^*(dt)$, $t$ being the usual
coordinate on $\R$. This means that $\Omega$ is a closed $2$-form
and that
$\eta\wedge\Omega^n=\eta\wedge\Omega\wedge\dots^{(n}\dots\wedge\Omega$
is a volume element on $M$. Thus, it is easy to prove that the
restriction to $V\tau$ of $\Omega$ is a symplectic section of
$(\pi_M)_{|V\tau}:V\tau\to M$ and, therefore, the Lie affgebroid
$\tau_{1,0}:J^1\tau\to M$ is symplectic.

$(ii)$ Let $\tau_A:A\to M$ be a Lie affgebroid modelled on the Lie
algebroid $\tau_V:V\to M$. Denote by $\rho_A^*(TV^*)$ the
pull-back of the vector bundle $T\tau_V^*:TV^*\to TM$ over the
anchor map $\rho_A:A\to TM$. Then, $\rho_A^*(TV^*)$ is a Lie
affgebroid over $V^*$ with affine bundle projection
$\widetilde{\pi_{V^*}}:\rho_A^*(TV^*)\to V^*$ given by
$\widetilde{\pi_{V^*}}(a,X)=\pi_{V^*}(X)$, for
$(a,X)\in\rho_A^*(TV^*)$ (see Section \ref{sec5}). Moreover, the
Lie affgebroid $\widetilde{\pi_{V^*}}:\rho_A^*(TV^*)\to V^*$ is
modelled on the Lie algebroid $\tau_V^{\tau_V^*}:{\mathcal
L}^{\tau_V^*}V\to V^*$ which admits a canonical symplectic section
$\Omega_V$ (see Section \ref{sec1.1.1}). Therefore, the Lie
affgebroid $\widetilde{\pi_{V^*}}:\rho_A^*(TV^*)\to V^*$ is
symplectic.\hfill$\triangle$ }
\end{examples}

Now, we will see that the prolongation of a symplectic Lie
affgebroid over the affine bundle projection is also a symplectic
Lie affgebroid. We recall that if $\tau_A:A\to M$ is a Lie
affgebroid modelled on the Lie algebroid $\tau_V:V\to M$ then the
prolongation
$${\mathcal J}^AA=\{(a,v)\in A\times TA/\rho_A(a)=(T\tau_A)(v)\}$$
is a Lie affgebroid modelled on the Lie algebroid
$\tau_V^{\tau_A}:{\mathcal L}^{\tau_A}V\to A$ (see Section
\ref{sec3}). Furthermore, we have that
\begin{theorem}\label{th6.3} Let $\tau_A:A\to M$ be a symplectic Lie affgebroid
modelled on the Lie algebroid $\tau_V:V\to M$ and $\Omega$ be a
symplectic section of $\tau_V:V\to M$. Then, the prolongation
${\mathcal J}^AA$ of the Lie affgebroid $A$ over the projection
$\tau_A:A\to M$ is a symplectic Lie affgebroid and the complete
lift $\Omega^{\textbf{c}}$ of $\Omega$ to the prolongation
${\mathcal L}^{\tau_A}V$ is a symplectic section of
$\tau_V^{\tau_A}:{\mathcal L}^{\tau_A}V\to A$.
\end{theorem}
\begin{proof}
Using (\ref{lcf}) and proceeding as in the proof of the Theorem
6.5 in \cite{LMM}, we deduce the result.
\end{proof}

\begin{example}\label{ex6.6}{\em
Let $\tau_A:A\to M$ be a Lie affgebroid modelled on the Lie
algebroid $\tau_V:V\to M$. Then, as we know, the pull-back
$\widetilde{\pi_{V^*}}:\rho_A^*(TV^*)\to V^*$ of the vector bundle
$T\tau^*_V:TV^*\to TM$ over the anchor map $\rho_A:A\to TM$ is a
Lie affgebroid modelled on the symplectic Lie algebroid
$\tau_V^{\tau_V^*}:{\mathcal L}^{\tau_V^*}V\to V^*$.

Now, suppose that $(x^i)$ are local coordinates on $M$ and that
$\{e_{\alpha}\}$ is a local basis of $\Gamma(\tau_V)$. Then, we
may consider the corresponding local coordinates
$(x^i,y_{\alpha})$ of $V^*$ and the corresponding local basis
$\{\tilde{e}_{\alpha},\bar{e}_{\alpha}\}$ of
$\Gamma(\tau_V^{\tau_V^*})$ (see Section \ref{sec1.1.1}). This
local basis induces a system of local coordinates
$(x^i,y_{\alpha};z^{\alpha},v_{\alpha})$ on ${\mathcal
L}^{\tau_V^*}V$. Moreover, if ${\mathcal
L}^{\widetilde{\pi_{V^*}}}({\mathcal L}^{\tau_V^*}V)$ is the
prolongation of ${\mathcal L}^{\tau_V^*}V$ over the projection
$\widetilde{\pi_{V^*}}:\rho_A^*(TV^*)\to V^*$ then
$\{\tilde{e}_{\alpha}^{\textbf{c}},\bar{e}_{\alpha}^{\textbf{c}},\tilde{e}_{\alpha}^{\textbf{v}},\bar{e}_{\alpha}^{\textbf{v}}\}$
is a local basis of
$\Gamma((\tau_V^{\tau_V^*})^{\widetilde{\pi_{V^*}}})$ and if
$\{\tilde{e}^{\alpha},\bar{e}^{\alpha}\}$ is the dual basis to
$\{\tilde{e}_{\alpha},\bar{e}_{\alpha}\}$ then
$$\{(\tilde{e}^{\alpha})^{\textbf{v}},(\bar{e}^{\alpha})^{\textbf{v}},(\tilde{e}^{\alpha})^{\textbf{c}},(\bar{e}^{\alpha})^{\textbf{c}}\}$$
is the dual basis of
$\{\tilde{e}_{\alpha}^{\textbf{c}},\bar{e}_{\alpha}^{\textbf{c}},\tilde{e}_{\alpha}^{\textbf{v}},\bar{e}_{\alpha}^{\textbf{v}}\}$.
Note that if $\rho_{\alpha}^i$, $C_{0\alpha}^{\gamma}$ and
$C_{\alpha\beta}^{\gamma}$ are the structure functions of the Lie
algebroid
$(\widetilde{A},\lcf\cdot,\cdot\rcf_{\widetilde{A}},\rho_{\widetilde{A}})$
with respect to the coordinates $(x^i)$ and to the basis
$\{e_0,e_{\alpha}\}$ then, from (\ref{coretilde}), (\ref{lcvs})
and (\ref{lvcn}), we deduce that
\begin{equation}\label{bcv}
\begin{array}{rcl}
\tilde{e}_{\alpha}^{\textbf{c}}(x^i,y_{\gamma};z^{\gamma},v_{\gamma})&=&(\tilde{e}_{\alpha}(x^j,y_{\beta}),\rho_{\alpha}^i(x^j)\displaystyle\frac{\partial}{\partial
x^i}-(C_{0\alpha}^\gamma(x^j)-C_{\alpha\beta}^{\gamma}(x^j)z^{\beta})\frac{\partial}{\partial
z^{\gamma}}),\\[8pt]
\bar{e}_{\alpha}^{\textbf{c}}(x^i,y_{\gamma};z^{\gamma},v_{\gamma})&=&(\bar{e}_{\alpha}(x^j,y_{\beta}),\displaystyle\frac{\partial}{\partial
y_{\alpha}}),\\[8pt]
\tilde{e}_{\alpha}^{\textbf{v}}(x^i,y_{\gamma};z^{\gamma},v_{\gamma})&=&(0,\displaystyle\frac{\partial}{\partial
z^{\alpha}}),\;\;\;\;\bar{e}_{\alpha}^{\textbf{v}}(x^i,y_{\gamma};z^{\gamma},v_{\gamma})=(0,
\frac{\partial}{\partial v_{\alpha}}).
\end{array}\end{equation}

In addition, using (\ref{formas}) and (\ref{lcf}), we deduce that
the local expression of the complete lifts
$\lambda_V^{\textbf{c}}$ and $\Omega_V^{\textbf{c}}$ are
\begin{equation}\label{OmegaVC}
\begin{array}{lcl}
\lambda_V^{\textbf{c}}(x^i,y_{\alpha};z^{\alpha},v_{\alpha})&=&v_{\alpha}(\tilde{e}^{\alpha})^{\textbf{v}}+y_{\alpha}(\tilde{e}^{\alpha})^{\textbf{c}},\\[8pt]
\Omega_V^{\textbf{c}}(x^i,y_{\alpha};z^{\alpha},v_{\alpha})&=&(\tilde{e}^{\alpha})^{\textbf{c}}\wedge(\bar{e}^{\alpha})^{\textbf{v}}+(\tilde{e}^{\alpha})^{\textbf{v}}\wedge(\bar{e}^{\alpha})^{\textbf{c}}+C_{\alpha\beta}^{\gamma}y_{\gamma}(\tilde{e}^{\alpha})^{\textbf{c}}\wedge(\tilde{e}^{\beta})^{\textbf{v}}\\[8pt]
&+&\displaystyle\frac{1}{2}((\rho_0^i +
\rho^i_{\mu}z^\mu)\frac{\partial
C_{\alpha\beta}^{\gamma}}{\partial
x^i}y_{\gamma}+C_{\alpha\beta}^{\gamma}v_{\gamma})(\tilde{e}^{\alpha})^{\textbf{v}}\wedge(\tilde{e}^{\beta})^{\textbf{v}}.
\end{array}
\end{equation}

\hfill$\triangle$}
\end{example}

\setcounter{equation}{0}
\section{Lagrangian Lie subaffgebroids in symplectic Lie
affgebroids}\label{sec7}

First of all, we will introduce the notion of a Lie subaffgebroid
of a Lie affgebroid.

\begin{definition}
Let $\tau_A:A\to M$ be a Lie affgebroid modelled on the Lie
algebroid $\tau_V:V\to M$. A Lie subaffgebroid of $A$ is a Lie
affgebroid morphism $((j:A'\to A,i:M'\to M),(j^l:V'\to V,i:M'\to
M))$ such that $i$ is an injective inmersion and $j:A'\to A$ is
also injective.
\end{definition}

\begin{examples}\label{ex7.2}{\em Let $\tau_A:A\to M$ be a Lie affgebroid modelled on the Lie
algebroid $\tau_V:V\to M.$

$i)$ Denote by
$(\lcf\cdot,\cdot\rcf_{\widetilde{A}},\rho_{\widetilde{A}})$ the
Lie algebroid structure on the bidual bundle
$\tau_{\widetilde{A}}:\widetilde{A}\to M$ to $A$. Now, suppose
that $X$ is a section of $\tau_A:A\to M$ and consider the map
$(Id,TX\circ\rho_{\widetilde{A}}):\widetilde{A}\to{\mathcal
L}^{\tau_A}\widetilde{A}$ given by
$$(Id,TX\circ\rho_{\widetilde{A}})(\tilde{a})=(\tilde{a},(TX)(\rho_{\widetilde{A}}(\tilde{a}))),\makebox[1cm]{for}\tilde{a}\in\widetilde{A}.$$
Using the definition of the anchor map of the Lie algebroid
$\tau^{\tau_A}_{\widetilde{A}}:{\mathcal
L}^{\tau_A}\widetilde{A}\to A$ and the fact that
$\rho_{\widetilde{A}}:\Gamma(\tau_{\widetilde{A}})\to\mathfrak
X(M)$ is a Lie algebra morphism, we deduce that the pair
$((Id,TX\circ\rho_{\widetilde{A}}),X)$ is a Lie algebroid
morphism. Moreover, it follows that
$((Id,TX\circ\rho_{\widetilde{A}}),X)^*(\phi_0)=1_A,$ where
$\phi_0$ is the $1$-cocycle on
$\tau^{\tau_A}_{\widetilde{A}}:{\mathcal
L}^{\tau_A}\widetilde{A}\to A$ defined by (\ref{T0}). Thus, the
pair $((Id,TX\circ\rho_{\widetilde{A}}),X)$ defines a Lie
subaffgebroid of $(\tau_A^{\tau_A}:{\mathcal J}^AA\to
A,\tau_V^{\tau_A}:{\mathcal L}^{\tau_A}V\to A)$

\begin{picture}(375,60)(130,40)
\put(200,20){\makebox(0,0){$M$}} \put(250,25){$X$}
\put(215,20){\vector(1,0){80}} \put(310,20){\makebox(0,0){$A$}}
\put(185,50){$\tau_{A}$} \put(200,70){\vector(0,-1){40}}
\put(320,50){$\tau_{A}^{\tau_A}$} \put(310,70){\vector(0,-1){40}}
\put(200,80){\makebox(0,0){$A$}}
\put(220,85){$(Id,TX\circ\rho_A)$} \put(215,80){\vector(1,0){80}}
\put(310,80){\makebox(0,0){${\mathcal J}^AA$}}
\end{picture}

\begin{picture}(375,60)(-50,-20)
\put(200,20){\makebox(0,0){$M$}} \put(250,25){$X$}
\put(210,20){\vector(1,0){80}} \put(310,20){\makebox(0,0){$A$}}
\put(185,50){$\tau_{V}$} \put(200,70){\vector(0,-1){40}}
\put(320,50){$\tau_{V}^{\tau_A}$} \put(310,70){\vector(0,-1){40}}
\put(200,80){\makebox(0,0){$V$}}
\put(220,85){$(Id,TX\circ\rho_V)$} \put(210,80){\vector(1,0){80}}
\put(310,80){\makebox(0,0){${\mathcal L}^{\tau_A}V$}}
\end{picture}

\vspace{-1cm}

$ii)$ Denote by $\rho_A^*(TV^*)$ the pull-back of the vector
bundle $T\tau_V^*:TV^*\to TM$ over the anchor map $\rho_A:A\to
TM$. $\rho_A^*(TV^*)$ is a vector bundle over $A$ with vector
bundle projection $pr_{1|\rho_A^*(TV^*)}:\rho_A^*(TV^*)\to A$
defined by
$pr_{1|\rho_A^*(TV^*)}(a,Y)=a,\makebox[1cm]{for}(a,Y)\in\rho_A^*(TV^*).$
Now, suppose that $\tilde{X}:A\to\rho_A^*(TV^*)$ is a section of
$pr_{1|\rho_A^*(TV^*)}:\rho_A^*(TV^*)\to A$ and let ${\mathcal
L}^{\tau_V^*}\widetilde{A}$ be the prolongation of the Lie
algebroid $\tau_{\widetilde{A}}:\widetilde{A}\to M$ over the map
$\tau_V^*:V^*\to M$ and $\widetilde{\pi_{V^*}}:\rho_A^*(TV^*)\to
V^*$ be the canonical projection. Then, we may consider the map
$(Id,T(\widetilde{\pi_{V^*}}\circ\tilde{X})):{\mathcal
L}^{\tau_A}\widetilde{A}\to{\mathcal L}^{\tau_V^*}\widetilde{A}$
given by
$$(Id,T(\widetilde{\pi_{V^*}}\circ\tilde{X}))(\tilde{a},Y)=(\tilde{a},T(\widetilde{\pi_{V^*}}\circ\tilde{X})(Y)),\makebox[1cm]{for}(\tilde{a},Y)\in{\mathcal
L}^{\tau_A}\widetilde{A}.$$ It is easy to prove that the pair
$((Id,T(\widetilde{\pi_{V^*}}\circ\tilde{X})),\widetilde{\pi_{V^*}}\circ\tilde{X})$

\begin{picture}(375,60)(60,40)
\put(200,20){\makebox(0,0){$A$}}
\put(240,25){$\widetilde{\pi_{V^*}}\circ\tilde{X}$}
\put(215,20){\vector(1,0){90}} \put(320,20){\makebox(0,0){$V^*$}}
\put(180,50){$\tau_{\widetilde{A}}^{\tau_A}$}
\put(200,70){\vector(0,-1){40}}
\put(330,50){$\tau_{\widetilde{A}}^{\tau_V^*}$}
\put(320,70){\vector(0,-1){40}}
\put(195,80){\makebox(0,0){${\mathcal L}^{\tau_A}\widetilde{A} $}}
\put(220,85){$(Id,T(\widetilde{\pi_{V^*}}\circ\tilde{X}))$}
\put(215,80){\vector(1,0){90}}
\put(325,80){\makebox(0,0){${\mathcal
L}^{\tau_V^*}\widetilde{A}$}}
\end{picture}
\vspace{1cm}

is a Lie algebroid morphism.

Next, denote by ${\mathcal L}^{\widetilde{\pi_{V^*}}}({\mathcal
L}^{\tau_V^*}\widetilde{A})$ the prolongation of the Lie algebroid
$\tau_{\widetilde{A}}^{\tau_V^*}:{\mathcal
L}^{\tau_V^*}\widetilde{A}\to V^*$ over the projection
$\widetilde{\pi_{V^*}}:\rho_A^*(TV^*)\to V^*$ and by
$((Id,T(\widetilde{\pi_{V^*}}\circ\tilde{X})),T\tilde{X}):{\mathcal
L}^{\tau_A}\widetilde{A}\to{\mathcal
L}^{\widetilde{\pi_{V^*}}}({\mathcal L}^{\tau_V^*}\widetilde{A})$
the map given by
$$((Id,T(\widetilde{\pi_{V^*}}\circ\tilde{X})),T\tilde{X})(\tilde{a},Y)=((\tilde{a},T(\widetilde{\pi_{V^*}}\circ\tilde{X}))(Y),(T\tilde{X})(Y)),$$
for $(\tilde{a},Y)\in{\mathcal L}^{\tau_A}\widetilde{A}$. Then,
using the above facts, we conclude that

\begin{picture}(375,60)(100,40)
\put(160,20){\makebox(0,0){$A$}} \put(235,25){$\tilde{X}$}
\put(175,20){\vector(1,0){120}}
\put(325,20){\makebox(0,0){$\rho_A^*(TV^*)$}}
\put(140,50){$\tau_{\widetilde{A}}^{\tau_A}$}
\put(160,70){\vector(0,-1){40}}
\put(320,50){$(\tau_{\widetilde{A}}^{\tau_V^*})^{\widetilde{\pi_{V^*}}}$}
\put(315,70){\vector(0,-1){40}}
\put(155,80){\makebox(0,0){${\mathcal L}^{\tau_A}\widetilde{A} $}}
\put(180,85){$((Id,T(\widetilde{\pi_{V^*}}\circ\tilde{X})),T\tilde{X})$}
\put(175,80){\vector(1,0){120}}
\put(390,80){\makebox(0,0){${\mathcal
L}^{\widetilde{\pi_{V^*}}}({\mathcal
L}^{\tau_V^*}\widetilde{A})\subseteq{\mathcal
L}^{\tau_V^*}\widetilde{A}\times T(\rho_A^*(TV^*))$}}
\end{picture}
\vspace{1cm}

is a Lie algebroid morphism.

On the other hand, as we know (see Section \ref{sec5}), the affine
bundle $\widetilde{\pi_{V^*}}:\rho_A^*(TV^*)\to V^*$  is a Lie
affgebroid and the bidual bundle to
$\widetilde{\pi_{V^*}}:\rho_A^*(TV^*)\to V^*$ may be identified
with the Lie algebroid $\tau_{\widetilde{A}}^{\tau_V^*}:{\mathcal
L}^{\tau_V^*}\widetilde{A}\to V^*$. Therefore, one may consider
the $1$-cocycle $\tilde{\phi}_0:{\mathcal
L}^{\widetilde{\pi_{V^*}}}({\mathcal
L}^{\tau_V^*}\widetilde{A})\to\R$ of ${\mathcal
L}^{\widetilde{\pi_{V^*}}}({\mathcal
L}^{\tau_V^*}\widetilde{A})\to\rho_A^*(TV^*)$ given by
$\tilde{\phi}_0((\tilde{a},Y),Z)=1_A(\tilde{a}),$ for
$((\tilde{a},Y),Z)\in{\mathcal
L}^{\widetilde{\pi_{V^*}}}({\mathcal L}^{\tau_V^*}\widetilde{A})$
(see (\ref{T0})). Moreover, it follows that
$$(((Id,T(\widetilde{\pi_{V^*}}\circ\tilde{X})),T\tilde{X}),\tilde{X})^*(\tilde{\phi}_0)=\phi_0,$$
where $\phi_0:{\mathcal L}^{\tau_A}\widetilde{A}\to\R$ is the
$1$-cocycle of $\tau_{\widetilde{A}}^{\tau_A}:{\mathcal
L}^{\tau_A}\widetilde{A}\to A$ defined by (\ref{T0}). Note that
$$\begin{array}{lclcl}
\tilde{\phi}_0^{-1}\{1\}&=&{\mathcal
J}^{\rho_A^*(TV^*)}\rho_A^*(TV^*),\;\;\tilde{\phi}_0^{-1}\{0\}&=&{\mathcal
L}^{\widetilde{\pi_{V^*}}}({\mathcal L}^{\tau_V^*}V$$),\\
\phi_0^{-1}\{1\}&=&{\mathcal
J}^AA,\;\;\;\;\;\;\;\;\;\;\;\;\;\;\;\;\;\;\;\;\;\;\;\phi_0^{-1}\{0\}&=&{\mathcal
L}^{\tau_A}V.
\end{array}
$$

Consequently, the pair
$(((Id,T(\widetilde{\pi_{V^*}}\circ\tilde{X})),T\tilde{X}),\tilde{X})$
defines a Lie subaffgebroid of the Lie affgebroid
$((\widetilde{\pi_{V^*}})^{\widetilde{\pi_{V^*}}}:{\mathcal
J}^{\rho_A^*(TV^*)}\rho_A^*(TV^*)
\to\rho_A^*(TV^*),(\tau_V^{\tau_V^*})^{\widetilde{\pi_{V^*}}}:{\mathcal
L}^{\widetilde{\pi_{V^*}}}({\mathcal
L}^{\tau_V^*}V)\to\rho_A^*(TV^*))$

\begin{picture}(375,60)(40,40)
\put(160,20){\makebox(0,0){$A$}} \put(235,25){$\tilde{X}$}
\put(175,20){\vector(1,0){120}}
\put(325,20){\makebox(0,0){$\rho_A^*(TV^*)$}}
\put(140,50){$\tau_{A}^{\tau_A}$} \put(160,70){\vector(0,-1){40}}
\put(320,50){$(\widetilde{\pi_{V^*}})^{\widetilde{\pi_{V^*}}}$}
\put(315,70){\vector(0,-1){40}}
\put(155,80){\makebox(0,0){${\mathcal J}^{A}A $}}
\put(180,85){$((Id,T(\widetilde{\pi_{V^*}}\circ\tilde{X})),T\tilde{X})$}
\put(175,80){\vector(1,0){120}}
\put(345,80){\makebox(0,0){${\mathcal
J}^{\rho_A^*(TV^*)}\rho_A^*(TV^*)$}}
\end{picture}
\vspace{1cm}

\begin{picture}(375,60)(40,40)
\put(160,20){\makebox(0,0){$A$}} \put(235,25){$\tilde{X}$}
\put(175,20){\vector(1,0){120}}
\put(325,20){\makebox(0,0){$\rho_A^*(TV^*)$}}
\put(140,50){$\tau_{V}^{\tau_A}$} \put(160,70){\vector(0,-1){40}}
\put(320,50){$(\tau_V^{\tau_V^*})^{\widetilde{\pi_{V^*}}}$}
\put(315,70){\vector(0,-1){40}}
\put(155,80){\makebox(0,0){${\mathcal L}^{\tau_A}V $}}
\put(180,85){$((Id,T(\widetilde{\pi_{V^*}}\circ\tilde{X})),T\tilde{X})$}
\put(175,80){\vector(1,0){120}}
\put(335,80){\makebox(0,0){${\mathcal
L}^{\widetilde{\pi_{V^*}}}({\mathcal L}^{\tau_V^*}V)$}}
\end{picture}
\vspace{1cm}

\hfill$\triangle$}
\end{examples}
Next, we will introduce the notion of a Lagrangian Lie
subaffgebroid of a symplectic Lie affgebroid.

\begin{definition}
Let $\tau_A:A\to M$ be a symplectic Lie affgebroid modelled on the
Lie algebroid $\tau_V:V\to M$ with symplectic section $\Omega$ and

\begin{picture}(375,60)(150,40)
\put(200,20){\makebox(0,0){$M'$}} \put(250,25){$i$}
\put(215,20){\vector(1,0){80}} \put(310,20){\makebox(0,0){$M$}}
\put(185,50){$\tau_{A'}$} \put(200,70){\vector(0,-1){40}}
\put(320,50){$\tau_{A}$} \put(310,70){\vector(0,-1){40}}
\put(200,80){\makebox(0,0){$A'$}} \put(250,85){$j$}
\put(215,80){\vector(1,0){80}} \put(310,80){\makebox(0,0){$A$}}
\end{picture}

\begin{picture}(375,60)(-40,-20)
\put(200,20){\makebox(0,0){$M'$}} \put(250,25){$i$}
\put(210,20){\vector(1,0){80}} \put(310,20){\makebox(0,0){$M$}}
\put(185,50){$\tau_{V'}$} \put(200,70){\vector(0,-1){40}}
\put(320,50){$\tau_{V}$} \put(310,70){\vector(0,-1){40}}
\put(200,80){\makebox(0,0){$V'$}} \put(250,85){$j^l$}
\put(210,80){\vector(1,0){80}} \put(310,80){\makebox(0,0){$V$}}
\end{picture}

\vspace{-1cm}

be a Lie subaffgebroid. Then, the Lie subaffgebroid is said to be
Lagrangian if $j^l(V'_{x'})$ is a Lagrangian subspace of the
symplectic vector space $(V_{i(x')},\Omega_{i(x')})$, for all
$x'\in M'$. In other words, we have that: $(i)$
$rank\;V'=\displaystyle\frac{1}{2}rank\;V$ and $(ii)$
$(\Omega_{i(x')})_{|j^l(V'_{x'})\times j^l(V'_{x'})}=0$, for all
$x'\in M'$.
\end{definition}
Now, let $\tau_A:A\to M$ be a symplectic Lie affgebroid modelled
on the Lie algebroid $\tau_V:V\to M$ with symplectic section
$\Omega$. Then, as we know (see Theorem \ref{th6.3}), the Lie
affgebroid $(\tau_A^{\tau_A}:{\mathcal J}^AA\to A$,
$\tau_V^{\tau_A}:{\mathcal L}^{\tau_A}V\to A)$ is symplectic and
the complete lift $\Omega^{\textbf{c}}$ of $\Omega$ is a
symplectic section of $\tau_V^{\tau_A}:{\mathcal L}^{\tau_A}V\to
A$.

Denote by $(\lcf\cdot,\cdot\rcf_V,D,\rho_A)$ the Lie affgebroid
structure of $A$ and suppose that $X:M\to A$ is a section of
$\tau_A:A\to M$. The section $X$ allows us to define the Lie
subaffgebroid  of $(\tau_A^{\tau_A}:{\mathcal J}^AA\to A,
\tau_V^{\tau_A}:{\mathcal L}^{\tau_A} V\to A)$ considered in
Examples \ref{ex7.2} $i)$ and, in addition, we will obtain a
necessary and sufficient condition for such a Lie subaffgebroid to
be Lagrangian.

For this purpose, we will introduce the operator
$D_X:\Gamma(\wedge^k\tau _V^*)\to\Gamma(\wedge^k\tau _V^*)$
defined as follows. If $\alpha\in\Gamma(\wedge^k\tau _V^*)$ then
\begin{equation}\label{Dalpha}
(D_X\alpha)(X_1,\dots,X_k)=\rho_A(X)(\alpha(X_1,\dots,X_k))-\displaystyle\sum_{i=1}^k\alpha(X_1,\dots,D_XX_i,\dots,X_k)
\end{equation}

for $X_1,\dots,X_k\in\Gamma(\tau_V)$. Note that
$$D_X(fX_i)=\rho_A(X)(f)X_i+fD_XX_i,\makebox[1cm]{for}f\in
C^{\infty}(M),$$ and, thus, $D_X\alpha\in\Gamma(\wedge^kV^*)$.

\begin{remark}\label{obs7.4}{\em The canonical inclusion $i_V:V\to\widetilde{A}$ is
a Lie algebroid monomorphism over the identity of $M$ (see
(\ref{+})) and, moreover, one may choose
$\tilde{\alpha}\in\Gamma(\wedge^k \tau _{A^+})$ such that
$\alpha=i_V^*\tilde{\alpha}.$ Then, from (\ref{+}) and
(\ref{Dalpha}), we deduce that
\begin{equation}\label{DXderlie}
D_X\alpha=i_V^*({\mathcal L}^{\widetilde{A}}_{(i_A\circ
X)}\tilde{\alpha}), \end{equation}
 where $i_A:A\to\widetilde{A}$
is the canonical inclusion and ${\mathcal
L}^{\widetilde{A}}_{(i_A\circ X)}$ is the Lie derivative in the
Lie algebroid $\tau_{\widetilde{A}}:\widetilde{A}\to M$ with
respect to the section $i_A\circ
X:M\to\widetilde{A}$.\hfill$\diamondsuit$}
\end{remark}

Using the operator $D_X$, we have that
\begin{proposition}\label{prop7.5}
Let $(\tau_A:A\to M, \tau_V:V\to M)$ be a symplectic Lie
affgebroid with symplectic section $\Omega$ and $X:M\to A$ be a
section of $A$. Then, the Lie subaffgebroid of
$(\tau_A^{\tau_A}:{\mathcal J}^AA\to A$,
$\tau_V^{\tau_A}:{\mathcal L}^{\tau_A}V\to A)$ considered in
Examples \ref{ex7.2} $i)$ is Lagrangian if and only if
$D_X\Omega=0.$
\end{proposition}
\begin{proof}
Let $Y$ be a section of $\tau_V:V\to M$ and denote by $\tilde{X}$
and $\tilde{Y}$ the sections of
$\tau_{\widetilde{A}}:\widetilde{A}\to M$ given by
$\tilde{X}=i_A\circ X$ and $\tilde{Y}=i_V\circ Y.$ Then, using
some results in \cite{LMM} (see Examples $7.5$ in \cite{LMM}), we
obtain that
$$(T\tilde{X})(\rho_{\widetilde{A}}(\tilde{Y}))=(\tilde{Y}^c-\lcf\tilde{X},\tilde{Y}\rcf_{\widetilde{A}}^v)\circ\tilde{X},$$
where $\tilde{Y}^c\in\mathfrak X(\widetilde{A})$ and
$\lcf\tilde{X},\tilde{Y}\rcf_{\widetilde{A}}^v\in\mathfrak
X(\widetilde{A})$ are the complete and vertical lift of
$\tilde{Y}\in\Gamma(\tau_{\widetilde{A}})$ and
$\lcf\tilde{X},\tilde{Y}\rcf_{\widetilde{A}}\in\Gamma(\tau_{\widetilde{A}})$,
respectively. Since the restriction to $A$ of $\tilde{Y}^c$ and
$\lcf\tilde{X},\tilde{Y}\rcf_{\tilde{A}}^v$ are tangent to $A$ and
$\rho_{\widetilde{A}}(\tilde{Y})=\rho_V(Y),$ it follows  that
$TX(\rho_V(Y))=((\tilde{Y}^c)_{|A}-(\lcf\tilde{X},\tilde{Y}\rcf_{\widetilde{A}}^v)_{|A})\circ
X,$ and thus, using (\ref{lcvs}) and (\ref{lvcn}), we deduce that
$$(Id,TX\circ\rho_V)\circ Y=(Y^{\textbf{c}}-(D_XY)^{\textbf{v}})\circ X,$$
where $Y^{\textbf{c}}\in\Gamma(\tau_V^{\tau_A})$ and
$(D_XY)^{\textbf{v}}\in\Gamma(\tau_V^{\tau_A})$ are the complete
and vertical lift of $Y\in\Gamma(\tau_V)$ and
$D_XY\in\Gamma(\tau_V)$ (see Section \ref{sec6}).

On the other hand, if $Y,Z\in\Gamma(\tau_V)$ then, from
Proposition \ref{prop6.1}, we have that
$$\Omega^{\textbf{c}}(Y^{\textbf{c}}-(D_XY)^{\textbf{v}},Z^{\textbf{c}}-(D_XZ)^{\textbf{v}})=\Omega(Y,Z)^c-\Omega(D_XY,Z)^v-\Omega(Y,D_XZ)^v.$$
Therefore, using (\ref{6.0}), it follows that
$$\begin{array}{lcl}
\Omega^{\textbf{c}}(Y^{\textbf{c}}-(D_XY)^{\textbf{v}},Z^{\textbf{c}}-(D_XZ)^{\textbf{v}})\circ
X&=&\rho_A(X)(\Omega(Y,Z))-\Omega(D_XY,Z)-\Omega(Y,D_XZ)\\
&=&(D_X\Omega)(Y,Z) \end{array}$$ and
$\Omega^{\textbf{c}}(Y^{\textbf{c}}-(D_XY)^{\textbf{v}},Z^{\textbf{c}}-(D_XZ)^{\textbf{v}})\circ
X=0,\makebox[1.5cm]{for all}Y,Z\in\Gamma(\tau_V)$, if and only if
$D_X\Omega=0.$

Consequently, taking into account that the rank of $V$ is
$\displaystyle\frac{1}{2}rank({\mathcal L}^{\tau_A}V)$, we deduce
the result.
\end{proof}

In the particular case when the symplectic section $\Omega$ is a
$1$-coboundary, we obtain the following corollary.

\begin{corollary}\label{cor7.6}
Let $(\tau_A:A\to M, \tau_V:V\to M)$ be a symplectic Lie
affgebroid with symplectic section $\Omega=-d^V\lambda$, $\lambda$
being a section of the vector bundle $\tau_V^*:V^*\to M$ and
suppose that $X:M\to A$ is a section of $\tau_A:A\to M$. Then, the
Lie subaffgebroid of $(\tau_A^{\tau_A}:{\mathcal J}^AA\to A,
\tau_V^{\tau_A}:{\mathcal L}^{\tau_A}V\to A)$ considered in
Examples \ref{ex7.2} $i)$ is Lagrangian if and only if the section
$D_X\lambda$ of $\tau_V^*:V^*\to M$ is a $1$-cocycle of
$\tau_V:V\to M$.
\end{corollary}
\begin{proof}
If $\alpha\in\Gamma(\wedge^k\tau _V^*)$ then, from Remark
\ref{obs7.4} and since $i_V:V\to\widetilde{A}$ is a Lie algebroid
morphism, it follows that $D_X(d^V\alpha)=d^V(D_X\alpha).$ In
particular, this implies that $D_X\Omega=-d^V(D_X\lambda).$ Thus,
using this fact and Proposition \ref{prop7.5}, we deduce the
result.
\end{proof}

Let $\tau_A:A\to M$ be a Lie affgebroid modelled on the Lie
algebroid $\tau_V:V\to M$ and $\rho_A^*(TV^*)$ be the pull-back of
the vector bundle $T\tau_V^*:TV^*\to TM$ over the anchor map
$\rho_A:A\to TM$. Then, $\rho_A^*(TV^*)$ is a Lie affgebroid over
$V^*$ modelled on the vector bundle $\tau_V^{\tau_V^*}:{\mathcal
L}^{\tau_V^*}V\to V^*$ and with affine bundle projection
$\widetilde{\pi_{V^*}}:\rho_A^*(TV^*)\to V^*$. Denote by
$\Omega_V$ the canonical symplectic section of ${\mathcal
L}^{\tau_V^*}V$. As we know (see Example \ref{ex6.6}), the Lie
affgebroid
$(\widetilde{\pi_{V^*}})^{\widetilde{\pi_{V^*}}}:{\mathcal
J}^{\rho_A^*(TV^*)}\rho_A^*(TV^*)\to\rho_A^*(TV^*)$ modelled on
the Lie algebroid
$(\tau_V^{\tau_V^*})^{\widetilde{\pi_{V^*}}}:{\mathcal
L}^{\widetilde{\pi_{V^*}}}({\mathcal
L}^{\tau_V^*}V)\to\rho_A^*(TV^*)$ is symplectic and the complete
lift $\Omega_V^{\textbf{c}}$ of $\Omega_V$ is a symplectic section
of $(\tau_V^{\tau_V^*})^{\widetilde{\pi_{V^*}}}:{\mathcal
L}^{\widetilde{\pi_{V^*}}}({\mathcal
L}^{\tau_V^*}V)\to\rho_A^*(TV^*)$.

Now, suppose that $\tilde{X}:A\to\rho_A^*(TV^*)$ is a section of
the vector bundle $pr_{1|\rho_A^*(TV^*)}:\rho_A^*(TV^*)\to A$.
Then, $\tilde{X}$ allows us to define the Lie subaffgebroid of
$((\widetilde{\pi_{V^*}})^{\widetilde{\pi_{V^*}}}:\;\;$ ${\mathcal
J}^{\rho_A^*(TV^*)}\rho_A^*(TV^*)\to\rho_A^*(TV^*),(\tau_V^{\tau_V^*})^{\widetilde{\pi_{V^*}}}:{\mathcal
L}^{\widetilde{\pi_{V^*}}}({\mathcal
L}^{\tau_V^*}V)\to\rho_A^*(TV^*))$ considered in E\-xam\-ples
\ref{ex7.2} $ii)$.

Next, we will obtain a necessary and sufficient condition for such
a Lie subaffgebroid to be Lagrangian. For this purpose, we will
introduce a section of $({\mathcal L}^{\tau_A}V)^*\to A$ as
follows.

Let $A_A:\rho_A^*(TV^*)\rightarrow({\mathcal L}^{\tau_A}V)^*$ be
the canonical isomorphism, over the identity of $A$, between the
vector bundles $\rho_A^*(TV^*)\rightarrow A$ and $({\mathcal
L}^{\tau_A}V)^*\rightarrow A$ considered in Section \ref{sec5}
(see (\ref{AA})) and $\alpha_{\tilde{X}}$ be the section of
$({\mathcal L}^{\tau_A}V)^*\rightarrow A$ given by
\begin{equation}\label{alphaX}
\alpha_{\tilde{X}}=A_A\circ\tilde{X}.
\end{equation}
Then, we have the following result.

\begin{proposition}\label{prop7.7}
Let $\tau_A:A\to M$ be a Lie affgebroid modelled on the Lie
algebroid $\tau_V:V\to M$ and $\tilde{X}:A\to\rho_A^*(TV^*)$ be a
section of $pr_{1|\rho_A^*(TV^*)}:\rho_A^*(TV^*)\to A$. Then, the
Lie subaffgebroid of
$((\widetilde{\pi_{V^*}})^{\widetilde{\pi_{V^*}}}:\;\;$ ${\mathcal
J}^{\rho_A^*(TV^*)}\rho_A^*(TV^*)\to\rho_A^*(TV^*),(\tau_V^{\tau_V^*})^{\widetilde{\pi_{V^*}}}:{\mathcal
L}^{\widetilde{\pi_{V^*}}}({\mathcal
L}^{\tau_V^*}V)\to\rho_A^*(TV^*))$ considered in Examples
\ref{ex7.2} $ii)$ is Lagrangian if and only if
$\alpha_{\tilde{X}}$ is a $1$-cocycle of the Lie algebroid
$\tau_V^{\tau_A}:{\mathcal L}^{\tau_A}V\to A$.
\end{proposition}
\begin{proof} Let $(\Psi_{\tilde{X}},\tilde{X})$ be the
monomorphism between the Lie algebroids $\tau_V^{\tau_A}:{\mathcal
L}^{\tau_A}V\to A$ and
$(\tau_V^{\tau_V^*})^{\widetilde{\pi_{V^*}}}:{\mathcal
L}^{\widetilde{\pi_{V^*}}}({\mathcal
L}^{\tau_V^*}V)\to\rho_A^*(TV^*)$ considered in Examples
\ref{ex7.2} $ii)$.

Suppose that $(x^i)$ are local coordinates on $M$ and that
$\{e_0,e_{\alpha}\}$ is a local basis of
$\Gamma(\tau_{\tilde{A}})$ adapted to $1_A$. Denote by
$(x^i\kern-0.5pt,y^{\alpha})$ (respectively,
$(x^i,y_{\alpha};z^{\alpha}\kern-0.5pt,v_{\alpha})$) the
corresponding coordinates on $A$ (respectively, $\rho_A^*(TV^*)$)
and by $\{\tilde{T}_{\alpha},\tilde{V}_{\alpha}\}$ the
corresponding local basis of $\Gamma(\tau_V^{\tau_A})$ (see
(\ref{tildeT})). If the local expression of
$\tilde{X}:A\to\rho_A^*(TV^*)$ is
\begin{equation}\label{Xtildeloc}
\tilde{X}(x^i,y^{\alpha})=(x^i,\tilde{X}_{\alpha};y^{\alpha},\tilde{X}'_{\alpha})
\end{equation}
then, using (\ref{tildeT}) and (\ref{bcv}), we deduce that
\begin{equation}\label{comput}
\kern-15pt\begin{array}{rcl}
\Psi_{\tilde{X}}(\tilde{T}_{\alpha}(x^i,y^\gamma))&=&\tilde{e}^{\textbf{
c}}_\alpha(\tilde{X}(x^i,y^\gamma)) +
\rho_\alpha^j(x^i)\displaystyle\frac{\partial
\tilde{X}_\beta}{\partial
x^j}_{|(x^i,y^\gamma)}\bar{e}_\beta^{\textbf{
c}}(\tilde{X}(x^i,y^\gamma))\\[8pt]
&-&C_{0\alpha}^{\gamma}(x^i)\tilde{e}_{\gamma}^{\textbf{
v}}(\tilde{X}(x^i,y^\gamma))+ C_{\alpha\beta}^\gamma(x^i)y^\beta
\tilde{e}_\gamma^{\textbf{v}}(\tilde{X}(x^i,y^\gamma))\\[8pt]
&+& \rho_\alpha^i(x^i)\displaystyle\frac{\partial
\tilde{X}'_\beta}{\partial
x^i}_{|(x^i,y^\gamma)}\bar{e}_\beta^{\textbf
v}(\tilde{X}(x^i,y^\gamma)),\\[8pt]
\Psi_{\tilde{X}}(\tilde{V}_{\alpha}(x^i,y^\gamma))&=&\displaystyle\frac{\partial
\tilde{X}_\beta}{\partial
y^\alpha}_{|(x^i,y^\gamma)}\bar{e}_\beta^{\textbf
c}(\tilde{X}(x^i,y^\gamma))+ \tilde{e}_\alpha^{\textbf
v}(\tilde{X}(x^i,y^\gamma))
\\[8pt]
&+& \displaystyle\frac{\partial \tilde{X}'_\beta}{\partial
y^\alpha}_{|(x^i,y^\gamma)}\bar{e}_\beta^{\textbf
v}(\tilde{X}(x^i,y^\gamma)),
\end{array}
\end{equation}

where $\tilde{e}_{\alpha}^{\textbf{c}}$ and
$\bar{e}_{\alpha}^{\textbf{c}}$ (respectively,
$\tilde{e}_{\alpha}^{\textbf{v}}$ and
$\bar{e}_{\alpha}^{\textbf{v}}$) are the complete lifts
(respectively, vertical lifts) of $\tilde{e}_{\alpha}$ and
$\bar{e}_{\alpha}$ to
$\Gamma((\tau_V^{\tau_{V^*}})^{\widetilde{\pi_{V^*}}})$ and
$\rho_0^i,$ $\rho_{\alpha}^i$, $C_{0\alpha}^{\gamma}$ and
$C_{\alpha\beta}^{\gamma}$ are the structure functions of the Lie
algebroid $\tau_{\widetilde{A}}:\widetilde{A}\to M$ with respect
to the local coordinates $(x^i)$ and to the basis
$\{e_0,e_{\alpha}\}$.

Thus, if $\lambda_V$ is the Liouville section of ${\mathcal
L}^{\tau_V^*}V$ then, from (\ref{OmegaVC}) and (\ref{comput}), we
obtain that
$$(\Psi_{\tilde{X}},\tilde{X})^*(\lambda_V^{\textbf
c})=(\tilde{X}'_{\alpha} + \tilde{X}_{\gamma} C_{\alpha
\beta}^{\gamma} y^{\beta}+\tilde{X}_{\gamma}C_{\alpha 0}^{\gamma}
)\tilde{T}^{\alpha} + \tilde{X}_{\alpha} \tilde{V}^{\alpha},$$
$\{\tilde{T}^{\alpha},\tilde{V}^{\alpha}\}$ being the dual basis
of $\{\tilde{T}_{\alpha},\tilde{V}_{\alpha}\}$.

 Therefore, using (\ref{AA}), (\ref{alphaX}) and
(\ref{Xtildeloc}), it follows that
$$(\Psi_{\tilde{X}},\tilde{X})^*(\lambda_V^{\textbf
c})=\alpha_{\tilde{X}}.$$ Now, since $\Omega_{V}^{\textbf c} =
(-d^{{\mathcal L}^{\tau_V^*}V} \lambda_{V})^{\textbf c} =
-d^{{\mathcal L}^{\widetilde{\pi_{V^*}}}({\mathcal
L}^{\tau_V^*}V)} \lambda_{V}^{\textbf c}$ (see Proposition
\ref{prop6.1}), we have that
\begin{equation}\label{e7.79}
(\Psi_{\tilde{X}},\tilde{X})^*(\Omega_{V}^{\textbf c}) =
-d^{{\mathcal L}^{\tau_A}V} \alpha_{\tilde{X}}.
\end{equation}

Consequently, using (\ref{e7.79}) and the fact that
$rank({\mathcal L}^{\tau_A}V)=\frac{1}{2}rank({\mathcal
L}^{\widetilde{\pi_{V^*}}}({\mathcal L}^{\tau_V^*}V))$, we deduce
the result.
\end{proof}

\setcounter{equation}{0}
\section{Lagrangian submanifolds, Tulczyjew's triple and
Euler-La\-gran\-ge (Hamilton) equations}\label{sec8}

Let $(\tau_A:A\to M,\tau_V\to M,
(\lcf\cdot,\cdot\rcf_V,D,\rho_A))$ be a symplectic Lie affgebroid
with symplectic section $\Omega$. Then, as we know (see Theorem
\ref{th6.3}), the Lie affgebroid $(\tau_A^{\tau_A}:{\mathcal
J}^AA\to A, \tau_V^{\tau_A}:{\mathcal L}^{\tau_A}V\to A)$ is
symplectic and the complete lift $\Omega^{\textbf{c}}$ of $\Omega$
is a symplectic section of $\tau_V^{\tau_A}:{\mathcal
L}^{\tau_A}V\to A$.

\begin{definition} Let $S$ be a submanifold of the symplectic Lie affgebroid $A$ and $i:S\to A$ be the canonical inclusion.
Denote by $\tau_A^S:S\to M$ the map given by $\tau_A^S=\tau_A\circ
i$ and suppose that
$\rho_V(V_{\tau_A^S(a)})+(T_a\tau_A^S)(T_aS)=T_{\tau_A^S(a)}M$,
for all $a\in S$. Then, the submanifold $S$ is said to be
Lagrangian if the corresponding Lie subaffgebroid
$(\widetilde{\pi_S}:\rho_A^*(TS)\to S, \tau_V^{\tau_A^S}:{\mathcal
L}^{\tau_A^S}V\to S)$ of the symplectic Lie affgebroid
$(\tau_A^{\tau_A}:{\mathcal J}^AA\to A, \tau_V^{\tau_A}:{\mathcal
L}^{\tau_A}V\to A)$ is Lagrangian.
\end{definition}

Now, we have the following result.

\begin{corollary}\label{cor8.2}
Let $(\tau_A:A\to M, \tau_V:V\to M)$ be a symplectic Lie
affgebroid with symplectic section $\Omega=-d^V\lambda$, $\lambda$
being a section of the vector bundle $\tau_V^*:V^*\to M$ and
suppose that $X:M\to A$ is a section of $\tau_A:A\to M$. Then, the
submanifold $S=X(M)$ of $A$ is Lagrangian if and only if the
section $D_X\lambda$ of $\tau_V^*:V^*\to M$ is a $1$-cocycle of
$\tau_V:V\to M$.
\end{corollary}
\begin{proof}
Let $(Id,TX\circ \rho_A):A\to {\mathcal J}^AA$ (respectively,
$(Id,TX\circ \rho_V):V\to{\mathcal L}^{\tau_A}V$) be the morphism
defined as in Examples \ref{ex7.2} $i).$ Then, a direct
computation proves that $(Id, TX\circ \rho_A)(A)=\rho_A^*(TS)$ and
$(Id,TX\circ \rho_V)(V)={\mathcal L}^{\tau_A^S}V. $ Thus, using
Corollary \ref{cor7.6}, we deduce that $S$ is Lagrangian if and
only if  $D_X\lambda$ is a $1$-cocycle.
\end{proof}
 If $\tau_A:A\to M$ is a Lie
affgebroid modelled on the Lie algebroid $\tau_V:V\to M$, we will
denote by $\rho_A^*(TV^*)$ the pull-back of the vector bundle
$T\tau_V^*:TV^*\to TM$ over the anchor map $\rho_A:A\to TM$ and by
$A_A:\rho_A^*(TV^*)\rightarrow({\mathcal L}^{\tau_A}V)^*$ the
canonical isomorphism, over the identity of $A$, between the
vector bundles $\rho_A^*(TV^*)\rightarrow A$ and $({\mathcal
L}^{\tau_A}V)^*\rightarrow A$ considered in Section \ref{sec5}
(see (\ref{AA})).

\begin{corollary}\label{cor8.3}
Let $\tau_A:A\to M$ be a Lie affgebroid modelled on the Lie
algebroid $\tau_V:V\to M$ and $\tilde{X}:A\to\rho_A^*(TV^*)$ be a
section of $pr_{1|\rho_A^*(TV^*)}:\rho_A^*(TV^*)\to A$. Denote by
$S$ the submanifold $\tilde{X}(A)$ of $\rho_A^*(TV^*)$ and by
$\alpha_{\tilde{X}}$ the section of $({\mathcal
L}^{\tau_A}V)^*\rightarrow A$ given by
$\alpha_{\tilde{X}}=A_A\circ\tilde{X}$. Then, $S$ is a Lagrangian
submanifold of the symplectic Lie affgebroid
$(\widetilde{\pi_{V^*}}:\rho_A^*(TV^*)\rightarrow
V^*,\tau_V^{\tau_V^*}:{\mathcal L}^{\tau_V^*}V\to V^*)$ if and
only if $\alpha_{\tilde{X}}$ is a $1$-cocycle of the Lie algebroid
$\tau_V^{\tau_A}:{\mathcal L}^{\tau_A}V\to A$.
\end{corollary}

\begin{proof}
Let $((Id,T(\widetilde{\pi_{V^*}}\circ
\tilde{X})),T\tilde{X}):{\mathcal J}^AA\to {\mathcal
J}^{\rho_A^*(TV^*)}\rho_A^*(TV^*)$ (respectively,
$((Id,T(\widetilde{\pi_{V^*}}\circ
\tilde{X})),T\tilde{X}):{\mathcal L}^{\tau_A}V\to {\mathcal
L}^{\widetilde{\pi_{V^*}}}({\mathcal L}^{\tau^*_V}V))$ be the
morphism defined as in Examples \ref{ex7.2} $ii)$. Then, it is
easy to prove that
\[
\begin{array}{lcr}
((Id,T(\widetilde{\pi_{V^*}}\circ
\tilde{X})),T\tilde{X})({\mathcal
J}^AA)&=&\rho^*_{\rho^*_A(TV^*)}(TS),\\
((Id,T(\widetilde{\pi_{V^*}}\circ
\tilde{X})),T\tilde{X})({\mathcal L}^{\tau_A}V))&=& {\mathcal
L}^{\tau^S_{\rho^*_A(TV^*)}} ({\mathcal L}^{\tau^*_V}V).
\end{array}
\]
Thus, using Proposition \ref{prop7.7}, we deduce that $S$ is
Lagrangian if and only if $\alpha_{\tilde{X}}$ is a $1$-cocycle.
\end{proof}

 Now, let $(\tau_A:A\to M,\tau_V:V\to M)$ be
a Lie affgebroid and $h:V^*\to A^+$ be a Hamiltonian section, that
is, $h$ is a section of $\mu:A^+\to V^*$. Then, we can consider
the cosymplectic structure $(\Omega_h,\eta)$ on ${\mathcal
L}^{\tau_V^*}\widetilde{A}$ given by (\ref{Omegah}) and
(\ref{eta}) and the Reeb section
$R_h\in\Gamma(\tau_{\widetilde{A}}^{\tau_V^*})$ of
$(\Omega_h,\eta)$ (see (\ref{reeb})). Since $\eta(R_h)=1$, we have
that  $R_h(V^*)\subseteq\rho_A^*(TV^*)$. Moreover, from
(\ref{reeb}), it follows that
\begin{equation}\label{Llambda}
{\mathcal L}^{{\mathcal
L}^{\tau_V^*}\widetilde{A}}_{R_h}\lambda_h=d^{{\mathcal
L}^{\tau_V^*}\widetilde{A}}(\lambda_h(R_h)).
\end{equation}
Thus, using (\ref{3.2'}), (\ref{DXderlie}), (\ref{Llambda}) and
Corollary \ref{cor8.2}, we deduce that $S_h=R_h(V^*)$ is a
Lagrangian submanifold of $\rho_A^*(TV^*)$.

On the other hand, it is clear that there exists a bijective
correspondence $\Psi_h$ between the set of curves in $S_h$ and the
set of curves in $V^*$. In fact, if $c: I \to V^*$ is a curve in
$V^*$ then the corresponding curve in $S_{h}$ is $R_{h} \circ c: I
\to S_{h}$.

A curve $\gamma$ in $S_h$,
\[
\gamma:I\to S_h\subseteq\rho_A^*(TV^*)\subseteq {\mathcal
L}^{\tau_V^*}\widetilde{A}\subseteq \widetilde{A}\times TV^*,\;\;
t\mapsto (\gamma_1(t),\gamma_2(t)),
\]
is said to be {\it admissible } if the curve $\gamma_2:I\to TV^*$
is a tangent lift, that is, $\gamma_2(t)=\dot{c}(t)$, where
$c:I\to V^*$ is the curve in $V^*$ given by $\pi_{V^*}\circ
\gamma_2,$ $\pi_{V^*}:TV^*\to V^*$ being the canonical projection.

\begin{theorem}
Under the bijection $\Psi_h$, the admissible curves in the
Lagrangian submanifold $S_h$ correspond with the solutions of the
Hamilton equations for $h$.
\end{theorem}
\begin{proof}
Suppose that $\gamma:I\to S_h\subseteq {\mathcal
L}^{\tau_V^*}\widetilde{A}\subseteq \widetilde{A}\times TV^*$,
$\gamma(t)=(\gamma_1(t),\gamma_2(t))$ is an admissible curve in
$S_h.$ Then, $\gamma_2(t)=\dot{c}(t)$, for all $t$, where $c:I\to
V^*$ is the curve in $V^*$ given by $c=\pi_{V^*}\circ \gamma_2.$

Now, since $R_h$ is a section of the vector bundle
$\tau_{\widetilde{A}}^{\tau_V^*}:{\mathcal
L}^{\tau_V^*}\widetilde{A}\to V^*$ and $\gamma(I)\subseteq
S_h=R_h(V^*),$ it follows that
\begin{equation}\label{cgamma0}
R_h(c(t))=\gamma(t),\;\;\;\mbox{ for all } t,
\end{equation}
that is, $c=\Psi_h(\gamma).$ Thus, from (\ref{cgamma0}), we obtain
that $\rho^{\tau_V^*}_{\widetilde{A}}(R_h)\circ
c=\gamma_2=\dot{c}$. In other words, $c$ is an integral curve of
the vector field $\rho^{\tau_V^*}_{\widetilde{A}}(R_h)$ and,
therefore, $c$ is a solution of the Hamilton equations associated
with $h$ (see Section \ref{sec3.1}).

Conversely, assume that $c:I\to V^*$ is a solution of the Hamilton
equations associated with $h$, that is, $c$ is an integral curve
of the vector field $\rho^{\tau_V^*}_{\widetilde{A}}(R_h)$ or,
equivalently,
\begin{equation}\label{ccdot0}
\rho^{\tau_V^*}_{\widetilde{A}}(R_h)\circ c=\dot{c}.
\end{equation}
Then, $\gamma=R_h\circ c$ is a curve in $S_h$ and, from
(\ref{ccdot0}), we deduce that $\gamma$ is admissible.
\end{proof}

Next, suppose that $L:A\to \R$ is a Lagrangian function. Then,
from Corollary \ref{cor8.3}, we obtain that $S_L=(A_A^{-1}\circ
d^{{\mathcal L}^{\tau_A}V}L)(A)$ is a Lagrangian submanifold of
the symplectic Lie affgebroid $\rho_A^*(TV^*).$

On the other hand, we have a bijective correspondence $\Psi_L$
between the set of curves in $S_L$ and the set of curves in $A$.
In fact, if $\gamma: I \to S_{L}$ is a curve in $S_{L}$ then there
exists a unique curve $c: I \to A$ in $A$ such that
$A_{A}(\gamma(t)) = (d^{{\mathcal L}^{\tau_A}V}L)(c(t))$, for all
$t$. Note that
\[
pr_{1}(\gamma(t)) = (\tau^{\tau_A}_V)^*(A_A(\gamma(t))) =
(\tau^{\tau_A}_V)^*((d^{{\mathcal L}^{\tau_A}V}L)(c(t))) = c(t),
\makebox[.25cm]{} \mbox{ for all } t,
\]
where $pr_{1}: \rho_A^*(TV^*)\subseteq A \times TV^* \to A$ is the
canonical projection on the first factor  and
$(\tau^{\tau_A}_V)^*: ({\mathcal L}^{\tau_A}V)^* \to A$ is the
vector bundle projection. Thus,
\[
\gamma(t) = (c(t), \gamma_{2}(t)) \in \rho_A^*(TV^*) \subseteq
{\mathcal L}^{\tau_V^*}A \subseteq A \times TV^*, \makebox[.3cm]{}
\mbox{ for all } t.
\]
A curve $\gamma$ in $S_L$
\[
\gamma:I\to S_L\subseteq \rho_A^*(TV^*)\subseteq A\times
TV^*,\;\;\; t\mapsto (c(t),\gamma_2(t)),
\]
is said to be {\it admissible} if the curve $\gamma_2:I\to TV^*$
is a tangent lift, that is, $\gamma_2(t)=\dot{c}^*(t),$ where
$c^*:I\to V^*$ is the curve in $V^*$ given by $c^* =
\pi_{V^*}\circ \gamma_2.$

\begin{theorem}
Under the bijection $\Psi_L$, the admissible curves in the
Lagrangian submanifold $S_L$ correspond with the solutions of the
Euler-Lagrange equations for $L.$
\end{theorem}
\begin{proof}
Suppose that $(x^i)$ are local coordinates on $M$ and that
$\{e_0,e_\alpha\}$  is a local basis of $\Gamma(\tau_{\tilde{A}})$
adapted to $1_A$. Denote by $(x^i,y^\alpha)$ (respectively,
$(x^i,y_\alpha)$ and $(x^i,y_\alpha;z^\alpha,v_\alpha))$ the
co\-rres\-pon\-ding coordinates on $A$ (respectively, $V^*$ and
$\rho_A^*(TV^*)$). Then, using (\ref{AA}), it follows that the
submanifold $S_L$ is cha\-rac\-te\-rized by the following
equations
\begin{equation}\label{e8.80^0}
y_\alpha=\frac{\partial L}{\partial y^\alpha},\;\;\;
z^\alpha=y^\alpha,\;\;\; v_\alpha=\rho^i_\alpha\frac{\partial
L}{\partial x^i}+(C_{0\alpha}^{\gamma}+C_{\beta\alpha}^\gamma
y^\beta)\frac{\partial L}{\partial y^\gamma},\makebox[1.5cm]{for
all}\alpha.
\end{equation}

Now, let $\gamma: I \to S_{L}$ be an admissible curve in $S_{L}$
\[
\gamma(t) = (c(t), \gamma_{2}(t)) \in S_{L} \subseteq
\rho_A^*(TV^*)\subseteq A \times TV^*, \mbox{ for all } t,
\]
and denote by $c^*: I \to V^*$ the curve in $V^*$ satisfying
\begin{equation}\label{e8.80^2}
\gamma_{2}(t) = \dot{c}^*(t), \mbox{ for all } t,
\end{equation}
i.e.,
\begin{equation}\label{e8.80^3}
c^*(t) = \pi_{V^*}(\gamma_{2}(t)), \mbox{ for all } t.
\end{equation}
If the local expressions of $\gamma$ and $c$ are
\[
\gamma(t) = (x^{i}(t), y_{\alpha}(t); z^{\alpha}(t),
v_{\alpha}(t)), \makebox[.3cm]{} c(t) = (x^{i}(t), y^{\alpha}(t)),
\]
then we have that
\begin{equation}\label{e8.80^5}
y^{\alpha}(t) = z^{\alpha}(t), \mbox{ for all } \alpha.
\end{equation}
Moreover, from (\ref{e8.80^2}) and (\ref{e8.80^3}), we deduce that
\begin{equation}\label{e8.80^6}
c^*(t) = (x^i(t), y_{\alpha}(t)), \makebox[.3cm]{} \gamma_{2}(t) =
\displaystyle \frac{dx^i}{dt} \frac{\partial}{\partial
x^i}_{|c^*(t)} + \frac{dy_{\alpha}}{dt} \frac{\partial}{\partial
y_{\alpha}}_{|c^*(t)}.
\end{equation}
Thus,
\begin{equation}\label{e8.80^7}
v_{\alpha}(t) = \displaystyle \frac{dy_{\alpha}}{dt}, \mbox{ for
all } \alpha.
\end{equation}
Therefore, using (\ref{e8.80^0}), (\ref{e8.80^5}),
(\ref{e8.80^6}), (\ref{e8.80^7}) and the fact that $\rho_A(c(t)) =
(T\tau_V^*)(\gamma_{2}(t))$, it follows that
\[
\frac{d x^i}{dt}=\rho_0^i+\rho_\alpha^i y^\alpha,\;\;\;
\frac{d}{dt}(\frac{\partial L}{\partial
y^\alpha})=\rho_\alpha^i\frac{\partial L}{\partial
x^i}+(C_{0\alpha}^{\gamma}+C_{\beta\alpha}^\gamma
y^\beta)\frac{\partial L}{\partial y^\gamma},
\]
for all $i$ and $\alpha$, that is, $c$ is a solution of the
Euler-Lagrange equations for $L$.

Conversely, let $c: I \to A$ be a solution of the Euler-Lagrange
equations for $L$ and $\gamma: I \to S_{L}$ be the corresponding
curve in $S_{L}$, $c = \Psi_{L}(\gamma)$. Suppose that
\[
\gamma(t) = (c(t), \gamma_{2}(t)) \in S_L\subseteq
\rho_A^*(TV^*)\subseteq A \times TV^*, \mbox{ for all } t,
\]
and denote by $c^*: I \to V^*$ the curve in $V^*$ given by $c^* =
\pi_{V^*} \circ \gamma_{2}$. If the local expressions of $\gamma$
and $c$ are
\[
\gamma(t) = (x^{i}(t), y_{\alpha}(t); z^{\alpha}(t),
v_{\alpha}(t)), \makebox[.3cm]{} c(t) = (x^{i}(t), y^{\alpha}(t)),
\]
then $y^{\alpha}(t) = z^{\alpha}(t), \mbox{ for all } \alpha$, and
the local expressions of $c^*$ and $\gamma_{2}$ are
\[
c^*(t) = (x^{i}(t), y_{\alpha}(t)), \makebox[.3cm]{} \gamma_{2}(t)
= \displaystyle
(\rho_0^i(x^j(t))+z^{\alpha}(t)\rho_{\alpha}^{i}(x^j(t)))
\frac{\partial}{\partial x^i}_{|c^*(t)} + v_{\alpha}(t)
\frac{\partial}{\partial y_{\alpha}}_{|c^*(t)}.
\]
Thus, using (\ref{e8.80^0}) and the fact that $c$ is a solution of
the Euler-Lagrange equations for $L$, we deduce that
$\gamma_{2}(t) = \dot{c}^*(t)$, for all $t$, which implies that
$\gamma$ is admissible.
\end{proof}

Now, assume that the Lagrangian function $L:A\to \R$ is
hyperregular and denote by $\Theta_L$ and $\Omega_L=-d^{{\mathcal
L}^{\tau_A}\widetilde{A}}\Theta_L$ the Poincar\'{e}-Cartan sections
associated with $L$. We consider the map $(i_V,Id):{\mathcal
L}^{\tau_A}V\to{\mathcal L}^{\tau_A}\widetilde{A}$ given by
$$(i_V,Id)(v,X_a)=(i_V(v),X_a),\makebox[1.5cm]{for all}(v,X_a)\in({\mathcal
L}^{\tau_A}V)_a,\mbox{ with } a\in A,$$ $i_V:V\to\widetilde{A}$
being the canonical inclusion. We have that $(i_V,Id)$ is a Lie
algebroid morphism over the identity of $A$. Furthermore, if
$\phi_0$ is the section of the dual bundle to ${\mathcal
L}^{\tau_A}\widetilde{A}$ defined by (\ref{T0}), it follows that
the pair $(\Omega_L,\phi_0)$ is a cosymplectic structure on
${\mathcal L}^{\tau_A}\widetilde{A}$ (see (\ref{PCloc})) and it is
easy to prove that $(i_V,Id)^*\phi_0=0$. This implies that
$(i_V,Id)^*\Omega_L$ is a symplectic section of the Lie algebroid
$\tau_V^{\tau_A}:{\mathcal L}^{\tau_A}V\to A$ and, thus, the Lie
affgebroid $(\tau_A^{\tau_A}:{\mathcal J}^AA\to A,
\tau_V^{\tau_A}:{\mathcal L}^{\tau_A}V\to A)$ is symplectic. Note
that $(i_V,Id)^*\Omega_L=-d^{{\mathcal
L}^{\tau_A}V}((i_V,Id)^*(\Theta_L))$.

Next, denote by $R_L$ the Reeb section of the cosymplectic
structure $(\Omega_L,\phi_0)$. Since $\phi_0(R_L)=1$, we deduce
that $R_L$ is a section of $\tau_A^{\tau_A}:{\mathcal J}^AA\to A.$
Moreover, we have that ${\mathcal L}_{R_L}^{{\mathcal
L}^{\tau_A}\widetilde{A}}\Theta_L=d^{{\mathcal
L}^{\tau_A}\widetilde{A}}(\Theta_L(R_L))$. Therefore, from
(\ref{DXderlie}) and Corollary \ref{cor8.2}, we deduce that
$S_{R_L}=R_L(A)$ is a Lagrangian submanifold of the symplectic Lie
affgebroid ${\mathcal J}^AA.$

On the other hand, it is clear that there exists a bijective
correspondence $\Psi_{S_{R_L}}$ between the set of curves in
$S_{R_L}$ and the set of curves in $A$.

A curve $\gamma$ in $S_{R_L}$
\[
\gamma:I\to S_{R_L}\subseteq {\mathcal J}^AA\subseteq A\times
TA,\;\;\; t\mapsto (\gamma_1(t),\gamma_2(t)),
\]
is said to be {\it admissible} if the curve $\gamma_2:I\to TA$ is
a tangent lift, that is, $\gamma_2(t)=\dot{c}(t),\mbox{ for all }
t,$ where $c:I\to A$ the curve in $A$ defined by $c=\pi_A\circ
\gamma_2,$ $\pi_A:TA\to A$ being the canonical projection.

\begin{theorem}
If the Lagrangian $L$ is hyperregular then under the bijection
$\Psi_{S_{R_L}}$ the admissible curves in the Lagrangian
submanifold $S_{R_L}$ correspond with the solutions of the
Euler-Lagrange equations for $L.$
\end{theorem}
\begin{proof}
Let $\gamma:I\to S_{R_L}\subseteq {\mathcal J}^AA\subseteq A\times
TA$ be an admissible curve in $S_{R_L}$,
$\gamma(t)=(\gamma_1(t),\gamma_{2}(t))$, for all $t$. Then,
$\gamma_2(t)=\dot{c}(t),\mbox{ for all } t,$ where $c:I\to A$ is
the curve in $A$ given by $c=\pi_A\circ \gamma_2.$

Now, since $R_L$ is a section of the vector bundle
$\tau_A^{\tau_A}:{\mathcal J}^AA\to A$ and $\gamma(I)\subseteq
S_{R_L}=R_L(A),$ it follows that $R_L(c(t))=\gamma(t)$, for all
$t$, that is, $c=\Psi_{S_{R_L}}(\gamma).$ Thus, we obtain that
$\rho_{\widetilde{A}}^{\tau_A}(R_L)\circ c=\gamma_2=\dot{c}$, that
is, $c$ is an integral curve of the vector field
$\rho_{\widetilde{A}}^{\tau_A}(R_L)$. Therefore, $c$ is a solution
of the Euler-Lagrange equations associated with $L$ (see Section
\ref{sec3.2}).

Conversely, assume that $c:I\to A$ is a solution of the
Euler-Lagrange equations associated with $L$, that is, $c:I\to A$
is an integral curve of the vector field
$\rho_{\widetilde{A}}^{\tau_A}(R_L)$ or, equivalently,
\begin{equation}\label{ccdot1}
\rho_{\widetilde{A}}^{\tau_A}(R_L)\circ c=\dot{c}.
\end{equation}
Then, $\gamma=R_L\circ c$ is a curve in $S_{R_L}$ and, from
(\ref{ccdot1}), we deduce that $\gamma$ is admissible.
\end{proof}

 If $L:A\to \R$ is hyperregular then the Legendre
transformation $leg_L:A\to V^*$ associated with $L$ is a global
diffeomorphism. So, we may consider the Hamiltonian section
$h_L:V^*\to A^+$ defined by $h_L=Leg_L\circ leg_L^{-1}$,
$Leg_L:A\to A^+$ being the extended Legendre transformation, and
the Reeb section $R_L$ (respectively, $R_{h_L}$) of the
cosymplectic structure $(\Omega_L,\phi_0)$ (respectively,
$(\Omega_{h_L},\eta)$) on ${\mathcal L}^{\tau_A}\widetilde{A}$
(respectively, ${\mathcal L}^{\tau_V^*}\widetilde{A}$).

Thus, we have:
\begin{itemize}
\item The Lagrangian submanifolds $S_L$ and $S_{h_L}$ of the symplectic
Lie affgebroid $\rho_A^*(TV^*)$.
\item The Lagrangian submanifold $S_{R_L}$ of the symplectic Lie
affgebroid ${\mathcal J}^AA.$
\end{itemize}
Denote by $({\mathcal L}leg_L,leg_L)$ the Lie algebroid
isomorphism between the Lie algebroids ${\mathcal
L}^{\tau_A}\widetilde{A}$ and ${\mathcal
L}^{\tau_V^*}\widetilde{A}$ induced by the transformation
$leg_L:A\to V^*.$

\begin{theorem}
If the Lagrangian function $L:A\to \R$ is hyperregular and
$h_L:V^*\to A^+$ is the corresponding Hamiltonian section then the
Lagrangian submanifolds $S_L$ and $S_{h_L}$ are equal and
\begin{equation}\label{RelLag}
{\mathcal L}leg_L(S_{R_L})=S_L=S_{h_L}.
\end{equation}
\end{theorem}
\begin{proof}
Using (\ref{Related}), we obtain that
\begin{equation}\label{conrel}
A_A\circ R_{h_L}\circ leg_L=A_A\circ {\mathcal L}leg_L\circ R_L.
\end{equation}
Now, suppose that $(x^i)$ are local coordinates in $M$ and that
$\{e_0,e_\alpha\}$ is a local basis of $\Gamma(\tau_{\tilde{A}})$
adapted to $1_A.$ Denote by $(x^i,y^\alpha)$ the corresponding
coordinates on $A$ and by $(x^i,y^\alpha;z^\alpha,v^\alpha)$
(respectively, $(x^i,y_\alpha;z^\alpha,v_\alpha)$ and
$(x^i,y^\alpha;z_\alpha,v_\alpha)$) the corresponding ones on
${\mathcal J}^A A$ (respectively, $\rho_A^*(TV^*)$ and $({\mathcal
L}^{\tau_A}V)^*$). Then, from (\ref{RLloc}), (\ref{Llegloc}) and
(\ref{AA}), we deduce that
\[
(A_A\circ {\mathcal L}leg_L\circ
R_L)(x^i,y^\alpha)=(x^i,y^\alpha;\rho_\alpha^i\frac{\partial
L}{\partial x^i},\frac{\partial L}{\partial y^\alpha}).\] Thus, it
follows that $ (A_A\circ {\mathcal L}{leg_L}\circ
R_L)(x^i,y^\alpha)=d^{{\mathcal L}^{\tau_A}V}L(x^i,y^\alpha), $
that is (see (\ref{conrel})), $A_A\circ R_{h_L}\circ
leg_L=d^{{\mathcal L}^{\tau_A}V}L.$ Therefore, $S_L=S_{h_L}.$

On the other hand, using (\ref{Related}), we obtain that
(\ref{RelLag}) holds.
\end{proof}

\setcounter{equation}{0}
\section{Applications}\label{sec9}
\subsection{The particular case of a Lie algebroid}\label{sec9.0}In this first
example,  we will show that when we consider a Lie algebroid as a
Lie affgebroid and we apply the results obtained in this paper, we
recover the constructions made in \cite{LMM}.

Let $(E,\lcf\cdot,\cdot\rcf,\rho)$ be a Lie algebroid  on $M$ with
vector bundle projection  $\tau_E:E\to M.$ In this case,
$\tau_E:E\to M$ can be considered as an affine bundle with
associated vector bundle $\tau_E:E\to M.$ Then, its dual bundle
${E^+}$ is $E^*\times\R$ and therefore, its bidual bundle
${\widetilde{E}}$ is just $E\times\R.$ Moreover, we can identify
the distinguished section
$1_E\in\Gamma(\tau_{E^+})=\Gamma(\tau_E^*)\times C^{\infty}(M)$
with the section $(0,1)\in\Gamma(\tau_E^*)\times C^{\infty}(M)$
induced by the constant function $1$ on $M.$

On the other hand, a local basis $\{e_0,e_{\alpha}\}$ of sections
of $\widetilde{E}=E\times \R$ adapted to the $1$-cocycle
$1_E=(0,1)$ can be constructed as follows
$$e_0=(0,1)\in \Gamma(\tau_{E^+})=
\Gamma(\tau_E)\times
C^\infty(M),\;\;\;e_{\alpha}=(e'_{\alpha},0)\in
\Gamma(\tau_{E^+})= \Gamma(\tau_E)\times C^\infty(M),$$ where
$\{e'_{\alpha}\}$ is a local basis of sections of $E$.

Now, suppose that $(x^i)$ are local coordinates on an open subset
$U$ of $M$. Denote by $(x^i,y^{\alpha},y^0)$ the  local
coordinates on $\widetilde{E}=E\times\R$ induced by
$\{e_0,e_{\alpha}\}$. The local equation defining $E$ as affine
subbundle (respectively, as vector subbundle) of $\widetilde{E}$
is $y^0=1$ (respectively, $y^0=0$). Thus, $(x^i,y^{\alpha})$ may
be considered as local coordinates on $E$.

In such a case, the local expressions of the Lie bracket and the
anchor map on $\widetilde{E}=E\times \R$ are the following:
$$
\begin{array}{rclcrrcl}
\lcf e_0,e_{\alpha}\rcf_{E\times\R}&=&0,&&& \lcf
e_{\alpha},e_{\beta}\rcf_{E\times\R}&=&C_{\alpha\beta}^{\gamma}e_{\gamma},\\
\rho_{E\times\R}(e_0)&=&0,&&&
\rho_{E\times\R}(e_{\alpha})&=&\rho_{\alpha}^i\displaystyle\frac{\partial}{\partial
x^i}.
\end{array}
$$

On the other hand,  the prolongation, ${\mathcal
 L}^{\tau_{E}^*}\widetilde{E}$,
of the bidual Lie algebroid $\widetilde{E}$ over the dual
projection of the vector bundle $\tau_E:E\to M$ is just the
product vector bundle ${\mathcal L}^{\tau_E^*}E\times\R\to E^*.$
So, the space of sections of ${\mathcal
 L}^{\tau_{E}^*}\widetilde{E}$ can be
identified with $\Gamma(\tau_E^{\tau_E^*})\times C^\infty(E^*).$

Since the  map  $\mu:E^+=E^*\times \R\to E^*$ is the canonical
projection on the first factor, a Hamiltonian section $h:E^*\to
E^*\times\R$ may be identified with a Hamiltonian function
$H:E^*\to \R$ is such a way that $h(\beta)=(\beta,-H(\beta)),$ for
$\beta\in E^*.$ Moreover, the cosymplectic structure
$(\Omega_h,\eta)$ on the Lie algebroid ${\mathcal
L}^{\tau_{E}^*}\widetilde{E}$ can be expressed in terms of the
canonical symplectic section of ${\mathcal L}^{\tau_{E}^*}{E},$
$\Omega_E,$ and the section $(0,1)\in
\Gamma(\tau_E^{\tau_E^*})\times C^\infty(E^*)$ as follows
 \begin{equation}\label{9.1'}
\Omega_h=\Omega_E + d^{{\mathcal L}^{\tau_E^*}E} H\wedge
(0,1),\;\;\;\; \eta=(0,1).
\end{equation}
Thus, the local expression of the Reeb section of
$(\Omega_h,\eta)$, $R_h\in
\Gamma(\tau_{\widetilde{E}}^{\tau_E^*})\cong
\Gamma(\tau_{{E}}^{\tau_E^*})\times C^\infty(E^*),$  is
\begin{equation}\label{9.1''}
R_h=(\xi_H,1),
\end{equation}
$\xi_H$ being the unique section of
$\Gamma(\tau_{{E}}^{\tau_E^*})$ satisfying
$i_{\xi_H}\Omega_E=d^{{\mathcal L}^{\tau^*_E}E} H.$ This implies
that the vector fields $\rho_{E\times \R}^{\tau_E^*}(R_h)$ and
$\rho^{\tau_E^*}(\xi_H)$ on $E^*$ coincide. Therefore, one deduces
that the Hamilton equations associated
 with $h$ on $E$ (as a Lie
affgebroid)  are just the Hamilton equations associated with $H$
considering on $E$ the structure of Lie algebroid (see Section
$3.3$ in \cite{LMM}).

Now, let $L:E\to\R$ be a Lagrangian function. Then, if we write
the Euler-Lagrange equations associated with $L$ (see (\ref{EL})),
we obtain the Euler-Lagrange equations associated with $L$
considering $E$ as a Lie algebroid (see $(2.40)$ in \cite{LMM}).

On the other hand, if the Lagrangian function $L$ is hyperregular,
we can consider the corresponding Hamiltonian section $h_L:E^*\to
E^+=E^*\times\R$ defined by $h_L=Leg_L\circ leg_L^{-1}.$ Thus,
$h_L(\alpha)=(\alpha,-H_L(\alpha))$, where $H_L:E^*\to \R$ is a
Hamiltonian function on $E^*$. One may prove that $H_L=E_L\circ
leg_L^{-1}$, where $E_L:E\to \R$ is the Lagrangian energy
associated with $L$ (for the definition of $E_L$, see \cite{LMM}).
Thus, $H_L$ is the Hamiltonian function associated with $L$
considered in \cite{LMM}.

Next, we are going to describe the Tulczyjew' triple associated
with $E$ (as a Lie affgebroid) and the Hamilton section $h:E^*\to
E^*\times \R$:

\begin{equation}\label{d1}
\begin{picture}(480,90)(120,-10)
\put(245,20){\makebox(0,0){$E$}} \put(310,70){\vector(-3,-2){60}}
\put(480,70){\vector(-3,-2){60}}
\put(417,20){\makebox(0,0){$E^*$}}
\put(170,40){$(\tau_E^{\tau_E})^*$}
\put(360,40){$\widetilde{\pi_{E^*}}$}
\put(460,40){$(\tau_E^{\tau_E^*})^*$}\put(180,70){\vector(3,-2){60}}
\put(290,40){$pr_1$} \put(350,70){\vector(3,-2){60}}
\put(175,80){\makebox(0,0){$({\mathcal L}^{\tau_E}E)^* $}}
\put(240,85){$A_E$}\put(400,85){$\flat_{\Omega_{h}}$}\put(280,80){\vector(-1,0){85}}
\put(305,80){\makebox(0,0){$\rho^*(TE^*)$}}
\put(345,80){\makebox(0,0){$\cong{\mathcal L}^{\tau_E^*}E$}}
\put(365,80){\vector(1,0){95}}
\put(480,80){\makebox(0,0){$({\mathcal
L}^{\tau_E^*}E)^*$}}\end{picture}
\end{equation}
\vspace{-45pt}

Note that the space ${\mathcal J}^EE$ is the prolongation
${\mathcal L}^{\tau_E}E$ of $E$ over $\tau_E:E\to M$. Moreover,
the canonical involution $\sigma_E:{\mathcal J}^EE\to {\mathcal
J}^EE$ associated with the Lie affgebroid $E$ is just the
canonical involution associated with the Lie algebroid $E$ which
was introduced in \cite{LMM} (see Section $4$ in \cite{LMM}).
Thus, it is easy to prove that the map
$A_E:\rho_E^*(TE)\equiv{\mathcal L}^{\tau_E^*}E\to ({\mathcal
L}^{\tau_E}E)^*$ coincides with the isomorphism of vector bundles
also introduced in \cite{LMM} (see Section 5 in \cite{LMM}). On
the other hand, let $\flat_{E^*}:{\mathcal L}^{\tau_E^*}E\to
({\mathcal L}^{\tau_E^*}E)^*$ be the canonical isomorphism between
the vector bundles $\widetilde{\pi_{E^*}}:{\mathcal
L}^{\tau_E^*}E\to E^*$ and $(\tau_E^{\tau_E^*})^*:({\mathcal
L}^{\tau_E^*}E)^*\to E^*$ induced by the symplectic section
$\Omega_E.$ Then, one can check that (see (\ref{9.1'}))
$$\flat_{\Omega_h}=\flat_{E^*}-(d^{{\mathcal
L}^{\tau_E^*}E}H)\circ\widetilde{\pi_{E^*}}.$$

Finally, we will analyze  the Lagrangian submanifolds which
describe the dynamics on the Lie affgebroid $E.$

If $L:E\to\R$ is a Lagrangian function then it is clear that
$S_L=(A_E^{-1}\circ d^{{\mathcal L}^{\tau_E}E}L)(E)$ is just the
Lagrangian submanifold of the symplectic Lie algebroid
$\widetilde{\pi_{E^*}}: {\mathcal L}^{\tau_E^*}E\to E^*$ which was
considered in \cite{LMM} in order to describe the Lagrangian
dynamics (see Section $8$ in \cite{LMM}). On the other hand,
suppose that $h:E^*\to E^+=E^*\times \R$ is a Hamiltonian section
and that $H:E^*\to \R$ is the corresponding Hamiltonian function.
Then, the vector bundle $\widetilde{\pi_{E^*}}: {\mathcal
L}^{\tau_E^*}E\to E^*$ may be identified with the affine subbundle
${\mathcal L}^{\tau_E^*}E\times \{1\}\to E^*$ of ${\mathcal
L}^{\tau_E^*}E\times \R={\mathcal L}^{\tau_E^*}\widetilde{E}\to
E^*$ and, under this identification, the Reeb section $R_h$ is the
Hamiltonian section $\xi_H$ of $\widetilde{\pi_{E^*}}:{\mathcal
L}^{\tau_E^*}E\to E^*$ (see (\ref{9.1''})). Therefore,
$S_h=R_h(E^*)$ may be identified with the Lagrangian submanifold
$S_H=\xi_H(E^*)$ of the symplectic Lie algebroid
$\widetilde{\pi_{E^*}}:{\mathcal L}^{\tau_E^*}E\to E^*$ associated
with the Hamiltonian function $H$. This submanifold was considered
in \cite{LMM} in order to describe the Hamiltonian dynamics.

\subsection{Lie affgebroids and time-dependent Mechanics}\label{sec9.1} Let
$\tau:M\to\R$ be a fibration and $\tau_{1,0}:J^1\tau\to M$ be the
associated Lie affgebroid modelled on the vector bundle
$\pi=(\pi_M)_{|V\tau}:V\tau\to M$. As we know, the bidual vector
bundle $\widetilde{J^1\tau}$ to the affine bundle
$\tau_{1,0}:J^1\tau\to M$ may be identified with the tangent
bundle $TM$ to $M$ and, under this identification, the Lie
algebroid structure on $\pi_M:TM\to M$ is the standard Lie
algebroid structure and the $1$-cocycle $1_{J^1\tau}$ on
$\pi_M:TM\to M$ is just the $1$-form $\eta=\tau ^\ast (dt)$, $t$
being the coordinate on $\R$ (see Section \ref{sec1.2}). If
$(t,q^i)$ are local fibred  coordinates on $M$ then $\{
\frac{\partial}{\partial q^i} \}$ (respectively, $\{
\frac{\partial}{\partial t}, \frac{\partial}{\partial q^i} \}$) is
a local basis of sections of $\pi:V\tau\to M$ (respectively,
$\pi_M:TM\to M$). Denote by $(t,q^i,\dot{q}^i)$ (respectively,
$(t,q^i,\dot{t},\dot{q}^i)$) the corresponding local coordinates
on $V\tau$ (respectively, $TM$). Then,  the (local) structure
functions of $TM$ with respect to this local trivialization are
given by
$$
\begin{array}{l}
C^k_{ij}=0 \mbox{ and } \rho ^i_j=\delta _{ij}, \mbox{ for }i,j, k
\in \{ 0,1,\ldots ,n\}.
\end{array}
$$

Now, let $\pi^*:V^*\tau\to M$ be the dual vector bundle to
$\pi:V\tau\to M$ and suppose that $\pi_1^*:V^*\tau\to \R$ is the
fibration defined by $\pi_1^*=\tau\circ \pi^*$, that
$(\pi_1^*)_{1,0}:J^1\pi_1^*\to V^*\tau$ is the $1$-jet bundle of
local sections of $\pi_1^*:V^*\tau \to \R$ and that
$\rho_{J^1\tau}:J^1\tau\to TM$ is the anchor map of
$\tau_{1,0}:J^1\tau\to M$ ($\rho_{J^1\tau}$ is the canonical
inclusion of $J^1\tau$ on $TM$). Then, one may introduce an
isomorphism $F$ (over the identity of $V^*\tau$) between the Lie
affgebroids $(\pi_1^*)_{1,0}:J^1\pi_1^*\to V^*\tau$ and
$\widetilde{\pi_{V^*\tau}}:\rho_{J^1\tau}^*(T(V^*\tau))\to
V^*\tau$ as follows. If $j_t^1\psi\in J^1\pi_1^*$, with $t\in I$
and $\psi:I\subseteq \R\to V^*\tau$ is a local section of
$\pi_1^*:V^*\tau\to \R$, then there exists a unique $z\in J^1\tau$
such that $\rho_{J^1\tau}(z)=(T_{\psi(t)}\pi^*)(\dot\psi(t))$ and
we define
\[
F(j_t^1\psi)=(z,(T_{\psi(t)}\pi^*)(\dot\psi(t))).
\]
On the other hand, since the anchor map of $\pi_{M}:TM\to M$ is
the identity of $TM$, it follows that the Lie algebroids
$\pi_M^{\pi^*}:{\mathcal L}^{\pi^*}(TM)\to V^*\tau$ and
$\pi_{V^*\tau}:T(V^*\tau)\to V^*\tau$ are isomorphic.

Note that if $(t,q^i,p_i;\dot{t},\dot{q}^i,\dot{p}_i)$ are the
local coordinates on $T(V^*\tau)$ induced by $(t,q^i,p_i)$ then
the local equation defining $J^1\pi_1^*$ as an affine subbundle of
$\pi_{V^*\tau}:T(V^*\tau)\to V^*\tau$ is $\dot{t}=1$. Therefore,
$(t,q^i,p_i; \dot{q}^i,\dot{p}_i)$ is a system of local
coordinates on $J^1\pi_1^*.$

Next, let $h$ be a Hamiltonian section, that is,  $h:V^*\tau \to
(J^1\tau)^+\cong T^\ast M$ is a section of the canonical
projection of $\mu:(J^1\tau)^+\cong T^*M \to V^*\tau$. $h$ is
locally given by
\[
h(t,q^i,p_i)=(t,q^i,-H(t,q^j,p_j),p_i).
\]
Moreover, the cosymplectic structure $(\Omega_h,\eta)$ on the Lie
algebroid $\pi_M^{\pi^*}: {\mathcal L}^{\pi^*}(TM)\cong T(V^*\tau)
\to V^*\tau$ is, in this case, the standard cosymplectic structure
$(\Omega_h,\eta)$ on the manifold $V^*\tau$ locally given by
\[
\Omega_h=dq^i\wedge dp_i + \frac{\partial H}{\partial
q^i}dq^i\wedge dt + \frac{\partial H}{\partial p_i}dp_i\wedge
dt,\;\;\; \; \eta=dt.
\]
Thus, the Reeb section of $(\Omega_h,\eta)$ is the vector field
$R_h$ on $V^*\tau$ defined by
\[
R_h=\frac{\partial }{\partial t}+\frac{\partial H}{\partial
p_i}\frac{\partial }{\partial q^i}-\frac{\partial H}{\partial
q^i}\frac{\partial }{\partial p_i}.
\]
It is clear that the integral sections of $R_h$ \[ t\mapsto
(t,q^i(t),p_i(t))\] are just the solutions of the classical
non-autonomous Hamilton equations
\[
\frac{dq^i}{dt}=\frac{\partial H}{\partial p_i},\quad
\frac{dp_i}{dt}=-\frac{\partial H}{\partial q^i}.
\]
Moreover, if  $S_h$ is the Lagrangian submanifold of the
symplectic Lie affgebroid $\widetilde{\pi_{V^*\tau}}:\rho ^\ast
_{J^1\tau}(V^*\tau)$ $\cong J^1\pi^*_1\to V^*\tau$ given by
$S_h=R_h(V^*\tau )$ then the local equations defining $S_h$ are
\[
\begin{array}{l}
\dot{q}^i=\displaystyle\frac{\partial H}{\partial p_i},\;\;\;\;
\dot{p}_i=-\displaystyle\frac{\partial H}{\partial q^i},
\end{array}
\]
that is, the Hamilton equations for $h$.

\begin{remark}{\em As we know (see Section \ref{sec3.1}),
$(J^1\tau)^+\cong T^*M$ is an affine bundle over $V^*\tau$ of rank
$1$ modelled on the trivial vector bundle
$\tau_{V^*\tau\times\R}:V^*\tau\times\R\to V^*\tau$ and the affine
bundle projection is the map $\mu:(J^1\tau)^+\cong T^*M\to
V^*\tau$. Furthermore, the Hamiltonian section $h$ induces an
affine function $F_h$ on $T^*M$. In the particular case when the
fibration $\tau$ is trivial, that is, $M=\R\times Q$ and $\tau$ is
the canonical projection on the first factor then $T^*M\cong
T^*(\R\times Q)\cong(\R\times\R)\times T^*Q$,
$V^*\tau\cong\R\times T^*Q$ and under these identifications $\mu$
is the projection given by
$$\mu(t,p,\alpha_q)=(t,\alpha_q),\makebox[1cm]{for}\alpha_q\in
T^*_qQ\makebox[1cm]{and}(t,p)\in\R\times\R.$$ Thus, $h$ may be
considered as the Hamiltonian function $H$ on $\R\times T^*Q$. In
addition, the affine function $F_h\equiv F_H$ on
$T^*M\cong(\R\times\R)\times T^*Q$ is given by
$$F_h(t,p,\alpha_q)\equiv F_H(t,p,\alpha_q)=-H(t,\alpha_q)-p.$$
Note that $F_h\equiv F_H$ is, up to the sign, the classical
extension of $H$ to a Hamiltonian function $H^+$ on the extended
phase space $T^*(\R\times Q)$ (see \cite{LL} and the references
therein).\hfill$\diamondsuit$}
\end{remark}

Finally, if  $L:J^1\tau \to \R$ is a time-dependent Lagrangian
function then the local equations defining the corresponding
Lagrangian submanifold $S_L$ are (see (\ref{e8.80^0}))
\begin{equation}\label{def-momenta}
\displaystyle p_i=\frac{\partial L}{\partial \dot{q}^i},
\end{equation}
\begin{equation}\label{EL-eqs}
\dot{p}^i=\displaystyle\frac{\partial L}{\partial q^i} .
\end{equation}
Note that Eqs. (\ref{def-momenta}) give the definition of the
momenta and Eqs. (\ref{EL-eqs}) are just the Euler-Lagrange
equations for L.

\subsection{Atiyah affgebroids and nonautonomous Hamilton(Lagrange)-Poincar\'{e}
equations}\label{sec9.2}
\subsubsection{Atiyah affgebroids}\label{sec9.2.1}
Let $p:Q\to M$ be a principal $G$-bundle. Denote by $\Phi:G\times
Q\to Q$ the free action of $G$ on $Q$ and by $T\Phi:G\times TQ\to
TQ$ the tangent action of $G$ on $TQ.$ Then, one may consider the
quotient vector bundle $\pi_{Q}|G:TQ/G\to M=Q/G$ and the sections
of this vector bundle may be identified with the vector fields on
$Q$ which are invariant under the action $\Phi$. Using that every
$G$-invariant vector field on $Q$ is $p$-projectable and that the
usual Lie bracket on vector fields is closed with respect to
$G$-invariant vector fields, we can induce a Lie algebroid
structure on $TQ/G$. This Lie algebroid is called {\it the Atiyah
algebroid} associated with the principal $G$-bundle $p:Q\to M$
(see \cite{LMM,Ma}).

Now, we suppose that $\nu:M\to \R$ is a fibration of $M$ on $\R$.
Denote by $\tau:Q\to \R$ the composition $\tau=\nu\circ p.$ Then,
$\Phi$ induces an action $J^1\Phi:G\times J^1\tau\to J^1\tau$ of
$G$ on $J^1\tau$ such that
$$J^1\Phi(g,j_t^1\gamma)= j_t^1(\Phi_g\circ \gamma),$$ for $g\in
G$ and $\gamma:I\subset\R\to Q$ a local section of $\tau$, with
$t\in I.$ Moreover, the projection
\[
\tau_{1,0}|G:J^1\tau/G\to M,\;\;\;\; [j^1_t\gamma]\mapsto
p(\tau_{1,0}(j^1_t\gamma))=p(\gamma(t))
\]
defines an affine bundle on $M$ which is modelled on the quotient
vector bundle
\[
\pi|G:V\tau/G\to M,\;\;\; [u_q]\mapsto p(q), \mbox{ for } u_q\in
V_q\tau,
\]
$\pi:V\tau\to Q$ being the vertical bundle of the fibration
$\tau:Q\to \R.$ Here, the action of $G$ on $V\tau$ is the
restriction to $V\tau$ of the tangent action $T\Phi$ of $G$ on
$TQ.$

In addition, the bidual vector bundle of $J^1\tau/G\to M$ is
$\pi_{Q}|G:TQ/G\to M$.

On the other hand, if $t$ is the usual coordinate on $\R$, the
$1$-form $\tau^*(dt)$ is $G$-invariant and defines a non-zero
$1$-cocycle $\phi:TQ/G\to \R$ on the Atiyah algebroid $TQ/G$. Note
that $\phi^{-1}\{1\}\cong J^1\tau|G$ and therefore, one may
consider the corresponding Lie affgebroid structure on $J^1\tau/G$
(see \cite{MMeS}). $J^1\tau/G$ endowed with this structure is
called {\it the Atiyah affgebroid }  associated with the principal
$G$-bundle $p:Q\to M$ and the fibration $\nu:M\to \R$.

Now, let $K:TQ\to \mathfrak g$ be a connection in the principal
bundle $p:Q\to M$. This connection will allow us to determine an
isomorphism between the vector bundles $TQ/G \to M$ and $TM\oplus
\widetilde{\mathfrak g}\to M$, where ${\mathfrak g}$ is the Lie
algebra of $G$ and $\widetilde{\mathfrak g}=(Q\times {\mathfrak
g})/G$ is the adjoint bundle associated with the principal bundle
$p:Q\to M.$  In fact, this isomorphism is defined as follows (see
\cite{LMM,Ma})
$$I_{K}:TQ/G\to TM\oplus \tilde{\mathfrak
g},\;\;\;[{u}_q]\mapsto (T_qp)({u}_q)\oplus [(q,K({u}_q))].$$
Using this identification, we can induce a Lie algebroid structure
on $TM\oplus \widetilde{\mathfrak g}$ in such a way that $I_{K}$
is a Lie algebroid isomorphism. In addition, the 1-cocycle  $\phi$
induces a 1-cocycle $\phi'$ on $TM\oplus \widetilde{\mathfrak g}$
given by  $\phi ' (v_{p(q)}\oplus [(q,\xi)])=(\pi ^\ast dt
)(p(q))(v_{p(q)})$, for $v_{p(q)}\in T_{p(q)}M$ and $\xi\in {\frak
g}.$ It is clear that $(\phi')^{-1}\{1\}$ may be identified with
the affine bundle (over $M$) $J^1\nu\oplus \widetilde{\frak g}.$
Thus, the affine bundle $J^1\nu\oplus \widetilde{\mathfrak g}\to
M$ is a Lie affgebroid and, from Proposition \ref{prop2.2},
$I_{K}$ induces an isomorphism between the Lie affgebroids
$J^1\tau/G$ and $J^1\nu\oplus \widetilde{\mathfrak g}$.

Denote by $B:TQ\oplus TQ\to \mathfrak g$ the curvature of $K$.
Then, we will obtain a local basis of $\Gamma (\pi_Q|G)$ as
follows. Let $U\times G$ be a local trivialization of $p :Q\to M$,
where $U$ is an open subset of $M$ and let $e$ be the identity
element of $G$. Assume that there are local fibred coordinates
$(t,x^i)$ in $U$ and that $\{ \xi _a\}$ is a basis of $\mathfrak
g$. Denote by $\{\xi _a^L\}$ the corresponding left-invariant
vector fields on $G$, that is, $\xi ^L_a(g)=(T_eL_g)(\xi
_a),\mbox{ for }g\in G,$ where $L_g:G\to G$ is the left
translation by $g$. Moreover, suppose that
\[
\begin{array}{l}
K\big (\frac{\partial }{\partial t}_{|(x,e)}\big )=K^a_0(x)\xi
_a,\;\;\;\;\; K\big (\frac{\partial }{\partial x^i}_{|(x,e)}\big
)=K^a_i(x)\xi
_a,\,\;\;\;\;i\in \{ 1,\ldots ,m\},\\[10pt]
 B\big (\frac{\partial }{\partial
t}_{|(x,e)},\frac{\partial }{\partial x^i}_{|(x,e)}\big
)=B^a_{0i}(x)\xi _a,\,\;\;\;\;\; B\big (\frac{\partial }{\partial
x^i}_{|(x,e)},\frac{\partial }{\partial x^j}_{|(x,e)}\big
)=B^a_{ij}(x)\xi _a,\,\;\;\;i,j\in \{ 1,\ldots ,m\},
\end{array}
\]
for $x\in U$. Note that if $\{ c_{ab}^c \}$ are the structure
constants of $\mathfrak g$ with respect to the basis $\{ \xi _a
\}$ then
\[
B_{0i}^c=\frac{\partial K_i^c}{\partial t}- \frac{\partial
K_0^c}{\partial x^i}-K_0^aK_i^bc_{ab}^c,\;\;\;  B_{ij}^c=
\frac{\partial K^c_j}{\partial x^i}- \frac{\partial
K^c_i}{\partial x^j}-c^c_{ab}K^a_iK^b_j.
\]
Now, the horizontal lift of the vector fields $\{ \frac{\partial
}{\partial t}, \frac{\partial }{\partial x^i} \}$ to
$p^{-1}(U)\cong U\times G$ is given by
\[
\big ( \frac{\partial }{\partial t} \big ) ^h = \frac{\partial
}{\partial t}- K^a_0\xi ^L_a,\;\;\;\; \big (\frac{\partial
}{\partial x^i}\big )^h =\frac{\partial }{\partial x^i} - K^a_i\xi
^L _a.
\]
Therefore, the vector fields on $U\times G$
\begin{equation}\label{b'}
 \{ e_0=\frac{\partial }{\partial t}-
K^a_0\xi ^L_a , e_i=\frac{\partial }{\partial x^i} - K^a_i\xi ^L
_a ,e_b=\xi ^L_b\} \end{equation}
 are
$G$-invariant and they define a local basis $\{ e'_0,e'_i,e'_a\}$
of $\Gamma (\pi_Q|G)$. We will denote by $(
t,x^i;\dot{t},\dot{x}^i,\bar{v}^a)$ (res\-pec\-tively, $(
t,x^i;\dot{x}^i,\bar{v}^a)$) the corresponding fibred coordinates
on $TQ/G$ (res\-pec\-tively, $J^1\tau/G$). Then, if
$(\lcf\cdot,\cdot\rcf_{TQ/G},\rho_{TQ/G})$ is the Lie algebroid
structure on  $TQ/G$, we deduce that
\[
\begin{array}{c}
\lcf e'_0,e'_j\rcf_{TQ/G}  = - B^c _{0j} e'_c,\;\; \lcf
e'_i,e'_j\rcf_{TQ/G} = - B^c _{ij} e'_c,\;\;\lcf e'_0,e'_a
\rcf_{TQ/G} = c^c_{ab} K^b_0
e'_c,\\[8pt]
\;\;\lcf e'_i,e'_a \rcf_{TQ/G}  = c^c_{ab} K^b_i e'_c,\;\;\lcf
e'_a,e'_b
\rcf_{TQ/G}  = c^c_{ab} e'_c,\\[8pt] \rho_{TQ/G} (e'_0)=\displaystyle\frac{\partial }{\partial
t},\quad \rho_{TQ/G} (e'_i)=\displaystyle\frac{\partial }{\partial
x^i},\quad \rho_{TQ/G} (e'_a)=0.
\end{array}
\]
Thus, the (local) structure functions of the Lie algebroid
$\pi_{Q}|G:TQ/G\to M$ with respect to a local trivialization are
zero except the following
\begin{equation}\label{reduced-constant}
\begin{array}{c}
C^a_{0j}=-B^a_{0j},\quad C^c_{0a}=-C^c_{a0}= c^c_{ab}K^b_0,\\\quad
C^a_{ij}=-B^a_{ij},\quad C^c_{ia}=-C^c_{ai}= c^c_{ab}K^b_i, \quad
C^c_{ab}=c^c_{ab},\\\quad \rho ^0_0=1,\quad \rho _j^i = \delta
_{ij}.
\end{array}
\end{equation}

\subsubsection{Nonautonomous Hamilton-Poincar\'e equations} Let
$\pi^*:V^*\tau\to Q$ be the dual vector bundle to the vertical
bundle $\pi:V\tau\to Q$. Then, the Lie group $G$ acts on $V^*\tau$
and one may consider the corresponding quotient vector bundle
$\pi^*|G:V^*\tau/G\to M=Q/G.$ It is easy to prove that this vector
bundle is isomorphic to the dual vector bundle to
$\pi|G:V\tau/G\to M$. Moreover, using that the action of $G$ on
$V^*\tau$ is free, it follows that $V^*\tau$ is a principal
$G$-bundle over $V^*\tau/G$ with bundle projection
$p_{V^*\tau}:V^*\tau\to V^*\tau/G.$ Thus, we have the
corresponding Atiyah algebroid $\pi_{V^*\tau}|G:T(V^* \tau)/G\to
V^*\tau/G$ and, in addition, the exact $1$-form $(\pi_1^*)^*(dt)$
on $V^*\tau$, which is $G$-invariant, induces a $1$-cocycle
$\tilde{\eta}:T(V^*\tau)/G\to \R$ on
$\pi_{V^*\tau}|G:T(V^*\tau)/G\to V^*\tau/G.$  In fact,
$\tilde{\eta}$ is given by
\begin{equation}\label{AAe0}
\tilde{\eta}([X_{\alpha_q}])=X_{\alpha_q}(\pi^*_1),
\end{equation}
for $X_{\alpha_q}\in T_{\alpha_q}(V^*\tau)$ and $\alpha_q\in
V_q^*\tau.$

Now, denote by ${\mathcal L}^{\pi^*|G}(TQ/G)$ the prolongation of
the Atiyah algebroid $\pi_Q|G:TQ/G\to M$ over the fibration
$\pi^*|G:V^*\tau/G\to M$ and by $\eta:{\mathcal
L}^{\pi^*|G}(TQ/G)\to \R$ the $1$-cocycle defined by (\ref{eta}).
Then, we introduce the map
\[
(p_{TQ}\circ T\pi^*,Tp_{V^*\tau}):T(V^*\tau)\to {\mathcal
L}^{\pi^*|G}(TQ/G)\subseteq TQ/G\times T(V^*\tau/G)
\]
given by
\begin{equation}\label{AAe1}
(p_{TQ}\circ
T\pi^*,Tp_{V^*\tau})(X_{\alpha_q})=(p_{TQ}((T_{\alpha_q}\pi^*)(X_{\alpha_q})),
(T_{\alpha_q}p_{V^*\tau})(X_{\alpha_q})),
\end{equation}
for $X_{\alpha_q}\in T_{\alpha_q}(V^*\tau)$ and $\alpha_q\in
V_q^*\tau,$ where $p_{TQ}:TQ\to TQ/G$ is the canonical projection.

Next, we will consider the fibration $\nu\circ
\pi^*|G:V^*\tau/G\to \R$ and the Atiyah affgebroid associated with
the principal $G$-bundle $p_{V^*\tau}:V^*\tau\to (V^*\tau)/G$ and
the projection $\nu\circ \pi^*|G.$ Since $\nu\circ \pi^*|G\circ
p_{V^*\tau}=\pi_1^*$, it follows that the Atiyah affgebroid is the
quotient affine bundle $(\pi_1^*)_{1,0}|G:J^1\pi_1^*/G\to
V^*\tau/G$, where $(\pi_1^*)_{1,0}|G$ is defined by
\[
((\pi_1^*)_{1,0}|G)([j^1_t\gamma])=[(\pi_1^*)_{1,0}(j_t^1\gamma)]=[\gamma(t)],
\]
for $\gamma:I\subseteq \R\to V^*\tau$ a local section of $\pi^*_1$
and $t\in I.$ Note that the Lie affgebroids
$\tilde{\eta}^{-1}\{1\}$ and  $(\pi_1^*)_{1,0}|G:J^1\pi_1^*/G\to
V^*\tau/G$ may be identified in such a way that the bidual Lie
algebroid to $(\pi^*_1)_{1,0}|G:J^1\pi_1^*/G\to V^*\tau/G$ also
may be identified with $\pi_{V^*\tau}|G:T(V^*\tau)/G\to
V^*\tau/G.$ Under the above identifications, the $1$-cocycle
$1_{(J^1\pi_1^*/G)}$ is just $\tilde{\eta}$.

\begin{theorem}\label{t9.1}
$(i)$ The map \kern-1.5pt$(p_{TQ}\kern-1.5pt\circ\kern-1.5pt
T\pi^*\kern-1.5pt,Tp_{V^*\tau})$ induces an isomorphism
\kern-1.5pt(\ltilde{70}{$p_{TQ}\circ T\pi^*,Tp_{V^*\tau})$,} over
the identity of $V^*\tau/G$, between the Lie algebroids
$\pi_{V^*\tau}|G:T(V^*\tau)/G\to V^*\tau/G$ and
$\tau^{\pi^*|G}_{TQ/G}:{\mathcal L}^{\pi^*|G}(TQ/G)\to V^*\tau/G.$
Furthermore,
\begin{equation}\label{AAe2}
((\mbox{\ltilde{70}{$p_{TQ}\circ
T\pi^*,Tp_{V^*\tau}$}}),Id)^*\eta=\tilde{\eta}.
\end{equation}

$(ii)$ The restriction of the map (\ltilde{70}{$p_{TQ}\circ
T\pi^*,Tp_{V^*\tau}$}) to $J^1\pi^*_1/G$ induces an isomorphism,
over the identity of $V^*\tau/G$, between the Atiyah affgebroid
$(\pi_1)^*_{1,0}|G:J^1\pi^*_1/G\to V^*\tau/G$ and the Lie
affgebroid
$\widetilde{\pi_{V^*\tau|G}}:\rho_{J^1\tau/G}^*(T(V^*\tau/G))\to
V^*\tau/G.$
\end{theorem}
\begin{proof}
The cotangent lift of $\Phi$ defines an action of $G$ on $T^*Q$
and it is clear that the ca\-no\-ni\-cal projection $\mu:T^*Q\to
V^*\tau$ is $G$-equivariant. Thus, $\mu$ induces an epimorphism
$\mu|G:T^*Q/G\to V^*\tau/G$ between the quotient vector bundles
$T^*Q/G$ and $V^*\tau/G$. Therefore, if $(\pi_Q|G)^*:T^*Q/G\to M$
is the dual vector bundle to the Atiyah algebroid
$\pi_{Q}|G:TQ/G\to M$ and ${\mathcal L}^{(\pi_Q|G)^*}(TQ/G)$ is
the prolongation of $\pi_Q|G:TQ/G\to M$ over the fibration
$(\pi_Q|G)^*:T^*Q/G\to M$, one may introduce the epimorphism
$$(Id,T(\mu|G)):{\mathcal L}^{(\pi_Q|G)^*}(TQ/G)\to{\mathcal
L}^{\pi^*|G}(TQ/G)$$ between the Lie algebroids ${\mathcal
L}^{(\pi_Q|G)^*}(TQ/G)\subseteq TQ/G\times T(T^*Q/G)$ and
${\mathcal L}^{(\pi^*|G)}(TQ/G)\subseteq TQ/G\times T(V^*\tau/G)$
defined by
\[
(Id,
T(\mu|G))([u_q],X_{[\alpha_{q'}]})=([u_q],(T_{[\alpha_{q'}]}(\mu|G))(X_{[\alpha_{q'}]})),
\]
for $u_q\in T_qQ$ and $X_{[\alpha_{q'}]}\in
T_{[\alpha_{q'}]}(T^*Q/G),$ with $\alpha_{q'}\in T^*_{q'}Q.$

The tangent map to $\mu$, $T\mu:T(T^*Q)\to T(V^*\tau),$ is also
$G$-invariant with respect to the tangent actions of $G$ on
$T(T^*Q)$ and $T(V^*\tau).$ This implies that $T\mu$ induces an
epimorphism between the vector bundles $\pi_{T^*Q}|G:T(T^*Q)/G\to
T^*Q/G$ and $\pi_{V^*\tau}|G:T(V^*\tau)/G\to V^*\tau/G.$ In
addition, since $T\mu$ is an epimorphism over $\mu$ between the
Lie algebroids $\pi_{T^*Q}:T(T^*Q)\to T^*Q$ and
$\pi_{V^*\tau}:T(V^*\tau)\to V^*\tau,$ we deduce that the map
$T\mu|G:T(T^*Q)/G\to T(V^*\tau)/G$ is also an epimorphism over
$\mu|G:T^*Q/G\to V^*\tau/G$ between the Lie algebroids $T(T^*Q)/G$
and $T(V^*\tau)/G.$

 Now, denote by $\pi_Q^*:T^*Q\to Q$ the bundle
projection, by $p_{T^*Q}:T^*Q\to T^*Q/G$ the canonical projection
and by $(p_{TQ}\circ T\pi_Q^*,Tp_{T^*Q}):T(T^*Q)\to {\mathcal
L}^{(\pi_Q|G)^*}(TQ/G)$ the map defined by
\[
(p_{TQ}\circ
T\pi_{Q}^*,Tp_{T^*Q})(X_{\alpha_q})=(p_{TQ}((T_{\alpha_q}\pi^*_Q)(X_{\alpha_q})),(T_{\alpha_q}p_{T^*Q}))(X_{\alpha_q})),
\]
for $X_{\alpha_q}\in T_{\alpha_q}(T^*Q),$ with $\alpha_q\in
T_q^*Q.$

Then, this map induces an isomorphism \ltilde{70}{$(p_{TQ}\circ
T\pi_Q^*,Tp_{T^*Q}),$} over the identity of $T^*Q/G,$ between the
Lie algebroids $T(T^*Q)/G\to T^*Q/G$ and ${\mathcal
L}^{(\pi_Q|G)^*}(TQ/G)$ $\to T^*Q/G$ (see Theorem 9.3 in
\cite{LMM}).

On the other hand, using (\ref{AAe1}), it follows that the map
$(p_{TQ}\circ T\pi^*,Tp_{V^*\tau}):T(V^*\tau)\to {\mathcal
L}^{\pi^*|G}(TQ/G)$ induces a map \ltilde{70}{$(p_{TQ}\circ
T\pi^*,Tp_{V^*\tau})$}$:T(V^*\tau)/G\to {\mathcal
L}^{\pi^*|G}(TQ/G)$ between the spaces $T(V^*\tau)/G$ and
${\mathcal L}^{\pi^*|G}(TQ/G).$  Furthermore, it is easy to prove
that the following diagram is commutative

\begin{picture}(375,60)(40,40)
\put(135,20){\makebox(0,0){${\mathcal L}^{(\pi_Q|G)^*}(TQ/G)$}}
\put(215,25){$(Id,T(\mu|G))$} \put(175,20){\vector(1,0){120}}
\put(330,20){\makebox(0,0){${\mathcal L}^{\pi^*|G}(TQ/G)$}}
\put(70,50){\ltilde{70}{$(p_{TQ}\circ T\pi_Q^*,Tp_{T^*Q})$}}
\put(160,70){\vector(0,-1){40}}
\put(320,50){\ltilde{70}{$(p_{TQ}\circ T\pi^*,Tp_{V^*\tau})$}}
\put(315,70){\vector(0,-1){40}}
\put(150,80){\makebox(0,0){$T(T^*Q)/G$}} \put(230,85){$T\mu|G$}
\put(175,80){\vector(1,0){120}}
\put(325,80){\makebox(0,0){$T(V^*\tau)/G$}}
\end{picture}

\vspace{1.2cm}

Thus, the map \ltilde{70}{$(p_{TQ}\circ T\pi^*,Tp_{V^*\tau})$} is
an epimorphism over the identity of $V^*\tau/G$ between the Lie
algebroids $T(V^*\tau)/G\to V^*\tau/G$ and ${\mathcal L
}^{\pi^*|G}(TQ/G)\to V^*\tau/G$. Therefore, using that the ranks
of these Lie algebroids are equal, we deduce that the map
\ltilde{70}{$(p_{TQ}\circ T\pi^*,Tp_{V^*\tau})$} is an
isomorphism, over the identity of $V^*\tau/G$, between the Lie
algebroids. Moreover, from (\ref{eta}), (\ref{AAe0}) and
(\ref{AAe1}), it follows (\ref{AAe2}). This proves $(i).$

On the other hand, using $(i)$ and since  $\eta^{-1}\{1\}$ is the
total space of the Lie affgebroid
$\widetilde{\pi_{V^*\tau|G}}:\rho_{J^1\tau|G}(T(V^*\tau/G))$ $\to
V^*\tau/G$ (see Section \ref{sec5}), we deduce $(ii).$
\end{proof}
As we know, the bidual Lie algebroid to $\tau_{1,0}|G$ may be
identified with the Atiyah algebroid $\pi_Q|G:TQ/G\to M.$ Thus,
the vector bundles $(J^1\tau/G)^+\to M$ and $T^*Q/G\to M$ are
isomorphic.

Now, suppose that $h:V^*\tau/G\to (J^1\tau/G)^+\cong T^*Q/G$ is a
Hamiltonian section. Then, we have that $h$ induces a Hamiltonian
section $\bar{h}:V^*\tau\to T^*Q$ with respect to the Lie
affgebroid $\tau_{1,0}:J^1\tau\to Q$ (note that
$(T^*Q/G)_{p(q)}\cong T_q^*Q,$ for all $q\in Q$). Moreover, using
the $G$-equivariant character of $\bar{h}$, we deduce that the
standard cosymplectic structure $(\Omega_{\bar{h}},\bar\eta)$ on
$V^*\tau$ (see Section \ref{sec9.1}) is also $G$-invariant. Thus,
it induces a cosymplectic structure on the Lie algebroid
${\mathcal L}^{\pi^*|G}(TQ/G)\cong T(V^*\tau)/G\to V^*\tau/G.$
This cosymplectic structure is just $(\Omega_h,\eta)$. In
addition, it is clear that the Reeb vector field $R_{\bar{h}}$ of
$(\Omega_{\bar{h}},\bar{\eta})$ is also $G$-invariant and,
therefore, it induces a section of the Lie algebroid ${\mathcal
L}^{\pi^*|G}(TQ/G)\cong T(V^*\tau)/G\to V^*\tau/G$ which is just
the Reeb section $R_h$ of $(\Omega_h,\eta)$. Consequently, the
solutions of the Hamilton equations for $h$ (that is, the integral
curves on $V^*\tau/G$ of the vector field
$\rho_{TQ/G}^{\pi^*|G}(R_h))$ are just the solutions of the {\it
nonautonomous Hamilton-Poincar\'{e} equations } for $\bar{h}.$

Next, we will obtain the expression of these equations. Let
$\{e_0',e_i',e_a'\}$ be the local basis of $\Gamma(\pi_Q|G)$
considered in Section \ref{sec9.2.1}. Then, $\{e_i',e_a'\}$ is a
local basis of $\Gamma(\pi|G)$ and may consider the corresponding
local coordinates $(t,x^i,p_i,\bar{p}_a)$ (respectively,
$(t,x^i;p_t,p_i,\bar{p}_a)$) on $V^*\tau/G$ (respectively,
$T^*Q/G$). Suppose that $h$ is locally given by
$$h(t,x^i,p_i,\bar{p}_a)=(t,x^i,-H(t,x^j,p_j,\bar{p}_b),
p_i,\bar{p}_a).$$ Using (\ref{ech}) and (\ref{reduced-constant}),
we get the Hamilton equations for $h$ (the {\it
Hamilton-Poincar\'e equations} for $\bar{h}$),
\[
\begin{array}{l}
\displaystyle \frac{dx^i}{dt}=\frac{\partial H}{\partial
p_i},\;\;\; \displaystyle  \frac{dp_i}{dt}=-\frac{\partial
H}{\partial x^i}-\bar{p}_b\Big ( B_{0i}^b+B_{ki}^b\frac{\partial
H}{\partial p_k}+c^b_{ac}K_i^c
\frac{\partial H}{\partial \bar{p}_a}\Big ),\\[10pt]
\displaystyle  \frac{d\bar{p}_a}{dt}=\bar{p}_c \Big
(c^c_{ab}K_0^b+c^c_{ab}K_k^b\frac{\partial H }{\partial p_k}-
c^c_{ab}\frac{\partial H }{\partial \bar{p}_b} \Big ).
\end{array}
\]
Now, we will obtain the local equations defining the Lagrangian
submanifold $S_h=R_h(V^*\tau/G)$ of the symplectic Lie affgebroid
$\rho ^\ast _{J^1\tau/G}(T(V^*\tau/G) )\cong J^1\pi_1^*/G$. For
this purpose, we consider local coordinates $(t,x^i,p_i,\bar{p}_a;
\dot{x}^i,\bar{v}^a,\dot{p}_i,\dot{\bar{p}}_a)$ on
$\rho^*_{J^1\tau/G}(T(V^*\tau/G))\cong J^1\pi_1^*/G$.

 Using
(\ref{rh}) and (\ref{reduced-constant}), we get that the Reeb
section $R_h$ is given by
\[
\begin{array}{rcl}
R_h&=&\displaystyle \tilde{e}_0+\frac{\partial H}{\partial
p_i}\tilde{e}_{i}+\frac{\partial H}{\partial
\bar{p}_a}\tilde{e}_{a}-\Big ( \frac{\partial H}{\partial
x^i}+B_{0i}^b\bar{p}_b+ B_{ki}^b\bar{p}_b\frac{\partial
H}{\partial p_k}+c_{ac}^bK^c_i\bar{p}_b\frac{\partial H}{\partial
\bar{p}_a} \Big )\bar{e}_{i}\\[5pt]
&&\displaystyle +\Big(  c_{ab}^c K^b_0\bar{p}_c + c_{ab}^c
K^b_k\bar{p}_c\frac{\partial H}{\partial p_k}-
c_{ab}^c\bar{p}_c\frac{\partial H}{\partial \bar{p}_b}\Big
)\bar{e}_a,
\end{array}
\]
where
$\{\tilde{e}_0,\tilde{e}_i,\tilde{e}_a,\bar{e}_i,\bar{e}_a\}$ is
the local basis of sections of ${\mathcal L}^{\pi^*|G}(TQ/G)\to
V^*\tau/G$ induced by $\{e_0',e_i',e_a'\}$. Thus, the local
equations defining the submanifold $S_h$ of
$\rho^*_{J^1\tau/G}(T(V^*\tau/G))\cong J^1\pi_1^*/G$

\[
\begin{array}{l}
\displaystyle \bar{v}^a=\frac{\partial H}{\partial
\bar{p}_a},\;\;\;\; \displaystyle \dot{x}^i=\frac{\partial
H}{\partial p_i},\;\;\;\; \displaystyle  \dot{p}_i=-\frac{\partial
H}{\partial x^i}-\bar{p}_b\Big ( B_{0i}^b+B_{ki}^b\frac{\partial
H}{\partial p_k}+c^b_{ac}K_i^c
\frac{\partial H}{\partial \bar{p}_a}\Big ),\\[10pt]
\displaystyle \dot{\bar{p}}_a=\bar{p}_c \Big
(c^c_{ab}K_0^b+c^c_{ab}K_k^b\frac{\partial H }{\partial p_k}-
c^c_{ab}\frac{\partial H }{\partial \bar{p}_b} \Big ).
\end{array}
\]
In other words,
$$\displaystyle \bar{v}^a=\frac{\partial H}{\partial \bar{p}_a},$$
\begin{equation}\label{HP-eqs}
\begin{array}{l}
\displaystyle \frac{dx^i}{dt}=\frac{\partial H}{\partial
p_i},\;\;\;\;\; \displaystyle  \frac{dp_i}{dt}=-\frac{\partial
H}{\partial x^i}-\bar{p}_b\Big ( B_{0i}^b+B_{ki}^b\frac{\partial
H}{\partial p_k}+c^b_{ac}K_i^c
\frac{\partial H}{\partial \bar{p}_a}\Big ),\\[10pt]
\displaystyle  \frac{d\bar{p}_a}{dt}=\bar{p}_c \Big
(c^c_{ab}K_0^b+c^c_{ab}K_k^b\frac{\partial H }{\partial p_k}-
c^c_{ab}\frac{\partial H }{\partial \bar{p}_b} \Big ).
\end{array}
\end{equation}
Note that Eqs. (\ref{HP-eqs}) are just the Hamilton-Poincar\'e
equations for $\bar{h}$.

\subsubsection{Nonautonomous Lagrange-Poincar\'e equations}

The action $J^1\Phi$ of $G$ on $J^1\tau$ is free and $J^1\tau$ is
the total space of a principal $G$-bundle over $J^1\tau/G$ with
bundle projection $p_{J^1\tau}:J^1\tau\to J^1\tau/G$. Therefore,
we have the corresponding Atiyah algebroid $\pi_{J^1\tau}|G:
T(J^1\tau)/G\to J^1\tau/G$ and, in addition, the exact $1$-form
$\tau_1^*(dt)$ on $J^1\tau$, which is $G$-invariant, induces a
$1$-cocycle $\tilde{\phi_0}:T(J^1\tau)/G\to \R$ on
$\pi_{J^1\tau}|G:T(J^1\tau)/G\to J^1\tau/G$. Here,
$\tau_1:J^1\tau\to \R$ is the map $\tau\circ \tau_{1,0}.$

Now, denote by ${\mathcal L}^{\tau_{1,0}|G}(TQ/G)$ the
prolongation of the Atiyah algebroid $\pi_Q|G:TQ/G\to M$ over the
fibration $\tau_{1,0}|G:J^1\tau/G\to M$ and by $\phi_0:{\mathcal
L}^{\tau_{1,0}|G}(TQ/G)\to \R$ the $1$-cocycle defined by
(\ref{T0}). We recall that $\phi_0^{-1}\{1\}$ is the Lie
affgebroid $(\tau_{1,0}|G)^{(\tau_{1,0}|G)}:{\mathcal
J}^{J^1\tau/G}(J^1\tau/G)\to J^1\tau/G$ (see Section \ref{sec3}).
Moreover, we may introduce the map
\[
(p_{TQ}\circ T\rho_{J^1\tau},Tp_{J^1\tau}):T(J^1\tau)\to {\mathcal
L}^{\tau_{1,0}|G}(TQ/G)\subseteq TQ/G\times T(J^1\tau/G)
\]
given by
\begin{equation}\label{AAe3}
(p_{TQ}\circ
T\rho_{J^1\tau},Tp_{J^1\tau})(X_{j_t^1\gamma})=(p_{TQ}((T_{j_t^1\gamma}\rho_{J^1\tau})(X_{j_t^1\gamma})),(T_{j^1_t\gamma}p_{J^1\tau}
(X_{j_t^1\gamma}))),
\end{equation}
for $X_{j_t^1\gamma}\in T_{j_t^1\gamma}(J^1\tau)$ and
$j_t^1\gamma\in J^1\tau$, where $\rho_{J^1\tau}:J^1\tau\to TQ$
 is the anchor map of the Lie affgebroid
$\tau_{1,0}:J^1\tau\to Q.$

On the other hand, we consider the fibration $\nu\circ
\tau_{1,0}|G:J^1\tau/G\to \R$ and the Atiyah affgebroid associated
with the principal $G$-bundle $p_{J^1\tau}:J^1\tau\to J^1\tau/G$
and the projection $\nu\circ \tau_{1,0}|G$. Since $\nu\circ
\tau_{1,0}|G\circ p_{J^1\tau}=\tau_1$, it follows that the Atiyah
affgebroid is the quotient affine bundle
$(\tau_1)_{1,0}|G:J^1\tau_1/G\to J^1\tau/G$, where
$(\tau_1)_{1,0}|G$ is defined by
\[
((\tau_1)_{1,0}|G)([j^1_t\gamma])=[(\tau_1)_{1,0}(j^1_t\gamma)]=[\gamma(t)],
\]
for $\gamma:I\subseteq \R\to J^1\tau$ a local section of $\tau_1$
and $t\in I.$ Note that the Lie affgebroids
$\tilde{\phi}_0^{-1}\{1\}$ and $(\tau_1)_{1,0}|G:(J^1\tau_1)/G\to
M$ may be identified in such a way that the bidual Lie algebroid
to $(\tau_1)_{1,0}|G:(J^1\tau_1)/G\to J^1\tau/G$ may be identified
with $\pi_{J^1\tau}|G:T(J^1\tau)/G\to J^1\tau/G$.  Under the above
identifications, the $1$-cocycle $1_{(J^1\tau_1)/G}$ is just
$\tilde{\phi}_0.$

\begin{theorem}\label{t9.2}
$(i)$ The map $(p_{TQ}\circ
T\rho_{J^1\tau},Tp_{J^1\tau}):T(J^1\tau)\to {\mathcal
L}^{\tau_{1,0}|G}(TQ/G)$ induces an isomorphism
\ltilde{70}{$(p_{TQ}\circ T\rho_{J^1\tau},Tp_{J^1\tau})$}, over
the identity of $J^1\tau/G$, between the Lie algebroids
$\pi_{J^1\tau}|G:T(J^1\tau)/G\to J^1\tau/G$ and
$(\pi_Q|G)^{\tau_{1,0}|G}:{\mathcal L}^{\tau_{1,0}|G}(TQ/G)\to
J^1\tau/G.$ Furthermore,
\[
\mbox{\ltilde{70}{$((p_{TQ}\circ
T\rho_{J^1\tau},Tp_{J^1\tau})$}},Id)^*\phi_0=\tilde{\phi}_0.
\]

$(ii)$ The restriction of the map \ltilde{70}{$(p_{TQ}\circ
T\rho_{J^1\tau},Tp_{J^1\tau})$} to $(J_1\tau_1)/G$ induces an
isomorphism, over the identity of $(J^1\tau)/G,$ between the
Atiyah affgebroid $(\tau_1)_{1,0}|G:(J^1\tau_1)/G\to J^1\tau/G$
and the Lie affgebroid $(\tau_{1,0}|G)^{(\tau_{1,0}|G)}:{\mathcal
J}^{(J^1\tau/G)}(J^1\tau/G)\to J^1\tau/G.$

\end{theorem}
\begin{proof} $(i)$ From (\ref{AAe3}) it follows that $(p_{TQ}\circ
T\rho_{J^1\tau},Tp_{J^1\tau})$ induces a morphism, over the
identity of $J^1\tau/G,$ between the Lie algebroids
$\pi_{J^1\tau}|G:T(J^1\tau)/G\to J^1\tau/G$ and
$(\pi_Q|G)^{\tau_{1,0}|G}:{\mathcal L}^{\tau_{1,0}|G}(TQ/G)\to
J^1\tau/G.$ We denote this morphism by \ltilde{70}{$(p_{TQ}\circ
T\rho_{J^1\tau}, Tp_{j^1\tau})$}. Moreover, proceeding as in the
proof of Theorem \ref{t9.1} and using the $G$-equivariant
character of $\rho_{J^1\tau}$, Theorem 9.1 in \cite{LMM} and  the
fact that the tangent map to $\rho_{J^1\tau}$,
$T(\rho_{J^1\tau}):T(J^1\tau)\to T(TQ)$, is a $G$-equivariant Lie
algebroid monomorphism over $\rho_{J^1\tau}, $we deduce $(i).$
\medskip

$(ii)$ It follows using $(i).$
\end{proof}

Now, suppose that $l:J^1\tau/G\to \R$ is a Lagrangian function.
Then, we will denote by $L:J^1\tau\to \R$ the function given by
$L=l\circ p_{J^1\tau}.$ Since $L$ is a $G$-invariant Lagrangian
function, we deduce that the Poincar\'{e}-Cartan $2$-form
$\Omega_L$ on $J^1\tau$ is also $G$-invariant. Thus, it induces a
section of the vector bundle $\wedge^2({\mathcal
L}^{\tau_{1,0}|G}(TQ/G)^*)\cong \wedge^2(T^*(J^1\tau)/G)\to
J^1\tau/G.$ This section is just the Poincar\'{e}-Cartan
$2$-section $\Omega_l$ associated with $l$. Therefore, the
solutions of the Euler-Lagrange equations for $l$ are just the
solutions of the {\it nonautonomous Lagrange-Poincar\'{e}
equations } for $L$.

Next, we will obtain the expression of these equations. Let
$\{e_0,e_i,e_a\}$ be the local basis of $G$-invariant vector
fields on $Q$ defined as in (\ref{b'}) and
$(t,x^i,\dot{x}^i,\bar{v}^a)$ (respectively,
$(t,x^i,\dot{t},\dot{x}^i,\bar{v}^a)$) be the corresponding local
fibred coordinates on $J^1\tau/G$ (respectively,
$\widetilde{J^1\tau/G}=TQ/G$). Using (\ref{EL}) and
(\ref{reduced-constant}),
we obtain the Euler-Lagrange equations for $l$ (that is, the {\it
Lagrange-Poincar\'{e} equations} for $L$),
\[
\begin{array}{l}
\displaystyle  \frac{\partial l}{\partial x^i} -\frac{d}{dt}\Big (
\frac{\partial l}{\partial \dot{x}^i}\Big )=
 B_{0i}^b\frac{\partial l}{\partial \bar{v}^b}+B_{ji}^b\dot{x}^j\frac{\partial
l}{\partial \bar{v}^b}+c^b_{dc}K_i^c \bar{v}^d \frac{\partial
l}{\partial \bar{v}^b},\\[10pt]
\displaystyle  \dot{\bar{p}}_a= \bar{p}_b \Big (
c_{ac}^bK_0^c+c_{ac}^bK_j^c\dot{x}^j -c_{ad}^b\bar{v}^d \Big ).
\end{array}
\]



On the other hand, we consider the local coordinates
$(t,x^i,p_i,\bar{p}_a;\dot{x}^i,\bar{v}^a,\dot{{p}}_i,\dot{\bar{p}}_a)$
on $J^1\pi_1^*/G\cong$ $\rho ^\ast _{J^1\tau/G}(T (V^*\tau/G))$.
Then, the local equations defining the Lagrangian submanifold
$S_l$ are (see (\ref{e8.80^0}))
\[
\begin{array}{l}
\displaystyle p_i=\frac{\partial l}{\partial \dot{x}^i},\quad
\bar{p}_a=\frac{\partial l}{\partial \bar{v}^a},\;\;\;\;\;
\displaystyle \dot{p}_i=\frac{\partial l}{\partial x^i}+
 B_{i0}^b\frac{\partial l}{\partial \bar{v}^b}+B_{ij}^b\dot{x}^j\frac{\partial
l}{\partial \bar{v}^b}
-c^b_{dc}K_i^c \bar{v}^d \frac{\partial l}{\partial \bar{v}^b},\\[10pt]
\displaystyle  \dot{\bar{p}}_a=\frac{\partial l}{\partial
\bar{v}^b}\Big (c_{ac}^bK_0^c+c_{ac}^bK_j^c\dot{x}^j
-c_{ad}^b\bar{v}^d \Big ).
\end{array}
\]
or, in other words,
\begin{equation}\label{momenta-eq}
\begin{array}{l}
\displaystyle p_i=\frac{\partial l}{\partial \dot{x}^i},\quad
\bar{p}_a=\frac{\partial l}{\partial \bar{v}^a},
\end{array}
\end{equation}
\begin{equation}\label{LP-eqs}
\begin{array}{l}
 \displaystyle
\frac{\partial l}{\partial x^i} -\frac{dp_i}{dt}= \bar{p}_b\Big (
B_{0i}^b+B_{ji}^b\dot{x}^j+c^b_{dc}K_i^c
\bar{v}^d \Big ),\\[10pt]
\displaystyle  \dot{\bar{p}}_a= \bar{p}_b \Big (
c_{ac}^bK_0^c+c_{ac}^bK_j^c\dot{x}^j -c_{ad}^b\bar{v}^d \Big ).
\end{array}
\end{equation}
Eqs. (\ref{momenta-eq}) give the definition of the momenta and
Eqs. (\ref{LP-eqs}) are just the Lagrange-Poincar\'e equations for
L.

\noindent {\small\sc David Iglesias:} {\small \it Department of
Mathematics, Penn State University, University Park, State
College, PA 16802,  USA.} {\small\it E-mail:
iglesias@math.psu.edu}

\medskip

\noindent{\small \sc Juan Carlos Marrero,  Edith Padr\'{o}n and Diana
Sosa:} {\small\it Departamento de Matem\'{a}tica Fundamental,
Universidad de La Laguna, La Laguna, Tenerife, Canary islands,
SPAIN. } {\small\it E-mail: jcmarrer@ull.es, mepadron@ull.es,
dnsosa@ull.es}

\end{document}